\documentclass[11pt,a4paper]{article}
\usepackage[T1]{fontenc}
\usepackage{lmodern}
\usepackage[margin=27mm,headheight=14pt]{geometry}
\usepackage{amsmath,amssymb,amsthm,mathtools}
\usepackage{microtype,booktabs,longtable,array,enumitem}
\usepackage{fancyhdr}
\usepackage[hyphens]{url}
\usepackage[hidelinks]{hyperref}
\hypersetup{pdftitle={Tightness of bulk crossing heights for the planar discrete Gaussian free field},pdfauthor={Xiang Fang, Yi Xu},pdfsubject={Uniform upper and lower tails for nearest-neighbor two-sided level-set crossings}}
\setlist{itemsep=2pt,topsep=4pt,leftmargin=*}
\newtheorem{theorem}{Theorem}[section]
\newtheorem{lemma}[theorem]{Lemma}
\newtheorem{corollary}[theorem]{Corollary}

\theoremstyle{definition}

\theoremstyle{remark}

\newcommand{\PP}{\mathbb P}
\newcommand{\EE}{\mathbb E}
\newcommand{\NN}{\mathrm{NN}}
\newcommand{\LR}{\mathrm{LR}}

\newcommand{\Poi}{\operatorname{PPP}}
\newcommand{\GammaDist}{\operatorname{Gamma}}
\newcommand{\Var}{\operatorname{Var}}
\newcommand{\Cov}{\operatorname{Cov}}
\newcommand{\ind}{\mathbf 1}
\newcolumntype{L}[1]{>{\raggedright\arraybackslash}p{#1}}
\title{Tightness of bulk crossing heights for the planar discrete Gaussian free field}
\author{Xiang Fang\thanks{Department of Applied Mathematics, National Yang Ming Chiao Tung University, Hsinchu City 30010, Taiwan. Email: \texttt{xfang@nycu.edu.tw}.}\and Yi Xu\thanks{School of Mathematical Sciences, Dalian University of Technology, Dalian 116024, P. R. China. Email: \texttt{xyi@mail.dlut.edu.cn}.}}
\date{}
\begin{document}
\maketitle

\begin{abstract}
Let $\phi_N$ be the planar discrete Gaussian free field in $[-2N,2N]^2\cap\mathbb Z^2$, with zero boundary condition imposed by killing the associated simple random walk on exit. We consider the least height $H_N$ at which the two-sided level set $\{|\phi_N|\leq a\}$ contains a nearest-neighbor left--right crossing of $[-N,N]^2\cap\mathbb Z^2$. We prove that $(H_N)_{N\geq1}$ is tight, uniformly over all box sizes. More precisely, for every $a\geq0$ and $N\geq1$,
\[
 \Phi\!\left(-a\sqrt{\frac{2}{\log2}}\right)
 \leq\PP(H_N>a)
 \leq4\exp\!\left\{-2^{-29}\frac{a^2}{\log(2+a)}\right\},
\]
where $\Phi$ is the standard normal distribution function. The upper bound yields uniform polynomial moments and a uniform exponential moment of $H_N^2/\log(2+H_N)$. For the upper bound, we compare a marked loop-soup exploration with a graph in which independent site variables open incoming edges. We then estimate connections in a massive field using two auxiliary vertices and bound the Gaussian remainder when the mass is removed. The lower bound uses a Gaussian translation inequality and a square-shell estimate for the Dirichlet energy, without using the massive-field argument. In particular, a height independent of $N$ gives any prescribed crossing probability below one, whereas every fixed finite height has a uniformly positive failure probability.
\end{abstract}


\tableofcontents

\section{Introduction and results}\label{sec:introduction}

At a bulk vertex of the two-dimensional discrete Gaussian free
field, the standard deviation grows on the scale $\sqrt{\log N}$
as the domain size $N$ increases. Thus, for any fixed $a$, the
probability that the field at a given bulk vertex lies in
$[-a,a]$ tends to zero as $N$ grows. A crossing path, however,
can choose its route after the field is sampled. Can such a path
cross a large box while staying entirely in the two-sided level
set $\{|\phi_N|\leq a\}$, for a height $a$ that does not depend
on $N$?

We show that the answer is yes when the crossing takes place
in an inner box separated from the boundary by a distance
proportional to $N$. More precisely, for every $\varepsilon>0$,
there is a finite height $a_\varepsilon$, independent of $N$,
such that a nearest-neighbor left--right crossing exists with
probability at least $1-\varepsilon$ for every $N$. On the other
hand, at any fixed finite height, the probability that no such
crossing exists is bounded away from zero, uniformly in $N$.
Theorem~\ref{thm:main} gives explicit upper and lower bounds
for these probabilities.

\subsection{The model and the crossing height}

For an integer $r\geq0$, let
\[
 B_r=[-r,r]^2\cap\mathbb Z^2.
\]
Fix an integer $N\geq1$ and set $D=B_{2N}$. On $D$, write
\[
 P_D(x,y)=\frac14\ind_{\{|x-y|_1=1\}},\qquad
 \mathsf Q_D=I-P_D,\qquad
 G_D=\mathsf Q_D^{-1}=\sum_{n\geq0}P_D^n.
\]
Thus $G_D$ is the Green matrix of nearest-neighbor simple random walk killed on leaving $D$, with its visit at time zero included. All matrices act with respect to counting measure. The field is
\begin{equation}\label{eq:field}
 \phi_N\sim\mathcal N(0,G_D).
\end{equation}
Equivalently, every lattice edge has conductance $1/4$, and a vertex $x\in D$ has killing $1-\sum_{y\in D}P_D(x,y)$. These conventions fix all height constants below.

Let $\LR_{\NN}(B_N)$ be the finite set of simple nearest-neighbor vertex paths in $B_N$ joining its left and right sides. The endpoints are free on these sides. For a function $f:D\to\mathbb R$, define
\begin{equation}\label{eq:H}
 H_N(f)=\min_{\gamma\in\LR_{\NN}(B_N)}\max_{v\in\gamma}|f(v)|,
 \qquad H_N=H_N(\phi_N).
\end{equation}
Deleting loops from a path does not increase its maximum, so allowing nonsimple crossing paths gives the same quantity. For $a\geq0$, the event $\{H_N\leq a\}$ is exactly the existence of a nearest-neighbor left--right crossing in the closed band $\{|\phi_N|\leq a\}\cap B_N$. The original crossings always stay in $B_N$. Auxiliary loop paths and connections used in the proof may travel throughout $D$.

All logarithms are natural. We use
\[
 \Phi(t)=\frac1{\sqrt{2\pi}}\int_{-\infty}^t e^{-s^2/2}\,ds,
 \qquad F(a)=\frac{a^2}{\log(2+a)}\quad(a\geq0),
 \qquad c_*=2^{-29}.
\]

\subsection{Uniform tails and their consequences}

\begin{theorem}[Uniform bounds for the crossing height]\label{thm:main}
For every integer $N\geq1$ and every real $a\geq0$,
\begin{equation}\label{eq:main}
 \Phi\!\left(-a\sqrt{\frac{2}{\log2}}\right)
 \leq\PP(H_N>a)
 \leq4e^{-c_*F(a)}.
\end{equation}
\end{theorem}

The constants are explicit sufficient choices; we make no claim of optimality. The right side of~\eqref{eq:main} may be replaced by its minimum with one. The upper bound is proved in Section~\ref{sec:gaussian}, after the two massive-field estimates developed in Sections~\ref{sec:comparison}--\ref{sec:probes}. The lower bound is proved in Section~\ref{sec:lower} by a separate argument.

\begin{corollary}[Uniform tightness and fixed-height failure]\label{cor:tightness}
For every $\varepsilon>0$ there is a finite $a_\varepsilon$ such that
\begin{equation}\label{eq:tightness}
 \inf_{N\geq1}\PP(H_N\leq a_\varepsilon)\geq1-\varepsilon.
\end{equation}
In particular,
\[
 \lim_{a\to\infty}\sup_{N\geq1}\PP(H_N>a)=0.
\]
Conversely, every fixed finite $a\geq0$ satisfies
\begin{equation}\label{eq:fixed-height-failure}
 \sup_{N\geq1}\PP(H_N\leq a)
 \leq1-\Phi\!\left(-a\sqrt{\frac{2}{\log2}}\right)<1.
\end{equation}
\end{corollary}
\begin{proof}
The function $F$ tends to infinity. For $0<\varepsilon<1$, choose $a_\varepsilon$ so that $4e^{-c_*F(a_\varepsilon)}\leq\varepsilon$ and apply the upper bound in Theorem~\ref{thm:main}. The cases $\varepsilon\geq1$ are immediate. Rearranging its lower bound gives~\eqref{eq:fixed-height-failure}.
\end{proof}

The order of the quantifiers matters. A desired confidence determines a height that works for every box size. No finite height makes the crossing probability tend to one. The bounds also rule out convergence of $H_N$ in probability to a finite deterministic value. They do not establish convergence in distribution or identify a limiting law.

\paragraph{Moment consequences.}
For every $0<\eta<c_*$ and every $p>0$,
\begin{align}
 \sup_{N\geq1}\EE\exp\!\left(\eta\frac{H_N^2}{\log(2+H_N)}\right)
 &\leq1+\frac{4\eta}{c_*-\eta},\label{eq:intro-exp-moment}\\
 \sup_{N\geq1}\EE H_N^p
 &\leq2^p+4\Gamma(p+1)2^{29p}.\label{eq:intro-poly-moment}
\end{align}
In particular, the expectation in~\eqref{eq:intro-exp-moment} is at most five when $\eta=2^{-30}$.
These conclusions follow by integrating the upper tail; the proof, including the monotonicity of $F$ and the use of Tonelli's theorem before finiteness is known, is given in Corollary~\ref{cor:tails-moments}. The logarithm in the exponential moment is part of the assertion.

\subsection{Two-sided crossings and arm estimates}

Two-sided level-set crossings have a long history and connect several lines of work on the Gaussian free field, loop soups, and percolation. 
Since we consider it worthwhile to explain these connections and place our results in context, the next four subsections give a detailed account of the relevant literature.
%
%

\paragraph{Bulk crossing bounds.}

Gao, Nolin and Qian~\cite[Theorem 1.1, equation (1.3)]{GNQ} prove a positive limiting lower bound for the probability of a bulk crossing of a two-sided level set at height $C\sqrt{\log\log n}$. We refer here to version 2 of their paper, dated 27 May 2026. Their field is defined on $B_n$, and this bulk crossing is in $B_{n/2}$. Their Green convention~\cite[Section 2.3, equation (2.4)]{GNQ} includes the time-zero visit and kills on exit, with conductance $1/4$. Consequently $n=2N$ gives exactly the model and bulk geometry considered here. Corollary~\ref{cor:tightness} strengthens this bulk conclusion by allowing a height independent of $N$ at every prescribed confidence below one. Their discussion of constant-order thresholds, following the conjectural direction described in~\cite{AGS}, provides context for this improvement. We prove tightness of the bulk crossing height, not convergence to a deterministic threshold or identification of a limiting law. Their full-box crossing in equation (1.2) and the carpet condition in Theorem 1.10 impose different geometric requirements; our theorem concerns only the stated bulk event.

For the loop-soup comparisons, recall that a random-walk loop soup is a Poisson collection of random-walk loops; its clusters are the maximal families joined by chains of intersections~\cite[Section 2.2]{GNQ}. The intensity multiplies the loop measure. In the marked version used here, occupation at a vertex sums all loop holding times there; Section~\ref{sec:soups} fixes these conventions.

The scope of Gao--Nolin--Qian's work is broader than the bulk comparison alone. Their full-box theorem gives a crossing probability tending to one at the stated height scale. Their carpet theorem further locates the crossing in the part of the box outside the interiors of the outermost loop-soup clusters, giving geometric information in addition to the bound on the field. For subcritical loop-soup intensities, they establish a nontrivial phase transition at a constant occupation level; using the generalized second Ray--Knight theorem, they also derive bounds on crossings formed by thick points of a random walk. These are separate contributions of their method, and the bulk conclusion proved here does not include these additional assertions; see~\cite[Theorems 1.1, 1.4, 1.5 and 1.10]{GNQ}.

For the event common to the two papers, the order of the quantifiers is useful to emphasize. Given a desired confidence strictly below one, our corollary supplies a finite band height, depending on that confidence but not on the box size, at which the bulk crossing has at least that confidence uniformly in $N$. This strengthens the positive limiting lower bound at a growing height in~\cite[Theorem 1.1, equation (1.3)]{GNQ}, under the domain identification $n=2N$ above. The family of crossing heights is therefore tight. The fixed-height failure bound also rules out convergence in probability to a finite deterministic threshold; the existence and form of a limiting distribution are not established here.

\paragraph{Arm estimates and continuum inputs.}

The arm-event tools developed by Gao, Nolin and Qian in the first paper of their series~\cite{GNQI} explain an important part of the earlier crossing method. Their bulk four-arm event records two distinct outermost loop clusters crossing an annulus. A separation argument for random-walk packets makes it possible to cut and reconnect loops while keeping the relevant clusters apart. Their locality estimate controls the arm probability using loops in a surrounding box, and their upper quasi-multiplicative inequality bounds a crossing over several scales by a product of probabilities for separated scale ranges. These steps allow information from Brownian loop soups to be used quantitatively for random-walk loop soups; see~\cite[Section 4 and Propositions 1.1--1.2]{GNQI}.

The resulting discrete four-arm estimate holds throughout the subcritical and critical intensity range, with an arbitrarily small loss in the continuum power-law exponent; the constant may depend on that loss. This qualification matters when comparing it with the sharp critical estimates discussed below. The estimate is an upper bound, and we do not attribute a matching lower bound or an exact-exponent estimate up to constants to~\cite[Theorem 1.3 and Theorem 6.2]{GNQI}. These tools support the cluster analysis in~\cite{GNQ}; they are methodological predecessors, rather than arm theorems imported into the present proof.

The continuum estimates behind this program are developed further in Gao--Nolin--Qian's work on simple conformal loop ensembles (CLE)~\cite{GNQCLE}. For CLE parameters in $(8/3,4]$, their four-arm events allow the four crossings of an annulus to be made by one loop or by two loops, without prescribing the order or topology of those crossings. Theorem~1.1 gives upper and lower bounds up to multiplicative constants, both in the interior and at the boundary, with the corresponding explicit arm exponents. Theorem~1.2 is a separate ingredient: a separation estimate allowing two crossing loops to be placed in a suitably separated configuration at a bounded multiplicative cost. It supplies geometric control needed to pass from estimates with separated strands to these more general loop events.

An important upstream input is Zhan's work on two-curve Green's functions for $2$-SLE (Schramm--Loewner evolution with two curves)~\cite{ZhanI,ZhanB}. Those papers quantify the probability that both curves approach a specified interior or boundary point, and establish the associated Green's functions and asymptotics. Gao, Nolin and Qian explicitly use these results in their CLE analysis~\cite[Section 1 and Theorems 2.4--2.5]{GNQCLE}. Thus the roles are complementary: Zhan supplies quantitative two-curve estimates, while the CLE paper treats loop configurations and their separation. The passage from such continuum statements to discrete arm bounds requires the additional discrete analysis described above; none of these continuum arm assertions is itself a theorem about nearest-neighbor crossings of a discrete field band.

\paragraph{Sharp critical arm estimates.}

At the critical loop-soup intensity $1/2$, Bi, Gao, Nolin and Qian~\cite{BGNQ} obtain sharp arm estimates directly on the metric graph, formed by replacing lattice edges with line segments. Their Theorems~1.1--1.2 give matching upper and lower bounds: the bulk four-arm probability is comparable to the square of the inner-to-outer radius ratio, and the half-plane two-arm probability is comparable to that ratio itself. Corollary~1.3 supplies the corresponding discrete upper bounds. These critical estimates remove the exponent loss relevant to the earlier crossing argument. They are distinct from the DGFF arm and circuit results of Bi, Gao and Li~\cite{BiGaoLi}.

There is a substantive transfer issue here. Adding loops can merge two clusters, so the bulk four-arm event is not monotone in the loop configuration. Bi, Gao, Nolin and Qian therefore adapt their upper-bound proof to obtain the discrete bulk estimate; simple stochastic domination by the metric-graph soup does not suffice. Their sharp bounds support the improved critical crossing scale recorded in~\cite[Theorem 1.1]{GNQ}. This explains the role of~\cite{BGNQ} in the earlier program while keeping the metric-graph two-sided bounds, the discrete upper bounds, and the final band-crossing conclusion separate.

\paragraph{Path lengths at fixed band height.}

Gao and Zhang~\cite[Theorem 1.1]{GaoZhang} study the lengths of macroscopic paths in a fixed two-sided band. Their conclusion is a lower bound on the length of every such path when it exists; their published introduction raises the associated fixed-band existence question~\cite[p.~730]{GaoZhang}. Their convention pins the inner boundary of $V_n$; taking $n=2N+1$ makes the unpinned domain $B_{2N}$ and the bulk box $B_N$. Thus Corollary~\ref{cor:tightness} answers the corresponding fixed-band bulk existence question with crossing probability bounded away from zero (indeed with arbitrarily high prescribed confidence at a suitable fixed height), while remaining compatible with their obstruction to short paths. We do not bound chemical distance.

The content of Gao--Zhang's obstruction is geometric: for each fixed band height, paths whose endpoints are separated on the scale of the box must, with high probability, have length exceeding a strictly superlinear power of that scale. Equivalently, the probability that there is any sufficiently short macroscopic band path tends to zero~\cite[Theorem 1.1]{GaoZhang}. This unconditional statement permits an existence theorem at the same fixed height: the crossings supplied by such a theorem may have to make long detours. Their published discussion~\cite[p.~730]{GaoZhang} explicitly motivates the fixed-band existence question through this distinction. With the pinned/unpinned domain identification above, our corollary supplies existence at a suitable confidence-dependent height, while their theorem continues to constrain the lengths of those crossings. No length upper bound follows from our existence argument.

\subsection{Related geometric questions}

Exit sets give a complementary way to formulate geometric questions for fields with several components. Aru, Garban and Sep\'ulveda~\cite[Section 1.2]{AGS} explore inward from a zero boundary, continuing through vertices where the field norm is at most a prescribed cutoff and including the stopping vertices in the explored set. They also allow a bounded jump size. For target dimension at least two, their finite-volume Theorem~1.1 shows that, with the cutoff and jump size fixed, the probability that this exit set reaches a fixed positive proportion of the box size away from the boundary tends to zero.

Their infinite-plane Theorem~1.2 proves exponential decay of connection probabilities in the bounded-norm set, for the field pinned at the origin and for the massive field when the mass is sufficiently small. The scalar discussion in Remark~1 and Figure~3 instead points toward different behavior for large cutoffs. These results make the number of field components central to the comparison. They provide context for our scalar bulk-crossing question, with their own boundary conditions and exploration event; we do not use the vector theorem as an input to the scalar proof.

The distinction between scalar and vector-valued fields also appears in the study of exit sets by Aru, Garban and Sep\'ulveda~\cite{AGS}. Their Theorem 1.2 concerns pinned or massive infinite-plane fields with target dimension at least two, while their scalar discussion in Remark 1 and Figure 3 motivates two-sided connectivity questions. The present theorem addresses a finite-box scalar crossing with a prescribed bulk geometry. One-sided connectivity results, such as~\cite{DingLi,DingWirth}, and the arm and circuit estimates in~\cite{BiGaoLi}, concern other events. In particular, Drewitz, Elias, Pr\'evost, Tykesson and Viklund~\cite[Theorem 2.7]{DEPTV} establish a nontrivial bounded critical parameter for one-sided DGFF excursion connectivity in a disk; such a one-sided event does not imply the existence of a path in the two-sided band $\{|\phi_N|\leq a\}$.

\paragraph{One-sided connectivity and the choice of graph.}
The scalar/vector distinction also leads to a useful separation between several scalar connectivity problems. Ding and Li~\cite{DingLi} study chemical distances in planar DGFF level sets and in random walk loop-soup clusters, using tools that include Makarov's theorem, isomorphism theorems and entropic repulsion. Their work shows how geometric information about paths can be extracted from the field and its associated loops. Ding and Wirth~\cite{DingWirth} develop this approach for the metric-graph GFF; in particular, their planar results control chemical distance between annulus boundaries, conditional on their being connected. These are important precedents for studying the geometry of connected field regions, whereas the present crossing height records the smallest absolute-value bound allowing a path and makes no assertion about its length.

Ding, Wirth and Wu~\cite[Theorem 1.1]{DWW} give a particularly clear reason to distinguish discrete paths from metric-graph paths. For a zero-boundary GFF in a rectangle of fixed aspect ratio, the probability of a positive discrete crossing stays bounded away from both zero and one as the mesh tends to zero, while the corresponding metric-graph crossing probability tends to zero. The crossing endpoints in this statement lie on the specified interior arcs adjacent to the boundary; they are not vertices pinned to zero. On a metric graph, positivity must persist along each continuous edge segment, whereas the discrete event only tests the vertex values. Thus even closely related one-sided crossing events can have different asymptotic behavior when the graph interpretation changes.

There is also a substantial literature on one-sided connections to distant sets. Drewitz, Elias, Pr\'evost, Tykesson and Viklund~\cite{DEPTV} study vacant-set percolation for discrete and Brownian excursion clouds in the disk and connect this problem to DGFF excursion sets. Their Theorem~2.7 places the critical parameter for the stated one-sided DGFF connectivity problem strictly between zero and a finite upper bound. Bi, Gao and Li~\cite{BiGaoLi} obtain point-to-bulk connection asymptotics and point-to-boundary estimates for metric-graph level sets, together with discrete estimates and a comparison between the two models. Their arm and circuit methods provide further quantitative information on one-sided geometry. In each case the connection event, its endpoints and the boundary regime are part of the result; these works should be read alongside, rather than identified with, the bulk two-sided event defining $H_N$.

A different use of the term two-sided occurs in higher dimensions. Rodriguez~\cite{Rodriguez} studies level-set percolation for random interlacements and the GFF in dimension at least three, including the set where the absolute value of the GFF is \emph{above} a specified level. This high-absolute-value event retains both signs but excludes a band around zero. Our low-band event instead requires every visited vertex to remain in that band. This distinction is useful when comparing the present scalar problem with the broader percolation literature.

\paragraph{Two-valued local sets and continuum connectivity.}
Aru, Sep\'ulveda and Werner~\cite{ASW} construct and study bounded-type thin local sets of the two-dimensional continuum GFF, including the two-valued local sets. A local set is a random closed set whose revelation permits a domain Markov decomposition: in the remaining domain, the field is a harmonic function plus an independent zero-boundary GFF. For a two-valued local set, the harmonic function takes one of two prescribed values on each complementary component. This formulation gives geometric meaning to information about field heights even though the continuum GFF itself is a random distribution and has no pointwise values.

Aru and Sep\'ulveda~\cite[Proposition 3.12]{ASlocal} establish a precise boundary-connectivity criterion for these sets. For the zero-boundary GFF in the unit disk, write the set as $A_{-a,b}$, with $a,b>0$ and $a+b\geq2\lambda$. Their Green function has singularity $(2\pi)^{-1}\log(1/|z-w|)$ and $\lambda=\sqrt{\pi/8}$. If $a,b\geq2\lambda$, almost surely every pair of boundary points can be joined by a continuous path whose interior lies in the local set inside the disk. If either parameter is smaller than $2\lambda$, any fixed pair of distinct boundary points is almost surely not joined in this way. The order of these quantifiers is essential: the second assertion does not exclude exceptional random pairs. As emphasized in their Remark~3.14, this is a local-set formulation of continuum connectivity between two heights, not a pointwise continuum band. It is a close conceptual precedent for the present question. Relating its boundary-to-boundary criterion to our discrete bulk event would require additional approximation and event-identification arguments; its numerical threshold is not asserted here to be a threshold for $H_N$.

\paragraph{Contour lines, first passage sets and loop representations.}
The continuum local-set viewpoint belongs to a wider theory of curves coupled with the GFF. Schramm and Sheffield~\cite{SSDiscrete} prove that appropriate contour lines of an interpolated discrete GFF, with specified boundary values, converge to Schramm--Loewner evolution with parameter four (SLE$_4$). Their continuum paper~\cite{SSContinuum} constructs an SLE$_4$ contour line coupled with the distribution-valued GFF and characterizes the relation through a spatial Markov property. These works establish precise discrete and continuum foundations for speaking of GFF level lines. Miller and Sheffield's imaginary geometry~\cite{MillerSheffield} develops a broader coupling of the GFF with interacting SLE flow lines. It provides further context for the rich path geometry encoded by the field, while the two-valued local-set criterion above is the more specific connectivity comparison relevant here.

First passage sets connect this geometric theory with one-sided explorations and loop representations. Aru, Lupu and Sep\'ulveda~\cite{ALSbase} construct and characterize continuum first passage sets, including in finitely connected domains. Their companion paper~\cite{ALSfps} proves, under its domain and boundary assumptions, convergence of metric-graph first passage sets to the continuum sets in the Hausdorff topology. It also represents the continuum sets through clusters built from Brownian loops and boundary excursions. The two papers therefore supply both an intrinsic continuum object and a bridge from metric graphs to that object, with an isomorphism interpretation that makes the connection to loop soups explicit.

These achievements explain why contour and first passage methods are natural reference points for a two-sided crossing problem. They also specify what an argument using such methods would have to identify. The event defining $H_N$ concerns a nearest-neighbor path whose \emph{vertex values} all satisfy both height bounds. An interpolated contour can pass between vertices with large values of opposite signs. Likewise, convergence of two separately constructed one-sided sets does not by itself establish convergence of their intersection, or provide a single path satisfying both constraints. We use the continuum results as geometric precedents and sources of possible connections; the finite-graph arguments below prove the stated band event directly.

\subsection{Loop soups and massive methods}

The proof uses Le Jan's marked loop measure and root-excursion theory~\cite{LeJan}, the classical representation recorded in~\cite[author version, Proposition 2.4]{ABJL}, and Lupu's augmented loop-cluster coupling~\cite{Lupu}. Werner's spatial-resampling and switching identities, together with A\"id\'ekon's metric-graph explorations, give related descriptions of conditional loop geometry~\cite{WernerSpatial,WernerSwitch,AidekonExploration}. Our comparison uses a pivotal reachability event in the discrete nontrivial loop trace with directed incoming activations, while the exploration records complete marked high clusters and retains occupation bounds at low queries. Section~\ref{sec:comparison} specifies this joint information, derives its conditional law, and proves the excursion-deletion estimate for that event. We combine the resulting comparison with a connectivity estimate obtained from two auxiliary vertices, and then remove the mass using an independent Gaussian remainder.

\paragraph{From path representations to loop-cluster geometry.}
The relation between Markov paths and Gaussian fields has a history preceding the loop-soup formulation. Brydges, Fr\"ohlich and Spencer~\cite{BFS}, building on Symanzik's representations, express correlations of classical spin systems through random walks and use these expressions to obtain correlation inequalities. Dynkin~\cite{Dynkin83,Dynkin84} develops isomorphism identities relating Markov occupation fields to Gaussian and more general fields. These contributions explain the path--field correspondence behind many later constructions; their role should be distinguished from the specific statement that the occupation of a Poissonian soup is a squared Gaussian field.

For that soup identity we use Le Jan's work~\cite{LeJanOriginal,LeJanRenorm,LeJan}. In the finite symmetric setting, intensity-$1/2$ loop-soup occupation has the law of half the square of the GFF with the corresponding Green covariance; see in particular~\cite[Theorem 9(a)]{LeJanRenorm}. Lupu~\cite{Lupu} adds a geometric coupling in which signs are constant on loop clusters, and identifies the augmented discrete clusters through a metric-graph construction. This geometric information allows occupation identities to be used for connectivity questions. In the planar continuum, Lawler and Werner~\cite{LawlerWerner} construct the conformally invariant Brownian loop measure and soup that underlie the scaling-limit developments discussed in this overview.

Two complementary accounts are particularly helpful for the finite marked conventions used here. Chang and Le Jan~\cite[Sections 3.1--3.3 and 4.1]{ChangLeJan} develop Markov loop measures, time changes, excursions and Poisson ensembles directly from the generator, without requiring symmetry. Lawler and Perlman~\cite[Section 4, Theorem 2]{LawlerPerlman} present the occupation isomorphism by adding exponential holding times to discrete loops and including the trivial-loop occupation. Their construction makes the relation between discrete words and continuous occupation especially transparent. We retain the finite calculations in Section~\ref{sec:soups} to make cyclic multiplicities, orientation and holding rates visible at the places where the conditional argument uses them. These calculations record classical inputs in our conventions.

\paragraph{Massive Brownian loops and massive CLE.}
Exponential decay of massive loop clusters has an important continuum precedent. Camia~\cite[Theorem 2.4, second bullet]{Camia} proves that, for a full-plane massive Brownian loop soup of intensity at most $1/2$ and a mass function bounded away from zero, the diameter of the filled loop cluster containing a fixed deterministic point has an exponentially decaying tail. The decay bound is uniform in that point; the filled cluster is taken to be empty when no such cluster contains it. This is a statement about clusters, including connections formed by many loops, and not merely the suppression of individual long loops. It concerns a cluster attached to a specified point, rather than the largest cluster anywhere in the plane.

Papon~\cite[Theorems 1.3 and 1.8]{Papon} establishes the associated massive GFF and loop-ensemble coupling. Her setting includes bounded simply connected domains with bounded continuous nonnegative mass functions, and extends to the domain/mass pairs allowed by conformal covariance in her theorem. Reweighting the massless GFF--CLE$_4$ coupling gives a coupling of the massive GFF with a massive version of the conformal loop ensemble CLE$_4$. She then identifies these loops as the outer boundaries of the outermost clusters of a massive Brownian loop soup of intensity $1/2$. This proves the coupling conjectured by Camia and gives a continuum construction linking massive field geometry and massive soup clusters. We therefore treat that identification as an established result.

The mass scale used here is compatible with this continuum perspective. After a box of side of order $N$ is rescaled to a domain of fixed size, lattice mass $m=K/N$ corresponds to continuum mass of order $K$, with the numerical conversion depending on the generator convention. Thus the fact that the lattice mass tends to zero at fixed $K$ does not remove the relevance of massive continuum results. To use them for the estimate needed here, one would still have to control the finite augmented graph, identify its crossing event under the approximation, and obtain the required numerical bounds uniformly in $N$ and in the stated range of $K$. These are quantitative application questions; the continuum theorems motivate the comparison without themselves being cited as the all-size finite-graph estimate proved below.

\paragraph{Discrete-to-continuum loop and cluster comparisons.}
There are also several distinct levels of approximation connecting random walk loops to Brownian geometry. Lawler and Trujillo Ferreras~\cite{LawlerTF} introduce the random walk loop soup and prove a strong approximation by Brownian loop soup using a coupling of random walk and Brownian bridges. Van de Brug, Camia and Lis~\cite{BCL} then address the geometry of clusters: with their cutoff on small lattice loops, they prove convergence of loop clusters and of the outer boundaries of outermost clusters to the corresponding Brownian objects and conformal loop ensembles. Lupu~\cite[Theorem 1.1]{LupuCLE} removes this small-loop cutoff in the massless half-plane setting. For intensities in $(0,1/2]$, he proves convergence of the outer boundaries of outermost clusters for both discrete and metric-graph soups to the appropriate conformal loop ensemble; at intensity $1/2$ the limit is CLE$_4$. These are cluster-boundary results, substantially beyond matching individual loops.

Qian~\cite[Theorems 1.3 and 1.4]{QianCoupling} strengthens the underlying loop-soup coupling to every fixed polynomial scale. We refer to version 4, dated 30 August 2026. In the planar comparison, the original restriction to rescaled loop durations above $N^{\theta-2}$ with $\theta\in(2/3,2)$ is extended to every fixed $\theta\in(0,2)$. Her results treat both continuous-time and discrete-time random walk loop soups, and more generally dimensions at least one. Theorems~1.3 and~1.4 are the soup statements; the bridge coupling used in their proof is Theorem~1.2. Thus there is now much finer control of mesoscopic loops than in the original approximation. The cutoff exponent is fixed when applying these statements; the result is not an assertion that a finite lattice resolves arbitrarily small continuum loops.

This sequence of works separates three useful inputs: matching loops, controlling the geometry of their clusters, and obtaining quantitative finite-graph connection bounds. The first two are already supported by substantial theorems, and could be ingredients in an alternative continuum-based approach. The estimate used in the present proof additionally concerns the augmented graph of a finite massive soup and must hold with the stated constants for every permitted box and mass parameter. Establishing that estimate from continuum convergence would require the corresponding event identification and quantitative uniformity. Section~\ref{sec:probes} gives the finite-network argument and discusses this comparison at the point where massive connectivity is used.

\subsection{Broader connections}

\paragraph{Forest-fire methods.}
The arm-based approach to two-sided levels also has an instructive antecedent in forest-fire percolation. Kiss, Manolescu and Sidoravicius~\cite{KMS} prove that, on the planar lattices they consider, a fixed small amount of independent recovery does not restore percolation after the infinite occupied clusters have been deleted. Their analysis organizes possible reconnections through distinguished passage sites and across many spatial scales. Gao, Nolin and Qian explicitly explain this methodological ancestry in \cite[Section 1.4.2]{GNQ}: their loop-soup argument organizes chains of clusters through passage edges, whose geometry is controlled by arm events. More specifically, van den Berg and Nolin~\cite[Theorem 5.3]{vdBNForest} develop a stability estimate for six-arm events with passage sites using an induction on two spatial scales. The opening of \cite[Section 4]{GNQ} explicitly identifies that argument as an inspiration for the loop-soup induction. This connection explains how an argument can sum over chains involving clusters of very different sizes, and gives credit both to the passage-site method and to its subsequent inductive development. The forest-fire recovery theorem itself is not an input to the finite conditional-exploration and massive-field estimates used in this paper.

\paragraph{Physics predictions.}
A complementary perspective comes from the physics analysis of Cao and Santachiara~\cite{CaoSantachiara}. Using a loop-model description, they derive predictions for one-sided planar GFF level clusters, including logarithmic corrections to cluster area and connectivity and an exponentially diverging correlation length, and test their predictions numerically. Their analytical discussion uses a regularized field in a simply connected domain with Dirichlet boundary conditions; the corresponding boundary-connected cluster is their analogue of an infinite cluster. Their simulations use a periodic square lattice. These choices make the observable and the limiting question explicit, and help explain why logarithmic effects and slow finite-size behavior arise in planar GFF percolation. We cite this work for its theoretical predictions and numerical evidence. Its one-sided clusters differ from a single nearest-neighbor path whose every vertex satisfies both bounds defining our band, and its predictions do not supply the uniform finite-box tail estimate proved here.

\paragraph{Planar Gaussian crossings.}
The dependence of crossing probabilities on shape and level is also a central theme of general planar percolation. K\"ohler-Schindler and Tassion~\cite{KST} prove a Russo--Seymour--Welsh comparison for planar percolation processes satisfying suitable invariance and positive-association assumptions. Their result relates crossing a rectangle in its long direction to crossing it in its short direction, through a comparison that depends only on the aspect ratio and is uniform over scale and model. It identifies a robust geometric principle behind many crossing arguments: once the relevant hypotheses and an initial crossing estimate are available, changing the shape need not destroy quantitative control. In the present problem, checking those hypotheses for the field and for the chosen band event is a separate issue from recognizing this common geometric question.

For smooth stationary Gaussian fields, this program has led to substantial results about sign regions, nodal lines and excursion sets. Beffara and Gayet~\cite{BeffaraGayet} establish Russo--Seymour--Welsh estimates for sign domains and nodal lines in the real Bargmann--Fock setting. Muirhead and Vanneuville~\cite{MuirheadVanneuville} prove a phase transition at level zero for a broad class of smooth centered planar fields under symmetry, positivity, regularity and covariance-decay hypotheses, together with exponential decay of subcritical crossing probabilities. Their use of white-noise representations connects the Gaussian field to methods for sharp thresholds and approximate independence. Muirhead, Rivera and Vanneuville, with an appendix by K\"ohler-Schindler~\cite{MRV}, develop a further approach that does not require nonnegative correlations. Their work combines threshold delocalization with a version of Russo--Seymour--Welsh reasoning using sprinkling, and identifies level zero as critical under their smoothness, symmetry and decay assumptions. These works show how crossing theory can adapt to different mechanisms of dependence.

The distinction between an excursion and a band matters when relating this theory to our setting. A one-sided excursion event is monotone in the signed field, whereas membership in a centered bounded band is not: raising a field value can either enter or leave the permitted interval. Thus positive association for increasing functions of Gaussian coordinates does not by itself verify the association hypothesis for the induced band configuration. We do not assert that this configuration fails positive association; an application of that property would require its own argument. Moreover, our fields are defined in finite boxes with a prescribed zero boundary condition, and their covariances depend on the box. The smooth stationary covariance hypotheses in the preceding results therefore require a separate comparison before those theorems can be used for this family. These distinctions explain the role of the finite exploration and covariance estimates below.

\paragraph{Higher-dimensional connectivity.}
There is a related sharpness theory for strongly correlated fields in higher dimensions. For the massless lattice GFF in dimensions at least three, Duminil-Copin, Goswami, Rodriguez and Severo~\cite{DGRS} prove that three natural critical parameters coincide: those describing a well-controlled subcritical regime, the appearance of an infinite cluster, and local uniqueness in the supercritical regime. Muirhead~\cite{MuirheadStrong} subsequently develops a general sharpness theory for classes of Gaussian fields with algebraically decaying correlations, without requiring positive correlations; the higher-dimensional lattice GFF is among its examples. That paper explicitly credits the earlier GFF result of Duminil-Copin, Goswami, Rodriguez and Severo and also obtains further near-critical information. The shared issue is how to control connectivity despite long-range dependence. The infinite-volume higher-dimensional excursion problem and the present two-dimensional Dirichlet band problem have different covariance behavior and different events, so this theory supplies a useful comparison of methods and questions rather than the finite-box estimate used here.

Rodriguez, Shi and Zhang~\cite{RSZ} give a related use of mass for one-arm probabilities on cable systems of transient weighted graphs. Their upper-bound argument compares critical and massive cluster laws by an entropy estimate and combines this with a Simon--Lieb-type inequality derived from the Lupu--Werner connection formula~\cite{LWmetric} and the field's Markov property~\cite[Propositions 1.3, 3.1 and 4.1]{RSZ}. Their lower-bound argument relates cluster capacity and volume. This provides another connection between mass perturbations, random-walk estimates and connectivity. Their point-to-set upper bound assumes ellipticity and uniform upper bounds on total vertex conductances and diagonal Green functions. A transfer to the present planar finite-box family would need control of these hypotheses and constants as the box and mass vary, together with identification of the required band-crossing event. That transfer is not used in the proof below.

\paragraph{Gaussian decompositions.}
The Gaussian comparison used to remove the mass also belongs to a wider theory of decomposing covariance across spatial scales. Brydges, Guadagni and Mitter~\cite{BGM} construct finite-range decompositions of lattice resolvents in their setting of dimension at least three: each covariance component is positive semidefinite and vanishes beyond its assigned range, so the associated Gaussian components can be taken independent. Bauerschmidt~\cite{Bauerschmidt} develops a general construction for Green forms associated with symmetric Dirichlet forms, producing independent smoother components with localized correlations. Such decompositions are useful when one wants to separate spatial scales while retaining positivity of each covariance. They place the use of independent Gaussian components in an established methodological setting.

For continuum log-correlated fields, Junnila, Saksman and Webb~\cite{JSW} make a related distinction between a regular comparison and an independent decomposition. Their Theorem A couples two Gaussian fields whose covariance difference has suitable Sobolev regularity so that their difference is a H\"older-continuous Gaussian process on interior subdomains; that process need not be independent of either field. Their Theorem B gives a local decomposition of a log-correlated field into an almost $\star$-scale invariant field and an independent H\"older-continuous remainder under the stated regularity assumptions. The two results serve different purposes: a regular coupled difference permits one kind of comparison, while independence supplies additional structure. In particular, smoothness of a covariance difference alone does not say that the difference is positive semidefinite and hence the covariance of an independent correction.

Our choice of decomposition is dictated by the input needed in the loop-soup argument. The first component must have the exact finite-domain massive GFF law used in Sections~\ref{sec:comparison}--\ref{sec:probes}. Replacing it by a finite-range component or by the scale-invariant field of a continuum decomposition would require a new comparison with that law. Instead, Section~\ref{sec:gaussian} verifies positivity of the difference between the massless and massive resolvents directly and then estimates the independent remainder. A close precedent for this massive-to-massless strategy is the work of Rodriguez and Zhang~\cite[version 2, equations (3.12)--(3.13) and Lemma 3.7]{RZ}; the local discussion in Section~\ref{sec:gaussian} explains the parameter range retained here. The elementary covariance identity and the subsequent Gaussian-process estimates are used as established tools, with their finite-box constants and mass dependence kept explicit.

\paragraph{Further reading.}
For readers approaching these connections from different backgrounds, Sheffield's account~\cite{Sheffield} develops the Gaussian Hilbert-space viewpoint, the domain Markov property and the discrete and massive free fields. Werner and Powell~\cite{WP} give a gradual introduction to the discrete and continuum GFF and their connections with Brownian loop soups and conformal loop ensembles. These expositions provide useful routes through the objects appearing in this overview. The original papers cited above, and the precise references placed beside the inputs in the later sections, identify the sources of the individual constructions and results.

\paragraph{The role of the present arguments.}
The main result supplies uniform quantitative control of the bulk crossing height. Its proof uses established occupation and root-excursion identities, Lupu's geometric coupling, and standard Poisson and Gaussian tools~\cite{LeJan,ABJL,Lupu,LastPenrose}. The finite comparison leading to Theorem~\ref{thm:routes-M1} checks how these inputs apply when the exploration retains complete marked high clusters and earlier occupation bounds. Its one-excursion certificate is verified for the specified incoming-activation pivotal event. The separate block construction supplies destination variables independent of the entire nontrivial marked soup, with the substitution made after the conditional comparison has been averaged.

The massive connectivity bound in Theorem~\ref{thm:probes-M2} uses a separate finite-network argument. Two auxiliary vertices attach to a path selected from the baseline configuration; the proof checks the independence of the visiting-loop classes after that selection and obtains quantitative attachment estimates. A Schur complement and Lupu's established two-point connection identity then control the probe connection. The mass-removal step estimates an independent Gaussian correction with the required mass dependence. The lower bound follows independently from planar duality, a known Gaussian translation inequality and a network-energy bound. Thus the detailed proofs identify the conditional constructions and quantitative applications while retaining the classical status of their inputs.

\paragraph{Note added after completion.}
After completing an initial draft of the present paper, we shared it with Yifan Gao. He subsequently informed us that he and his coauthors had also obtained a solution to the same problem, using different methods, and that their manuscript was in preparation. We have not yet seen their manuscript and therefore cannot comment on the relationship between the two approaches. We thank him for informing us of their work.

\section{Proof overview and reductions}\label{sec:overview}

The upper bound separates into three estimates. The first converts failure of a low crossing for a massive field into connectivity in its augmented loop graph. The second makes this connectivity small when the mass scale is short compared with the box. The third controls the independent Gaussian field added back when the mass is removed. We state these estimates here with their parameter ranges before constructing their ingredients.

\subsection{The massive comparison and connectivity estimates}

For $0<m\leq1$, let
\begin{equation}\label{eq:massive}
 \psi_{N,m}\sim\mathcal N\bigl(0,(\mathsf Q_D+m^2I)^{-1}\bigr).
\end{equation}
Let $\mathcal G_{N,m}$ denote the augmented loop graph on $D$ for this precision matrix, at soup intensity $1/2$. Its exact definition, including the additional edges conditional on the full marked soup, is given in Section~\ref{sec:soups}. Write
\[
 \operatorname{TOP}_N=\{(i,N):-N\leq i\leq N\},\qquad
 \operatorname{BOT}_N=\{(i,-N):-N\leq i\leq N\},
\]
and define
\begin{equation}\label{eq:theta}
 \theta_N(m)=\PP\bigl(\operatorname{TOP}_N
 \longleftrightarrow\operatorname{BOT}_N\text{ in }\mathcal G_{N,m}\bigr).
\end{equation}
This connection is allowed to leave $B_N$ and use all of $D$.

The two estimates are
\begin{align}
 \PP\bigl(H_N(\psi_{N,m})>256\bigr)
 &\leq\theta_N(m),
 &&N\geq1,\quad0<m\leq1,\tag{M1}\label{eq:M1-overview}\\
 \theta_N(K/N)
 &\leq\min\left\{1,\frac{2^{28}}3 K^3e^{-K/2}\right\},
 &&K\geq8,\quad N\geq K.\tag{M2}\label{eq:M2-overview}
\end{align}
Here and below $N$ is an integer, and the scale $K=mN$ may be real. Estimate~\eqref{eq:M1-overview} is Theorem~\ref{thm:routes-M1}; estimate~\eqref{eq:M2-overview} is Theorem~\ref{thm:probes-M2}.

For~\eqref{eq:M1-overview}, the basic planar alternative says that absence of a nearest-neighbor left--right crossing in a closed low band is equivalent to a top--bottom star crossing in its strict complement. Star adjacency means $|x-y|_\infty=1$ for distinct lattice vertices; nearest-neighbor adjacency means $|x-y|_1=1$. The exact finite-box alternative is stated and proved in Lemma~\ref{lem:routes-duality}. At intensity $1/2$, full loop occupation has the law of one half of the square of the matching field. It is the full occupation, including the trivial loops, that enters the threshold event.

To prove~\eqref{eq:M1-overview}, we explore clusters in the nontrivial loop trace, starting from the top side. At a new star-neighbor, a low observation records only that its occupation is at most the threshold. A high observation reveals the complete marked configuration of the new cluster. Thus the conditional law of the remaining marked loops must retain all earlier low-occupation inequalities. Sections~\ref{sec:soups} and~\ref{sec:comparison} establish this law. We then show that, on the pivotal event, keeping one selected root excursion suffices after further deletions. A Gamma tail of shape $3/2$ controls the corresponding occupation and yields a finite supermartingale comparison with independent site variables that add incoming edges.

Section~\ref{sec:routes} constructs the required independent site variables from Gamma occupation clocks and Poisson edge-opening processes. Since the construction uses the trivial occupation clocks, we make this substitution only after the conditional comparison has been integrated to an event involving the nontrivial loop trace and independent Bernoulli variables. The augmented graph keeps its exact law. This proves~\eqref{eq:M1-overview} without claiming that the constructed site variables are independent of the trivial clocks.

Estimate~\eqref{eq:M2-overview} uses a different argument. We add two auxiliary vertices to the massive network, reducing the old killing rates by the added conductances so that the old diagonal rates are unchanged. A selected crossing has suitable portions near the top and bottom, even if its path travels throughout $D$. Green function, capacity, and hitting estimates give each auxiliary vertex a controlled chance to connect to the appropriate portion. The two classes of loops used for these connections are disjoint, and their Poisson counts are independent conditional on the original configuration used to select the path. We compute the enlarged Gaussian covariance by a Schur complement and then use Lupu's augmented-cluster identity to bound the probability that the auxiliary vertices connect. Section~\ref{sec:probes} gives the details. This argument does not use~\eqref{eq:M1-overview}.

\subsection{Removing the mass}

The independent Gaussian correction is an established comparison device, with a closely related planar implementation in Rodriguez and Zhang~\cite{RZ}. The estimates here retain the mass dependence needed when the comparison scale is chosen according to the desired crossing height.

Since $\mathsf Q_D$ is positive definite, the matrix
\[
 \mathsf Q_D^{-1}-(\mathsf Q_D+m^2I)^{-1}
\]
is positive semidefinite. We may therefore construct independent centered Gaussian fields $\psi_{N,m}$ and $\zeta$, with covariance matrices
\begin{equation}\label{eq:overview-split}
 (\mathsf Q_D+m^2I)^{-1}
 \quad\text{and}\quad
 \mathsf Q_D^{-1}-(\mathsf Q_D+m^2I)^{-1},
\end{equation}
respectively. Their sum has the law of $\phi_N$. When $m=K/N$, we write $\zeta_{N,K}$ for this remainder, indexing it by the dimensionless scale $K$. For any two functions $f,g$ on $D$,
\begin{equation}\label{eq:deterministic-H}
 H_N(f+g)\leq H_N(f)+\|g\|_{\infty,B_N}.
\end{equation}
Indeed, choose a minimizing path for $f$ in the finite set $\LR_{\NN}(B_N)$ and use the triangle inequality along that path.

Put $\ell_K=\log(2K)$. The Gaussian estimate proved in Theorem~\ref{thm:gaussian-remainder} is
\begin{equation}\label{eq:remainder-overview}
 \PP\!\left(\|\zeta_{N,K}\|_{\infty,B_N}
 >800\{\ell_K+t\sqrt{\ell_K}\}\right)
 \leq2e^{-t^2},
 \qquad1\leq K\leq N,\quad t\geq0.
\end{equation}
The field is defined on $D$, while this supremum is over $B_N$. A finite sine expansion gives the point and increment bounds. A net starting at the mass scale, with actual integer grid vertices and explicit parent maps, then proves~\eqref{eq:remainder-overview}. All frequency cutoffs and the parameter $K$ may be real.

Using~\eqref{eq:M1-overview}, \eqref{eq:M2-overview} and~\eqref{eq:remainder-overview} with~\eqref{eq:deterministic-H} gives
\begin{align}\label{eq:overview-combined}
 \PP\!\left(H_N>256+800\{\ell_K+t\sqrt{\ell_K}\}\right)
 &\leq \frac{2^{28}}3 K^3e^{-K/2}+2e^{-t^2},\nonumber\\[-6pt]
 &\hspace{25mm}K\geq8,\quad N\geq K,\quad t\geq0.
\end{align}
For qualitative tightness as $N\to\infty$, one may first choose $K$ and then $t$. To obtain the uniform tail, Section~\ref{sec:gaussian} takes, for $a\geq2^{16}$,
\[
 K=\lceil a^2\rceil,\qquad
 t=\frac{a}{1600\sqrt{\log(2K)}}.
\]
The comparison applies when $N\geq K$. For $N<K$, a finite-box Green bound and a Gaussian union bound give the required estimate directly. The remaining bounded range of heights is covered by the probability bound one. These cases together give the assertion for every $a\geq0$ and $N\geq1$.

\subsection{A lower bound at fixed height}

For the lower bound, a strict positive top--bottom star crossing is a barrier to a nearest-neighbor left--right path in the band $[-a,a]$. Symmetry and finite planar duality give a positive star crossing at level zero with probability at least $1/2$. A deterministic shift equal to $a$ on $B_N$ converts that crossing into a strict barrier above $a$.

The cost of the shift is controlled by its energy in the original killed network. The finite Gaussian translation inequality of Li and Kuelbs~\cite[Theorem 1, equation (1.1)]{KuelbsLi}, also obtained here through the Neyman--Pearson likelihood-ratio principle, converts this energy into a normal tail. We use the equilibrium potential and the classical Nash--Williams cutset bound for the square shells. The bound includes every shell through the killing boundary and gives an energy upper bound less than $2/\log2$. Section~\ref{sec:lower} proves the resulting lower tail. This argument uses neither of the massive estimates nor the Gaussian remainder bound.

\paragraph{Organization.}
Sections~\ref{sec:soups}--\ref{sec:routes} develop the marked exploration and prove~\eqref{eq:M1-overview}. Section~\ref{sec:probes} contains the network extension, potential estimates, and probe calculation proving~\eqref{eq:M2-overview}. Section~\ref{sec:gaussian} combines these bounds and proves the upper tail and moments; Appendix~\ref{app:gaussian} supplies its finite spectral and net calculations. Section~\ref{sec:lower} can be read independently once the model and the planar crossing alternative are fixed.

\section{Finite marked loop soups}\label{sec:soups}\label{sec:finite-networks}

This section fixes the marked-loop model and the normalization of occupation times. We then state the classical root-excursion law in the finite form needed for the conditional comparison. The cyclic-loop measure, trivial occupations, and root Gamma--Poisson decomposition are standard results from \cite[Sections 2 and 7.2]{LeJan}, also summarized in \cite[author version, Proposition 2.4]{ABJL}. We include short calculations to fix the factors arising from periodic loops and holding rates. Section~\ref{sec:probes} later uses the same conventions when two auxiliary vertices are added.

Throughout, $\GammaDist(\alpha,\text{rate }\beta)$, for $\alpha,\beta>0$, has density $\beta^\alpha t^{\alpha-1}e^{-\beta t}/\Gamma(\alpha)$ on $(0,\infty)$.

\subsection{Finite networks and the marked loop measure}

All network data are deterministic. Let $V$ be a finite nonempty vertex set. Conductances satisfy $c_{uv}=c_{vu}\geq0$ for distinct vertices and $c_{vv}=0$. Killing rates satisfy $k_v\geq0$. Put
\[
 a_v=k_v+\sum_{u\in V}c_{vu},\qquad
 A=\operatorname{diag}(a_v),\qquad C=(c_{uv}),\qquad \mathsf Q=A-C.
\]
Assume $\mathsf Q$ is positive definite. Then every $a_v>0$, and $P=A^{-1}C$ is the transition matrix of the killed embedded jump chain. Its spectral radius is strictly less than one: $P$ is similar to the symmetric matrix $A^{-1/2}CA^{-1/2}$, whose largest eigenvalue is less than one by $\mathsf Q>0$; nonnegativity of $P$ identifies that eigenvalue with its spectral radius. The continuous-time generator on $V$ is $-\mathsf Q$, so its occupation Green matrix is $\mathsf Q^{-1}$.

A nontrivial oriented cyclic marked loop is written
\[
 [\mathbf x;\mathbf s]=[(x_0,\ldots,x_{n-1});(s_0,\ldots,s_{n-1})],
 \quad n\geq2,\quad x_n=x_0,\quad s_i>0.
\]
We identify simultaneous cyclic shifts of the vertex and holding-time sequences. Reversal is not an additional identification. The repeated endpoint $x_n=x_0$ is not an extra cyclic occurrence. Each occurrence contributes its holding time to the occupation of its vertex.

The following is Le Jan's marked cyclic measure; see \cite[author version, Sections 2.2--2.3, equations (2.4), (2.8)--(2.9)]{LeJan}. We spell out the normalization, including periodic words and zero-jump loops.

The generator-based treatment in Chang--Le Jan~\cite[Sections 3.1--3.3]{ChangLeJan} and the discrete-loop marking construction in Lawler--Perlman~\cite[Section 4]{LawlerPerlman} give complementary formulations of this background. The present calculation keeps the holding rates explicit and fixes the cyclic multiplicity and orientation conventions used throughout. The historical distinction between the occupation isomorphism and its geometric coupling is discussed in Section~\ref{sec:introduction}.

\begin{lemma}[Exact marked measure and trivial occupations]\label{lem:finite-marked}
For any nonnegative measurable cyclic-invariant function $F$, the nontrivial loop measure is
\begin{equation}\label{eq:finite-marked}
 \int F\,d\mu_{\mathsf Q}^{\mathrm{nt}}
 =\sum_{n\geq2}\frac1n
 \sum_{\substack{x_0,\ldots,x_{n-1}\in V\\x_n=x_0}}
 \int_{(0,\infty)^n}F([\mathbf x;\mathbf s])
 \prod_{i=0}^{n-1}\left(c_{x_ix_{i+1}}e^{-a_{x_i}s_i}\,ds_i\right).
\end{equation}
Equivalently, the product in the integral is the product of transition weights $P(x_i,x_{i+1})$ and independent exponential holding-time densities of rates $a_{x_i}$. An unmarked oriented loop which repeats its primitive word $r$ times has measure weight $r^{-1}\prod_iP(x_i,x_{i+1})$.

The measure of trivial loops of duration $t$ at $v$ is $e^{-a_vt}\,dt/t$. In an intensity-$\alpha$ soup, their summed occupation $\xi_v$ has law $\GammaDist(\alpha,\text{rate }a_v)$. These occupations are mutually independent and independent of the nontrivial soup.
\end{lemma}
\begin{proof}
We treat the nontrivial cyclic loops first and the zero-jump loops afterward. For the cyclic loops, we need to recover the factor coming from the choice of a starting point, including when a word is periodic.

Use the rooted continuous-time loop measure with duration factor $dt/t$, as in \cite[(1.1)]{Lupu}. For a loop with $n$ jumps and initial vertex $x_0$, its initial and final holding pieces lie at the same vertex. Combine them into $s_0$. For fixed total duration $t=\sum_i s_i$, integration over the position of the time origin inside that combined interval contributes $s_0$. Thus the rooted density, tested against $F$, has factor $s_0/t$ multiplying the product in~\eqref{eq:finite-marked}.

The sum over all closed words and the holding-time integral are invariant under simultaneous cyclic rotations. Averaging the factor $s_0/t$ over the $n$ positions gives
\[
 \frac1n\sum_{i=0}^{n-1}\frac{s_i}{\sum_j s_j}=\frac1n,
\]
which proves~\eqref{eq:finite-marked}. This averaging does not assume that the word is primitive. If its minimal period is $d$ and $n=rd$, its unmarked cyclic orbit has $d$ distinct rooted words, so its total factor is $d/n=1/r$. Holding marks are attached to occurrences, including repeated occurrences of the same vertex. No extra factor for reversal is introduced.

With the cyclic factor fixed, we can now read off the holding-time law from the same product density.

The identity $c_{uv}e^{-a_us}\,ds=P(u,v)a_ue^{-a_us}\,ds$ gives the stated exponential marking law. The total nontrivial measure is $\sum_{n\geq2}\operatorname{Tr}(P^n)/n<\infty$. Hence the nontrivial soup has finitely many loops almost surely, each with finite marks.

We next consider the zero-jump loops. Their total occupation gives the Gamma law in the statement.

A loop with no jump remains at $v$ with weight $e^{-a_vt}$, which gives the asserted trivial-loop measure. By the Poisson Laplace functional \cite[Theorem 3.9]{LastPenrose}, its occupation has Laplace transform, for $z\geq0$,
\begin{align*}
 \EE e^{-z\xi_v}
 &=\exp\!\left\{-\alpha\int_0^\infty
       (1-e^{-zt})e^{-a_vt}\frac{dt}{t}\right\}\\
 &=\left(\frac{a_v}{a_v+z}\right)^\alpha.
\end{align*}
The integral equals $\log((a_v+z)/a_v)$, as follows by integrating $\int_0^z(a_v+u)^{-1}\,du$. This is the Gamma transform. Independence follows by restricting a Poisson process to the disjoint trivial-vertex and nontrivial classes. The graph construction below uses only the total trivial-loop occupation at each vertex, not the individual trivial loops.
\end{proof}

\subsection{The augmented graph}

For an intensity-$1/2$ marked soup, let $X_v$ be full occupation and $n_e$ the number of loop traversals of the unordered edge $e$. For independent uniforms $U_e$, independent of the entire soup, define
\begin{equation}\label{eq:finite-augmentation}
 \eta_e=\ind_{\{n_e>0\}}\vee
 \ind_{\{U_e\leq1-\exp(-2c_{uv}\sqrt{X_uX_v})\}},
 \qquad e=\{u,v\},\ c_{uv}>0.
\end{equation}
Write $\mathsf{Aug}(\mathsf Q)$ for this graph law. The construction of \cite[Theorem 1 bis and Corollary 3.6]{Lupu} gives $X\overset d=\phi^2/2$ for a Gaussian field with covariance $\mathsf Q^{-1}$. Here conductances are jump rates, $X$ includes trivial-loop occupation, and the additional edge openings are independent conditional on the entire marked soup. These are the conventions used in~\eqref{eq:finite-augmentation}.

Lupu states a connected-network setup. For a disconnected finite network, apply the input to each non-singleton component, whose killing makes it transient. For a singleton with precision $a>0$, $\xi\sim\GammaDist(1/2,\text{rate }a)$ and an independent fair sign times $\sqrt{2\xi}$ has law $\mathcal N(0,1/a)$; the graph has no edges. This extends the construction to every finite positive definite precision network considered here. Lemma~\ref{lem:probes-arcsine} will apply the connection identity in \cite[proof of Proposition 5.2]{Lupu} to this same augmented partition.

\paragraph{Graph and intensity conventions.}
The discrete trace graph retains the edges traversed by nontrivial loops. The augmented graph also uses the full occupations and the additional edge randomness in~\eqref{eq:finite-augmentation}. Under Lupu's construction~\cite{Lupu}, its connectivity between original vertices agrees with the corresponding cable-graph loop-cluster connectivity. Thus raw trace connectivity, augmented connectivity and continuous paths on the cable graph refer to related but distinct objects. The incoming directed arcs of Section~\ref{sec:comparison} are a further auxiliary construction; Section~\ref{sec:routes} realizes them by routes in the augmented graph.

For comparison with the planar continuum literature, Sheffield and Werner~\cite{SheffieldWerner} establish the loop-soup construction and Markovian characterization of simple conformal loop ensembles. The intensity convention matters when using that connection. With the Brownian loop measure convention used in Lupu's account~\cite[Introduction]{LupuCLE}, central charge $c$ and soup intensity $\alpha$ satisfy $c=2\alpha$; in particular CLE$_4$ corresponds to $\alpha=1/2$. Lupu explicitly explains the historical factor-two correction. Throughout the finite calculations here, the measure is the oriented cyclic measure specified above, reversal is not an extra identification, and the occupation identity is $X\overset d=\phi^2/2$ at intensity $1/2$. These conventions, rather than an untranslated intensity parameter from another source, fix the factors in the augmentation and probe calculations.

\subsection{Uniform time rescaling}

Multiply every conductance and killing rate of a network by $r>0$. Its embedded jump probabilities are unchanged, and its marked soup can be coupled by dividing every holding time by $r$. Its occupations become $X/r$. The trivial-loop measure is transformed in exactly the same way. Moreover
\[
 2(rc_{uv})\sqrt{(X_u/r)(X_v/r)}=2c_{uv}\sqrt{X_uX_v}.
\]
Using the same uniforms therefore leaves the entire augmented graph unchanged. In the massive application, taking $r=(1+m^2)^{-1}$ produces the rate-one normalization for $\chi=\sqrt{1+m^2}\,\psi_{N,m}$. Thus the normalized graph used in Section~\ref{sec:routes} has the same law as the unnormalized massive graph defining $\theta_N(m)$ in~\eqref{eq:theta}.

\paragraph{General intensity and permanental fields.}
The occupation construction is meaningful at every positive soup intensity, although the geometric GFF application uses intensity $1/2$. At general intensity its finite-dimensional Laplace transforms are powers of determinants, placing the occupation field in the family of infinitely divisible permanental fields associated with Markov processes. Eisenbaum and Kaspi~\cite{EisenbaumKaspi} develop this connection and corresponding occupation identities beyond the squared symmetric Gaussian case. Their parameter is specified through its determinant exponent; throughout the present paper, $\alpha$ remains the intensity multiplying the loop measure. This broader viewpoint explains why the root calculation is naturally stated for general $\alpha$.

\subsection{Root excursions and their intensity}\label{subsec:root-excursions}

The root decomposition is classical; see \cite[author version, Section 7.2, p.~73]{LeJan} and \cite[author version, Proposition 2.4]{ABJL}. We give the marked calculation because the conditional comparison needs its full joint law, with the holding rates retained.

Let $R$ be a finite nonempty network with the conductances, killing, rates and positive definite precision $\mathsf Q_R$ specified above. Fix a root $u\in R$, and set
\[
 S=R\setminus\{u\},\qquad a=a_u,\qquad
 g=(\mathsf Q_R^{-1})_{uu},\qquad \mathbf c=(c_{uv})_{v\in S}.
\]
The principal submatrix $\mathsf Q_{SS}$ keeps the rates $a_v$ at $v\in S$ and kills the chain upon an attempted jump to $u$.

An excursion mark is
\[
 \epsilon=(u,x_1,\ldots,x_k,u;\eta_1,\ldots,\eta_k),
 \quad k\geq1,\quad x_i\in S,\quad \eta_i>0.
\]
The intermediate word may revisit any nonroot vertex. Each occurrence has its own holding mark. The path is oriented, and reversal is not an identification. No individual root holding time is included in $\epsilon$. The excursion space $\mathsf E_u$ is the countable disjoint union of these finite word spaces with their Euclidean holding coordinates and Borel sigma-fields. Words with a zero-conductance transition have measure zero. Define the finite measure $\rho$ by
\begin{equation}\label{eq:root-rho}
 \rho(d\epsilon)=c_{ux_1}
 \left(\prod_{i=1}^{k-1}c_{x_ix_{i+1}}\right)c_{x_ku}
 \prod_{i=1}^{k}e^{-a_{x_i}\eta_i}\,d\eta_i
\end{equation}
on each word component, and write $\rho(h)=\int h\,d\rho$. Each mark retains the complete nonroot holding data of one excursion occurrence; trivial root loops are represented separately.

\begin{lemma}[Excursion mass and Green conversion]\label{lem:root-rho}
Writing $b=\rho(\mathsf E_u)$ and $\beta=g^{-1}$, one has
\begin{equation}\label{eq:root-mass}
 b=\mathbf c^{\mathsf T}\mathsf Q_{SS}^{-1}\mathbf c=a-g^{-1},
 \qquad 0\leq b<a,\qquad \beta=a-b>0.
\end{equation}
For $S=\varnothing$, the quadratic form and $\rho$ are zero. If $\nu$ denotes the subprobability law of the first-return path to $u$ before killing, with its initial root hold omitted and its nonroot holds retained, then $\rho=a\nu$. In particular $b/a$, not $b$, is the return probability in general.
\end{lemma}
\begin{proof}
We first sum the weights of excursions that leave the root and return to it before killing. A Schur complement identifies their total mass. We then compare these weights with first-return probabilities for the killed jump chain.

Integrating the holding coordinates in~\eqref{eq:root-rho} and summing words gives
\[
 b=\mathbf c^{\mathsf T}\left[\sum_{j\geq0}
 (A_S^{-1}C_{SS})^jA_S^{-1}\right]\mathbf c
 =\mathbf c^{\mathsf T}\mathsf Q_{SS}^{-1}\mathbf c,
 \qquad A_S=\operatorname{diag}(a_v:v\in S).
\]
The series converges because $\mathsf Q_{SS}$ is positive definite and its embedded killed transition matrix has spectral radius less than one, by the argument in Section~\ref{sec:finite-networks}. The block matrix
\[
 \mathsf Q_R=\begin{pmatrix}a&-\mathbf c^{\mathsf T}\\-\mathbf c&\mathsf Q_{SS}\end{pmatrix}
\]
has positive Schur complement $a-b$, and block inversion yields $g=(a-b)^{-1}$. This also proves finiteness of $\rho$. It remains to compare the marked excursion weights with first-return probabilities. The difference occurs at the initial departure from the root, whose holding time is omitted from the excursion mark.

For a fixed excursion word, the initial transition probability is $c_{ux_1}/a$. At each intermediate vertex, the transition probability times its exponential holding density is
$(c_{x_ix_{i+1}}/a_{x_i})a_{x_i}e^{-a_{x_i}\eta_i}\,d\eta_i$.
The same formula applies to the final return jump. Multiplication shows that the first-return subprobability density is exactly $\rho/a$.
\end{proof}

For a rate-one root,~\eqref{eq:root-mass} becomes $b=1-g^{-1}$. The factor $a$ in $\rho=a\nu$ is essential when the root rate is not one.

\paragraph{Root cycles and the avoiding soup.}

Use the original intensity-$\alpha$ marked loop soup. Its nontrivial loops have finite total intensity by Lemma~\ref{lem:finite-marked}. Partition them into loops avoiding $u$ and loops visiting $u$. Let $\mathcal B$ consist of all avoiding nontrivial marked loops and all trivial occupation clocks at vertices of $S$. This background has the marked-soup law for $\mathsf Q_{SS}$: each allowed word keeps the same conductances and holding rates, and each old trivial clock keeps its rate. It is independent of the entire root-visiting nontrivial soup and of the trivial root clock
\[
 \xi_0\sim\GammaDist(\alpha,\text{rate }a).
\]
These assertions follow from restriction to deterministic disjoint Poisson classes, including the trivial classes.

Cut each root-visiting loop at every cyclic occurrence of $u$. A loop with $r$ such occurrences yields a cyclic list of $r$ excursions $\epsilon_1,\ldots,\epsilon_r$ and $r$ root holds $\tau_1,\ldots,\tau_r$, with $\tau_i$ immediately before $\epsilon_i$. The total root occupation and aggregate excursion point measure are
\[
 T=\xi_0+\sum_{\text{root-visiting loops }l}\sum_{i=1}^{r(l)}\tau_i(l),
 \qquad
 \mathcal E=\sum_{\text{root-visiting loops }l}\sum_{i=1}^{r(l)}\delta_{\epsilon_i(l)}.
\]
The multiset retains multiplicities. It forgets the cyclic grouping into original loops and the separate root holding times, while $T$ retains their sum together with $\xi_0$.

\begin{lemma}[Exact root-cycle factor]\label{lem:root-cycle}
For a nonnegative measurable function $F$ on cyclic marked loops, the contribution of loops with exactly $r\geq1$ root occurrences is
\begin{equation}\label{eq:root-root-cycle}
 \int_{\{r(l)=r\}}F(l)\,\mu_{\mathsf Q_R}^{\mathrm{nt}}(dl)
 =\frac1r\int F\bigl([\tau_1,\epsilon_1;\ldots;\tau_r,\epsilon_r]\bigr)
 \prod_{i=1}^{r}\bigl(e^{-a\tau_i}\,d\tau_i\,\rho(d\epsilon_i)\bigr).
\end{equation}
In particular, for $s\geq0$ and measurable $0\leq h\leq1$,
\begin{equation}\label{eq:root-cycle-test}
 \int_{\{r(l)=r\}}e^{-s\sum_i\tau_i}\prod_i h(\epsilon_i)
 \,\mu_{\mathsf Q_R}^{\mathrm{nt}}(dl)
 =\frac1r\left(\frac{\rho(h)}{a+s}\right)^r.
\end{equation}
Periodic words and repeated root visits use this same factor $1/r$.
\end{lemma}
\begin{proof}
The factor $1/r$ comes from starting the cyclic word at one of its root occurrences. Counting occurrences, rather than distinct rotations, also handles periodic words.

Start with the length-$n$ formula~\eqref{eq:finite-marked}, restricted to words with $r$ root occurrences. The word sum and holding integral are invariant under rotation. Averaging $\ind_{\{x_0=u\}}/r$ over all $n$ cyclic positions gives $1/n$, because exactly $r$ positions are $u$. Thus the original factor $1/n$ with all starting positions can be replaced by $1/r$ with $x_0=u$.

The resulting word, starting at $u$, has a unique ordered decomposition into $r$ first-return excursions and their preceding root holds. Its product density factors into the right side of~\eqref{eq:root-root-cycle}. Summing over all possible lengths proves that identity. This rotation argument counts positions rather than assuming a free cyclic action, so it also covers periodic words. Orientation is unchanged and no reversal factor is inserted. Finally, each root hold contributes $\int_0^\infty e^{-(a+s)\tau}\,d\tau=(a+s)^{-1}$, and each excursion contributes $\rho(h)$, proving~\eqref{eq:root-cycle-test}.
\end{proof}

Taking $s=0$ and $h=1$ shows that the mass of root-visiting nontrivial loops is
\begin{equation}\label{eq:root-root-mass}
 M_u=\sum_{r\geq1}\frac{(b/a)^r}{r}
 =-\log(1-b/a)=\log(a/\beta).
\end{equation}
The root holding factor in~\eqref{eq:root-root-cycle} is $e^{-a\tau}\,d\tau$, not a separately normalized exponential density; the corresponding rate factor is already included in $\rho=a\nu$.

\paragraph{The joint transform.}

The following theorem is a finite marked version of the classical root-excursion decomposition of Le Jan~\cite[Section 7.2]{LeJan}; see also~\cite[Proposition 2.4 in the author version]{ABJL}. The calculation below is retained to make the present holding-time and cyclic-word normalization explicit.

In particular, the Gamma root occupation, its conditional Poisson excursion representation and the independent avoiding soup are classical ingredients. Le Jan's excursion theory~\cite[Section 7.2]{LeJan} gives the general framework, and A\"id\'ekon, Berestycki, Jego and Lupu~\cite[author version, Proposition 2.4]{ABJL} record the discrete root representation. Superposing the excursions there gives the aggregate excursion measure used here. The calculation below keeps the nonroot holding marks and identifies the excursion intensity with our generator convention, so that it can subsequently be inserted into the precisely weighted history law.

\begin{theorem}[Joint marked-root law]\label{thm:root-root}
Under the intensity-$\alpha$ soup on $\mathsf Q_R$, for $s\geq0$ and measurable $0\leq h\leq1$ on $\mathsf E_u$,
\begin{equation}\label{eq:root-PGFL}
 \EE\left[e^{-sT}\prod_{\epsilon\in\mathcal E}h(\epsilon)\right]
 =\left(\frac{\beta}{a+s-\rho(h)}\right)^\alpha
 =\left[1+g\left(s+\int(1-h)\,d\rho\right)\right]^{-\alpha}.
\end{equation}
Consequently the joint law has the following realization:
\begin{equation}\label{eq:root-Gamma-Poisson}
 T\sim\GammaDist(\alpha,\text{rate }\beta)
       ,
 \qquad \mathcal E\mid T=t\sim\Poi(t\rho).
\end{equation}
The entire pair $(T,\mathcal E)$ is independent of $\mathcal B$, which has the marked-soup law for $\mathsf Q_{SS}$. Thus for bounded measurable $J(\mathcal B)$, the expectation in~\eqref{eq:root-PGFL} with the extra factor $J$ is multiplied by $\EE J$. The application below takes $\alpha=1/2$.
\end{theorem}
\begin{proof}
We compute the joint transform of the root occupation and its marked excursions. This will identify the conditional excursion law, rather than only the Gamma distribution of the total occupation.

The independent trivial root clock contributes $(a/(a+s))^\alpha$. The Poisson exponential formula and Lemma~\ref{lem:root-cycle} give the contribution of the root-visiting nontrivial loops. Hence
\begin{align*}
 \EE\left[e^{-sT}\prod_{\epsilon\in\mathcal E}h(\epsilon)\right]
 &=\left(\frac{a}{a+s}\right)^\alpha
 \exp\left\{\alpha\sum_{r\geq1}\frac1r
 \left[\left(\frac{\rho(h)}{a+s}\right)^r
             -\left(\frac ba\right)^r\right]\right\}\\
 &=\left(\frac{a}{a+s}\right)^\alpha
   \left(\frac{1-b/a}{1-\rho(h)/(a+s)}\right)^\alpha
 =\left(\frac{a-b}{a+s-\rho(h)}\right)^\alpha.
\end{align*}
Both series converge absolutely because $0\leq\rho(h)/(a+s)\leq b/a<1$. The formula remains valid when $b=0$. Using $a-b=\beta=g^{-1}$ proves~\eqref{eq:root-PGFL}.

The calculation above gives the transform under the original loop soup. We now check that the proposed Gamma--Poisson construction has the same transform.

Conversely, sample $T$ with the Gamma law in~\eqref{eq:root-Gamma-Poisson} and then sample a Poisson point measure of intensity $T\rho$. Its conditional probability generating functional is
\[
 \EE\left[\prod_{\epsilon\in\mathcal E}h(\epsilon)\mid T=t\right]
 =\exp\{-t\rho(1-h)\}.
\]
The Gamma Laplace transform gives exactly~\eqref{eq:root-PGFL}. This identifies the joint law, not just the marginal of $T$: taking $h$ constant on each set of a finite measurable partition identifies the joint count generating functions together with the Laplace transform of $T$; these tests determine the law of a finite point measure on the stated standard Borel mark space. The conditional Poisson law can therefore be taken as a version for every $t>0$. Independence from $\mathcal B$ was established before the calculation by deterministic loop partition, and all transformations of the root part preserve it.
\end{proof}

The conditional Poisson law describes the aggregate excursion multiset given $T$. Its independence from $\mathcal B$ is a property of the unweighted reference law. Section~\ref{sec:comparison} will impose the low-occupation inequalities recorded during exploration through an explicit density on this joint law.

\section{Conditional comparison with independent activations}\label{sec:comparison}\label{sec:conditional-comparison}

We next compare a high star crossing with a connection in the nontrivial loop trace after independent site variables have been allowed to open incoming edges. The exploration reveals the complete marked configuration of a new high cluster, but records only an occupation inequality at a low query. The conditional law of the unexplored marked loops therefore retains the inequalities from earlier low queries. We first establish this law, then control the occupation of a pivotal vertex by keeping one of its excursions under further deletions. Finally, we use these bounds in the finite exploration.

Restriction to deterministic classes of loops follows from classical Poisson factorization; see \cite[Section 2.4]{LeJanLemaire} and \cite[Theorem 5.2]{LastPenrose}. Here we must also account for the successive observations made during the exploration.

\paragraph{Related conditioning and switching identities.}
Loop-soup explorations have several precise descriptions in the literature. Werner studies spatial resampling of oriented and unoriented soups~\cite{WernerSpatial}. A\"id\'ekon develops cluster explorations on metric graphs, with conditional excursion laws linked to boundary local times and positivity on explored edge segments~\cite[Theorems 1.2 and 1.3]{AidekonExploration}. These results explain how a revealed portion of the geometry can be used to describe the remaining random paths. The observation rule is an essential part of each such description.

Werner's switching identity gives a particularly direct description of connection conditioning for the critical loop soup on a cable graph. After fixing positive occupations at two points and conditioning them to be connected, his Theorem~2 represents the occupation field using an avoiding soup, returning excursions and excursions between the two points whose count is conditioned to be odd~\cite[Theorem 2]{WernerSwitch}. The corresponding boundary-point limit uses one connecting excursion. This gives a close precedent for understanding a connection event through excursions. An identity for the occupation field also remains valid after reweighting both sides by the same nonnegative measurable function of that field, provided its normalizing expectation is finite and positive. Consequently, upper bounds on observed occupations alone do not distinguish the present problem from such identities.

To use an excursion identity in the exploration below, we must specify both the connection event and the information retained. Our later pivotal event compares reachability with and without an incoming activation in the discrete nontrivial loop trace. The exploration also keeps the complete cyclic marked loops of every discovered high cluster. This asks for more information than the occupation field and cable-cluster connection alone. We therefore record the joint conditional law needed for these observations and prove the deletion estimate for this particular pivotal event. The formulas below make explicit which information is fixed and which loops remain random.

\subsection{Remaining networks and conditional laws}

Fix a finite network $(D,c,k,\mathsf Q_D)$ as in Section~\ref{sec:finite-networks}, an intensity $\alpha>0$, and a threshold $\lambda>0$. Write $a_v=k_v+\sum_{w\in D}c_{vw}$ for its diagonal rates. For $R\subset D$, the residual network has precision $\mathsf Q_D[R,R]$, conductances restricted to $R$, and killing
\begin{equation}\label{eq:history-killing}
 k_v^{(R)}=k_v+\sum_{w\in D\setminus R}c_{vw},\qquad v\in R.
\end{equation}
Its diagonal rate is still $a_v$. Each nonempty principal precision is positive definite. No transition probabilities are renormalized after restriction.

Let $\mu_R$ be the nontrivial cyclic marked-loop measure for this network, with the convention of Lemma~\ref{lem:finite-marked}. Write $\PP_R$ for the law of
\[
 \omega_R=\bigl(\mathcal L_R,(\xi_v)_{v\in R}\bigr),\qquad
 \mathcal L_R\sim\Poi(\alpha\mu_R),\qquad
 \xi_v\sim\GammaDist(\alpha,\text{rate }a_v),
\]
where all clocks and the nontrivial soup are independent. The total nontrivial measure $M_R=\mu_R(\text{all nontrivial loops})$ is finite. The soup is a finite point measure of original cyclic marked loops, with multiplicities and all occurrence-wise holding times. Set $\PP_\varnothing$ to be the point mass on the empty configuration.

The configuration space $\Omega_R$ is standard Borel. Indeed, each finite word and its holding coordinates form a Borel space; quotienting by the finite cyclic group preserves this property. The countable union of word spaces, its finite point-measure space, and the finite vector of positive clocks are standard Borel. Histories made from finitely many such records are also standard Borel.

Write $L_v(l)$ for the sum of the holding marks at $v$ in the loop $l$. Then $X_v(\omega_R)=\xi_v+\sum_{l\in\mathcal L_R}L_v(l)$ is full occupation. The trace graph has vertex set $R$ and an edge precisely when some nontrivial loop traverses it. Isolated vertices are singleton clusters. For a label set ${S_{\rm low}}\subset R$, define
\begin{equation}\label{eq:history-weight}
 \begin{aligned}
 W_{S_{\rm low}}^R(\omega_R)
 &=\prod_{z\in {S_{\rm low}}}\ind_{\{X_z(\omega_R)\leq\lambda\}},
 \qquad \mathcal Z_R({S_{\rm low}})=\EE_R W_{S_{\rm low}}^R,\\
 \PP_R^{S_{\rm low}}(d\omega)
 &=\frac{W_{S_{\rm low}}^R(\omega)}{\mathcal Z_R({S_{\rm low}})}\PP_R(d\omega).
 \end{aligned}
\end{equation}
The superscript records a set of retained low labels; it does not change the residual network. Empty products and $\mathcal Z_\varnothing(\varnothing)$ equal one.

\begin{lemma}[Exact restriction and positive normalization]\label{lem:history-positive}
For every $C\subset R\subset D$, the restriction of $\mu_R$ to loops contained in $C$ is $\mu_C$, including all cyclic and holding marks. Moreover,
\begin{equation}\label{eq:history-positive}
 0<e^{-\alpha M_R}
       \prod_{z\in {S_{\rm low}}}\PP\bigl(\GammaDist(\alpha,\text{rate }a_z)\leq\lambda\bigr)
 \leq \mathcal Z_R({S_{\rm low}})\leq1.
\end{equation}
Thus every kernel in~\eqref{eq:history-weight} is well-defined, including for disconnected residual networks.
\end{lemma}
\begin{proof}
In the explicit marked density~\eqref{eq:finite-marked}, a word contained in $C$ retains each conductance and each holding rate $a_v$. Its density is therefore exactly that of the principal residual network on $C$. Each old trivial clock also retains its rate. For positivity, require no nontrivial loops in $R$, an event of probability $e^{-\alpha M_R}>0$, and then require the finitely many labelled Gamma clocks to lie in $(0,\lambda]$. Each of these probabilities is positive. On this event $W_{S_{\rm low}}^R=1$, proving the bound. The empty network has the stated convention.
\end{proof}

Positivity is the only bound on the normalizer needed below.

\paragraph{Two-block factorization.}

Fix $C\subset R$ and put $E=R\setminus C$. Divide the nontrivial loop space into loops contained in $C$, loops contained in $E$, and loops visiting both. Let
\[
 m_R(C,E)=\mu_R\{l:l\cap C\ne\varnothing,\ l\cap E\ne\varnothing\},
 \qquad q_R(C,E)=e^{-\alpha m_R(C,E)}>0.
\]
The three Poisson restrictions are independent. Write $\omega_C$ and $\omega_E$ for the internal configurations, including their respective clocks, and let $A_{C,E}$ be the event that the crossing restriction is empty. On this event, all occupations on either side equal their internal occupations, and
\[
 W_{S_{\rm low}}^R=W_{{S_{\rm low}}\cap C}^C W_{{S_{\rm low}}\cap E}^E.
\]

\begin{lemma}\label{lem:history-block}
For bounded nonnegative measurable functions $f$ on $\Omega_C$ and $g$ on $\Omega_E$,
\begin{equation}\label{eq:history-block}
 \EE_R^{S_{\rm low}}\!\left[\ind_{A_{C,E}}f(\omega_C)g(\omega_E)\right]
 =\frac{q_R(C,E)}{\mathcal Z_R({S_{\rm low}})}
    \EE_C[fW_{{S_{\rm low}}\cap C}^C]\,\EE_E[gW_{{S_{\rm low}}\cap E}^E].
\end{equation}
\end{lemma}
\begin{proof}
On $A_{C,E}$, there are no loops visiting both sides of the split. This allows us to factor the weighted expectation over the two internal configurations.
Multiply the expectation under $\PP_R$ by $W_{S_{\rm low}}^R/\mathcal Z_R({S_{\rm low}})$. On $A_{C,E}$ the weight separates as above. The crossing process, both internal soups, and their two clock vectors are independent. The crossing process is empty with probability $q_R(C,E)$, while Lemma~\ref{lem:history-positive} identifies the internal laws as $\PP_C$ and $\PP_E$. Multiplying these factors proves~\eqref{eq:history-block}.
\end{proof}

When we reveal a complete trace cluster, its vertex set is random. We first apply the identity to a fixed candidate $C$ and then sum over the possible revealed clusters.

\subsection{Conditional laws after revealing complete clusters}

Fix a fresh root $u\in R\setminus {S_{\rm low}}$. Let $C_u$ be its complete trace-cluster vertex set in the original configuration. For a deterministic $C\subset R$ containing $u$, let $\mathsf{Conn}_{C,u}$ be the event that every vertex of $C$ is connected to $u$ by the internal trace graph on $C$. For a singleton $C=\{u\}$ this event is certain. Then
\begin{equation}\label{eq:history-cluster-event}
 \{C_u=C\}=A_{C,E}\cap\mathsf{Conn}_{C,u},\qquad E=R\setminus C.
\end{equation}
On this event, every loop touching $C$ is contained in $C$. Thus $\omega_C$ is exactly the marked configuration revealed with this cluster.

\begin{lemma}[Conditional law after revealing a cluster]\label{lem:history-reveal}
Under $\PP_R^{S_{\rm low}}$, define the following finite measure on internal configurations:
\begin{equation}\label{eq:history-internal-measure}
 \nu_{R,{S_{\rm low}},u,C}(dy)=
 \frac{q_R(C,E)\mathcal Z_E({S_{\rm low}}\cap E)}{\mathcal Z_R({S_{\rm low}})}
 \ind_{\mathsf{Conn}_{C,u}}(y)\ind_{\{X_u(y)>\lambda\}}
 W_{{S_{\rm low}}\cap C}^C(y)\,\PP_C(dy).
\end{equation}
The joint law on the high branch with cluster $C$ is
\begin{equation}\label{eq:history-high-kernel}
 \PP_R^{S_{\rm low}}\bigl(C_u=C,\ X_u>\lambda,
                    \omega_C\in dy,\omega_E\in dz\bigr)
 =\nu_{R,{S_{\rm low}},u,C}(dy)\,\PP_E^{{S_{\rm low}}\cap E}(dz).
\end{equation}
This is an identity of measures, equivalently an identity tested against all bounded measurable internal and exterior functions. Hence, for almost every revealed internal marked datum on this branch, the exterior law is $\PP_E^{{S_{\rm low}}\cap E}$.

On the low branch,
\begin{equation}\label{eq:history-low-kernel}
 \PP_R^{S_{\rm low}}(X_u\leq\lambda)
 =\frac{\mathcal Z_R({S_{\rm low}}\cup\{u\})}{\mathcal Z_R({S_{\rm low}})},
 \qquad
 \mathcal L(\omega_R\mid X_u\leq\lambda;\PP_R^{S_{\rm low}})
 =\PP_R^{{S_{\rm low}}\cup\{u\}}.
\end{equation}

For initialization, let $F\subset D$ be a deterministic seed set and let $K_F$ be the union of all complete trace clusters meeting $F$. For $F\subset C\subset D$, let $\mathsf{Seed}_{C,F}$ mean that every internal trace component of $C$ meets $F$. Under the unweighted law,
\begin{equation}\label{eq:history-seed-kernel}
 \PP_D\bigl(K_F=C,\omega_C\in dy,\omega_{D\setminus C}\in dz\bigr)
 =q_D(C,D\setminus C)\ind_{\mathsf{Seed}_{C,F}}(y)
    \PP_C(dy)\PP_{D\setminus C}(dz).
\end{equation}
Thus the exterior remains unweighted after the complete seed-union revelation. The seed union need not be connected.
\end{lemma}
\begin{proof}
We first fix the candidate cluster $C$, since Lemma~\ref{lem:history-block} applies to a deterministic split of the network.
For~\eqref{eq:history-high-kernel}, use~\eqref{eq:history-cluster-event} in Lemma~\ref{lem:history-block}, replacing its internal test by
$f\ind_{\mathsf{Conn}_{C,u}}\ind_{\{X_u>\lambda\}}$.
On the no-crossing event, the high condition depends only on the marked configuration inside $C$. The exterior factor satisfies
\[
 \EE_E[gW_{{S_{\rm low}}\cap E}^E]
 =\mathcal Z_E({S_{\rm low}}\cap E)\,\EE_E^{{S_{\rm low}}\cap E}g.
\]
This gives the product of $\nu_{R,{S_{\rm low}},u,C}(f)$ and $\EE_E^{{S_{\rm low}}\cap E}g$, establishing the measure identity and its conditional-law consequence. If the internal factor is zero, the branch is null and no conditional assertion on that branch is needed.

Every earlier low-occupation condition at a vertex of $C$ appears in $W_{{S_{\rm low}}\cap C}^C(y)$. It is checked against the occupation revealed on $C$ and is included in the internal branch measure before that vertex is removed from the remaining low-label set. The absence-of-crossing factor and the exterior normalizer affect which marked cluster is revealed, but neither depends on the exterior configuration after normalization. In particular they cancel in its conditional law.

The low branch is simpler: no marked cluster is revealed, so the new information is only the upper bound on occupation at $u$.
For~\eqref{eq:history-low-kernel}, multiplication by $\ind_{\{X_u\leq\lambda\}}$ replaces $W_{S_{\rm low}}^R$ by $W_{{S_{\rm low}}\cup\{u\}}^R$. Positivity is supplied by Lemma~\ref{lem:history-positive}.

For the seed union, $K_F=C$ is exactly the event of no crossing loops and $\mathsf{Seed}_{C,F}$, provided $F\subset C$. The condition requires every internal trace component to meet $F$; it does not require all vertices of $C$ to be connected. Thus several components are allowed. Apply the unweighted block identity. If $F=\varnothing$, only $C=\varnothing$ has a nonzero seed predicate. The empty exterior and isolated seed cases follow from the same conventions.
\end{proof}

As a normalization check, the low branch and all possible high-cluster branches partition the sample space. Summing their masses yields
\begin{align}\label{eq:history-mass-balance}
 \mathcal Z_R({S_{\rm low}})
 &=\mathcal Z_R({S_{\rm low}}\cup\{u\})\nonumber\\
 &\quad+\sum_{\substack{C\subset R\\u\in C}}
 q_R(C,R\setminus C)\,
 \EE_C\!\left[\ind_{\mathsf{Conn}_{C,u}}\ind_{\{X_u>\lambda\}}
                         W_{{S_{\rm low}}\cap C}^C\right]
 \mathcal Z_{R\setminus C}({S_{\rm low}}\setminus C).
\end{align}
\subsection{The exploration and its conditional law}

Fix a deterministic query domain $J\subset D$ and seed set $F\subset D$. The observation rule is as follows.
\begin{enumerate}
\item Initially reveal $K_0=K_F$ and the complete original marked configuration on $K_0$, including its separate clock values. Set ${S_{{\rm low},0}}=\varnothing$ and $R_0=D\setminus K_0$.
\item At stage $j$, choose $u_j\in J\setminus(K_j\cup {S_{{\rm low},j}})$ by a measurable function of the existing history. Any eligibility or stopping rule must also be a function of that history. A fixed deterministic ordering of eligible vertices supplies such a choice.
\item Record the bit $\eta_j=\ind_{\{X_{u_j}>\lambda\}}$. If $\eta_j=0$, reveal only this bit, put ${S_{{\rm low},j+1}}={S_{{\rm low},j}}\cup\{u_j\}$, and leave $K_j$ unchanged.
\item If $\eta_j=1$, reveal the complete trace cluster $C_j$ of $u_j$ in $D$, all its original cyclic marked loops, and the separate trivial occupations at every vertex of $C_j$. Set
\[
 K_{j+1}=K_j\cup C_j,\qquad {S_{{\rm low},j+1}}={S_{{\rm low},j}}\setminus C_j,
 \qquad R_{j+1}=R_j\setminus C_j.
\]
\end{enumerate}
Let $\mathcal H_j$ be the sigma-field generated by the initial record and all ordered query vertices, answer bits, and high-branch marked revelations completed before stage $j$. It includes every marked configuration revealed after a high answer, the identities of previously queried low vertices later included in revealed clusters, and all deterministic consequences of the record. At a retained low vertex, the only direct occupation observation is its upper bound; no exact occupation or clock value is observed there, although the complete history may imply further information.

The loop data are intrinsic cyclic marked point measures, including multiplicity and grouping. A listing computed only from an already revealed restriction gives the same information. A global enumeration that uses unrevealed atoms, the total soup count, exterior clocks, or unrevealed occupation values is not part of this observation rule. The word ``full'' refers to every datum in this specified history, not arbitrary extra information about the exterior.

Because $K_j$ is always a union of complete original trace clusters, no original loop crosses between $K_j$ and $R_j$. Every original occupation on $R_j$ therefore equals the occupation of its residual configuration. The complete cluster of a newly queried vertex is disjoint from $K_j$ and coincides with its cluster in the remaining configuration. Once a vertex has been queried it stays in $K_j\cup {S_{{\rm low},j}}$, so it is not queried again. For the conditional-law statement below, we freeze the recorded history after a stopping time determined by that history.

The next statement records the conditional-law bookkeeping for the specified observation rule. Its ingredients are Poisson restriction and factorization~\cite{LeJanLemaire,LastPenrose}, applied successively to the complete marked high-cluster revelations and the low-answer events. The spatial-resampling and metric-graph exploration precedents discussed above~\cite{WernerSpatial,AidekonExploration} explain the broader framework. We keep the finite induction to verify the actual retained data and the residual density at each history, including at the allowed stopping rules.

\begin{theorem}[Conditional law of the remaining marked configuration]\label{thm:history-source}
Under the observation rule above, for every stage $j$ and almost every realized history $h$, a regular conditional law of the residual original marked configuration is
\begin{equation}\label{eq:source}
 \mathcal L(\omega_{R_j}\mid\mathcal H_j=h)
 =\PP_{R_j(h)}^{{S_{{\rm low},j}}(h)}
 =\frac{W_{{S_{{\rm low},j}}(h)}^{R_j(h)}}{\mathcal Z_{R_j(h)}({S_{{\rm low},j}}(h))}
       \,\PP_{R_j(h)}.
\end{equation}
The law retains all cyclic marked loops and the separate aggregate trivial clock at each residual vertex. It applies with general holding rates and arbitrary $\alpha>0$, in particular with $\alpha=1/2$. It remains valid at a history-measurable stop. It is not an assertion of independence between the realized history and its residual configuration.
\end{theorem}
\begin{proof}
The argument is an induction over the recorded observations. At each stage the conditional law given the preceding history is the starting point; we then check the low and high answers separately.

We prove the conditional-expectation identity for bounded measurable test functions on the disjoint union $\bigsqcup_{R\subset D}\Omega_R$. All state and history spaces are standard Borel, so regular conditional laws exist; see \cite[Theorem A.14]{LastPenrose}. The displayed candidate kernel is measurable in the history: the finitely valued sets $R_j,{S_{{\rm low},j}}$ are history-measurable, and each normalizer is strictly positive.

\emph{Initialization.} For each deterministic candidate $C=K_0$, equation~\eqref{eq:history-seed-kernel}, tested against any bounded function of the full internal datum, gives exterior conditional law $\PP_{D\setminus C}$. Summing over the finitely many sets $C$ proves~\eqref{eq:source} for $j=0$. This includes continuous numerical revelations without assigning positive probability to any exact marked configuration.

\emph{Predictable selection.} Assume the formula at stage $j$. Condition on its actual full history $h$, outside the null set where the inductive formula fails. Then $R=R_j(h)$, ${S_{\rm low}}={S_{{\rm low},j}}(h)$, and the selected fresh root $u=u_j(h)$ are fixed. The conditional residual law is $\PP_R^{S_{\rm low}}$. Recording this predictable choice adds no likelihood factor: the root identity was already a measurable function of $h$. If the exploration stops, both the recorded history and the remaining marked configuration are left unchanged, so the formula persists.

\emph{Low answer.} The only new observation is $X_u\leq\lambda$. Formula~\eqref{eq:history-low-kernel} gives precisely $\PP_R^{{S_{\rm low}}\cup\{u\}}$, which is the asserted law for the updated $R_{j+1},{S_{{\rm low},j+1}}$.

\emph{High answer and complete revelation.} For every deterministic candidate cluster $C\subset R$ containing $u$, equation~\eqref{eq:history-high-kernel} gives the joint conditional measure of its full marked datum and the exterior. The exterior kernel is $\PP_{R\setminus C}^{{S_{\rm low}}\setminus C}$ for almost every internal datum under the branch measure. The updated history records exactly $h$, the already determined query, its high answer, and that full internal datum. The high condition and every earlier low-occupation condition at a vertex of $C$ appear as internal factors in~\eqref{eq:history-internal-measure}. Old ordered decisions and earlier numerical revelations remain in $h$ and have already been accounted for by the inductive hypothesis. Thus no additional exterior condition is introduced or dropped.

We use the joint measure identity here because the new history includes the complete marked observation, not just the cluster vertex set.

More formally, integrate~\eqref{eq:history-high-kernel} against a bounded measurable function of $(h,C,\omega_C)$ and a residual test function, then integrate over the law of $h$ and sum over $C$. This is the defining conditional-expectation identity for the enlarged history. A monotone-class argument extends the product tests to its full sigma-field. Together with the low branch, this proves the induction.

It remains to check that the same formula holds when the exploration stops according to the recorded history.

There are at most $|J|$ queries by freshness and the monotonicity of $K_j\cup {S_{{\rm low},j}}$; after a stop the record is frozen. The finitely many executed stages therefore have consistent almost-everywhere kernel versions. Alternatively, for a history-measurable stopping stage, split according to its finitely many possible values and apply the same conditional-expectation identity on each event. This proves the conditional-law formula at a history-measurable stopping time. The argument does not use the later continuation or exhaustion estimates.
\end{proof}

The conditional law of the remaining configuration therefore depends on the recorded history only through $(R_j,{S_{{\rm low},j}})$. This is what allows us to study the occupation at one new query using the root-excursion decomposition.

\subsection{A conditional bound using one root excursion}\label{subsec:pivotal-posterior}

\paragraph{Occupation control and pivotal-loop surgery.}
Gao, Nolin and Qian~\cite[Lemmas 3.15 and 4.12]{GNQ} provide close methodological precedents for combining loop occupation with connectivity. Their Lemma~3.15 controls the occupation contributed by loops avoiding a nearby vertex, uniformly in the surrounding domain. Their Lemma~4.12 separates a high occupation event from a modified arm event by a surgery argument on a pivotal loop, together with independence of disjoint loop classes. These results show concretely how local occupation estimates can enter geometric percolation arguments. The comparison below has a related purpose: it bounds occupation after the recorded exploration history and the pivotality of an incoming activation have both been fixed. Its deletion certificate is designed for that directed trace event and the retained observations. The switching identities discussed earlier give a complementary excursion-based precedent for connection conditioning.

\paragraph{The reference law and forced graphs.}

For the comparison construction, take $D\subset\mathbb Z^2$; star adjacency means $|v-w|_\infty=1$ for distinct vertices in $D$.

Fix an almost-everywhere history from Theorem~\ref{thm:history-source}. Its revealed set is $K$, its residual is $R=D\setminus K$, and its surviving low-label set is ${S_{\rm low}}\subset R$. Let $J\subset D$ be the query domain and let $\mathcal T\subset D$ be a fixed target set. In the DGFF application,
\[
 D=B_{2N},\qquad J=B_N,\qquad \mathcal T=\operatorname{BOT}_N.
\]
Choose a fresh root $u\in J\setminus(K\cup {S_{\rm low}})$ star-adjacent to $K$. Every vertex of $K$ remains seeded in the reachability graph, including its vertices outside $J$. Let $\mathsf Q_R$ be the principal precision of the current residual network, and put
\[
 g=(\mathsf Q_R^{-1})_{uu}>0,\qquad S=R\setminus\{u\}.
\]

Under the unweighted reference law of Theorem~\ref{thm:root-root}, write $T=X_u$ for full root occupation, $\mathcal E$ for the finite multiset of marked first-return root excursions, and $\mathcal B$ for the exact root-avoiding marked background on $S$. Thus
\begin{equation}\label{eq:posterior-reference}
 T\sim\GammaDist(\alpha,\text{rate }g^{-1}),\qquad
 \mathcal E\mid T=t\sim\Poi(t\rho),\qquad
 (T,\mathcal E)\ \hbox{is independent of }\mathcal B.
\end{equation}
The finite measure $\rho$ has the normalization of Lemma~\ref{lem:root-rho}; in particular $\rho(1)=a_u-g^{-1}$. Every nonroot holding occurrence is part of its mark. The occupation $T$ includes the trivial root clock.

For $v\in S$, define $L_v(\epsilon)=\sum_{i:x_i=v}\eta_i$. For an unoriented edge $e$, let $n_e(\epsilon)$ count its traversals in the full excursion word, including the first and last edges. Let $X_v^{\mathcal B}$ and $n_e^{\mathcal B}$ be the corresponding background occupation and traversal count. These are data of $\mathcal B$; its trivial clocks contribute to occupation and create no edges. For a realization $b$ of $\mathcal B$, write these background quantities as $X_v^b$ and $n_e^b$.

Adjoin independent future incoming activation coins $\eta^{\rm aux}_w$ at $w\in J\setminus(K\cup {S_{\rm low}}\cup\{u\})$. They are independent of the entire reference marked configuration and have fixed parameters after conditioning on the history; the intended application uses a common parameter $\varepsilon\in[0,1]$. Vertices in ${S_{\rm low}}$ have activation zero. There are no activated destinations in $D\setminus J$. An activated destination $w$ adds arcs $v\to w$ from its star-neighbors $v\in D$. These arcs have no occupation-threshold qualification.

For fixed $y=(b,\eta^{\rm aux})$, let $\mathcal A_i(y)$, $i\in\{0,1\}$, be the directed graph on all of $D$ formed from:
\begin{itemize}
\item the known trace edges inside $K$ and the trace edges of $b$ on $S$, each used in both directions;
\item the fixed incoming arcs into activated nonroot destinations;
\item when $i=1$, every permitted incoming arc $v\to u$ with $v\in D$ star-adjacent to $u$.
\end{itemize}
All vertices of $K$ are sources. For an excursion multiset $\nu$, let $\operatorname{tr}(\nu)$ contain both directed versions of every edge traversed by one of its excursions, and set
\begin{equation}\label{eq:posterior-graphs}
 \mathcal G_i(\nu,y)=\mathcal A_i(y)\cup\operatorname{tr}(\nu),
 \qquad
 B_i(\nu,y)=\ind_{\{K\leadsto\mathcal T\text{ in }\mathcal G_i(\nu,y)\}},
 \qquad \Delta=B_1-B_0.
\end{equation}
Graph inclusion gives $\Delta\in\{0,1\}$. All arcs in $\mathcal A_1\setminus\mathcal A_0$ have head $u$. Moreover, the root is reachable from $K$ already in $\mathcal A_1$, since a star-neighbor in $K$ supplies a direct incoming arc. Its tail need not lie in $J$.

Define
\begin{equation}\label{eq:posterior-weight-event}
 W_{S_{\rm low}}(\nu,b)=\prod_{z\in {S_{\rm low}}}
 \ind_{\{X_z^b+\sum_{e\in\nu}L_z(e)\leq\lambda\}},
 \qquad A(\nu,y)=W_{S_{\rm low}}(\nu,b)\Delta(\nu,y).
\end{equation}
This is an indicator. It has no additional dependence on the numerical root time $t$: the root is fresh, so $u\notin {S_{\rm low}}$, and graph reachability uses trace edges and independent coins. The high-root condition $T>\lambda$ is \emph{not} included in $A$; it will be used as a tail test after the posterior law is controlled. These observables are measurable by the finite path characterization of reachability and the marked occupation sums.

Lemma~\ref{lem:continuation-branches} will relate these two graphs to the actual low and high observations, including the Bernoulli variables at vertices that become part of the revealed cluster and all marked data revealed on the high branch.

\paragraph{Deterministic observables and the conditional law.}

\begin{lemma}\label{lem:root-observables}
For every original marked configuration, the projection $(T,\mathcal E,\mathcal B)$ determines its occupations and all traversal counts by
\begin{align}
 X_u&=T,&
 X_v&=X_v^{\mathcal B}+\sum_{\epsilon\in\mathcal E}L_v(\epsilon)
       &&(v\in S),\nonumber\\
 n_e&=n_e^{\mathcal B}+\sum_{\epsilon\in\mathcal E}n_e(\epsilon).\label{eq:root-observables}
 \end{align}
It therefore determines the trace graph, with an edge precisely when its count is positive, and all its vertex clusters, including isolated vertices.

Fix a retained label set ${S_{\rm low}}\subset S$, a threshold $\lambda>0$, the revealed exterior history parameters, and a realization of the independent incoming activation coins just defined. Let
\[
 W=\prod_{z\in {S_{\rm low}}}\ind_{\{X_z\leq\lambda\}},
 \qquad \Delta=B_1-B_0,
\]
where $B_i$ is the seed-to-target reachability indicator in the trace graph with the root incoming activation forced to $i$. All other activation destinations and coins are fixed, with destinations in the fixed query set $J$ and permitted tails in $D$. The whole revealed set $K$ remains present and seeded in this directed graph, with its arcs to residual destinations preserved; it is removed only from the residual soup network. Retained low vertices have activation coin zero. Then $W$, $B_0$, $B_1$, and $\Delta$ are functions of $(\mathcal E,\mathcal B)$ and those fixed parameters. There is no additional direct dependence on $T$ because the fresh root is not in ${S_{\rm low}}$ and these reachability rules do not impose occupation thresholds on the future coins.
\end{lemma}
\begin{proof}
Cutting a loop at $u$ neither deletes nor duplicates any nonroot holding occurrence or any edge traversal. Summing these contributions gives~\eqref{eq:root-observables}; the removed root holds contribute exactly to $T$ together with $\xi_0$. All excursions are attached to $u$, so their original grouping into root loops cannot change the union trace graph or its connected components. The upper-label product uses only the displayed occupations at vertices of $S$. Adding the fixed incoming arcs and computing directed reachability uses only the trace graph. Forcing the root activation changes only the prescribed incoming arcs. This proves each asserted factorization as deterministic identities, before any conditional-law calculation.
\end{proof}

The history continues to record the original cyclic marked soup. The projection above is used only for occupations and reachability. Combining it with Theorem~\ref{thm:history-source} now gives the conditional root representation.

\begin{corollary}[The root representation under the history law]\label{cor:history-root-interface}
Fix an almost-everywhere history from Theorem~\ref{thm:history-source}, with residual $R$, retained labels ${S_{\rm low}}$, and a fresh root $u\in R\setminus {S_{\rm low}}$. Let $g=(\mathsf Q_D[R,R]^{-1})_{uu}$ and let $\rho$ be the excursion measure in~\eqref{eq:root-rho}. For any bounded measurable function $\Phi$ of the reduced root data $(T,\mathcal E,\mathcal B)$,
\begin{align}\label{eq:history-root-interface}
 \EE[\Phi\mid\mathcal H]
 =\frac1{\mathcal Z_R({S_{\rm low}})}
 \int \Phi(t,e,b)W_{S_{\rm low}}(e,b)
 \,\GammaDist(\alpha,\text{rate }g^{-1})(dt)
 \,\Poi(t\rho)(de)\,\PP_{R\setminus\{u\}}(db),
\end{align}
where the final law is the exact marked root-avoiding background and
\[
 W_{S_{\rm low}}(e,b)=\prod_{z\in {S_{\rm low}}}
 \ind_{\{X_z^b+\sum_{\epsilon\in e}L_z(\epsilon)\leq\lambda\}}.
\]
Independent future activation coins can be adjoined to both sides with their product law, provided they were not used to inspect or select the source history. The weight occurs exactly once. The unweighted root/background independence and Gamma marginal are reference-law properties, not generally properties of the weighted conditional law.
\end{corollary}
\begin{proof}
Theorem~\ref{thm:history-source} gives $\EE_R[W_{S_{\rm low}}^R\Phi]/\mathcal Z_R({S_{\rm low}})$. Apply Theorem~\ref{thm:root-root} to this unweighted expectation. Lemma~\ref{lem:root-observables} identifies the label factor as $W_{S_{\rm low}}(e,b)$, since $u\notin {S_{\rm low}}$. This proves~\eqref{eq:history-root-interface}. Adjoining an independent product of coins preserves the identity. The argument keeps the original marked soup for the history and projects it only after its conditional law has been proved.
\end{proof}

\paragraph{Keeping one root excursion after deletion.}

Write $\nu'\leq\nu$ for deletion of excursion occurrences, with multiplicity retained. Holding marks of a retained occurrence are unchanged.

The point of the following certificate is the compatibility of deletion with all parts of the present event. It keeps a successful route when the root activation is on, preserves failure when it is off, and preserves every retained upper-occupation bound. This is the event-specific step that permits the standard Poisson thinning argument to be used under the history law.

\begin{lemma}[One excursion suffices after deletion]\label{lem:posterior-certificate}
For every fixed $y=(b,\eta^{\rm aux})$ and every finite excursion multiset $\nu$ with $A(\nu,y)=1$, there is a measurably chosen submultiset $\kappa(\nu,y)\leq\nu$, of size at most one, such that
\begin{equation}\label{eq:posterior-certificate}
 \kappa(\nu,y)\leq\nu'\leq\nu
 \quad\Longrightarrow\quad A(\nu',y)=1.
\end{equation}
The certificate may be empty. The conclusion holds for every intermediate subconfiguration, not only for the certificate alone.
\end{lemma}
\begin{proof}
We argue in two parts. First we keep a route to the target using at most one root excursion. We then check that deleting other excursions preserves the remaining conditions of the event.

Since $\Delta(\nu,y)=1$, the target is reachable in $\mathcal G_1(\nu,y)$ and not in $\mathcal G_0(\nu,y)$. A successful simple path from $K$ must use an arc of $\mathcal A_1\setminus\mathcal A_0$; otherwise the same path would exist in $\mathcal G_0$. Each such arc enters $u$. There is therefore a directed path from $u$ to $\mathcal T$ in $\mathcal G_1$.

Choose a simple such path $\pi=(v_0=u,v_1,\ldots,v_m\in\mathcal T)$ by a fixed ordering of the finitely many simple paths. A path of length zero is allowed if $u\in\mathcal T$. First suppose that every directed step of $\pi$ is already in $\mathcal A_1(y)$. Then $K$ reaches $u$ in $\mathcal A_1$ and follows $\pi$ to the target, even if every excursion is deleted. Choose the empty certificate.

Otherwise, let $(v_j,v_{j+1})$ be the last step of $\pi$ not contained in $\mathcal A_1(y)$. It belongs to the trace of at least one excursion occurrence. Choose one such occurrence $e_*$ and retain it as the certificate. Every edge of the suffix from $v_{j+1}$ to the target is in $\mathcal A_1$. The trace of $e_*$ is connected to $u$ and is available in both directions. It therefore contains a walk from $u$ to $v_{j+1}$. For every $\nu'$ containing $e_*$, concatenate the fixed access from $K$ to $u$, that walk in $e_*$, and the fixed suffix. This is a directed walk to the target in $\mathcal G_1(\nu',y)$, so $B_1(\nu',y)=1$. A walk suffices for reachability; it need not remain simple after concatenation.

This definition uses the last step absent from the fixed graph. An edge that is both an excursion edge and a fixed background edge or activation arc may be treated as fixed. Overlapping excursion traces therefore cause no ambiguity.

The certificate chosen above preserves reachability when the root activation is forced to one, even when the certificate is empty. It remains to check the failure of reachability when it is forced to zero and the upper bounds on occupation.

For either certificate choice, deleting excursions only removes trace edges from $\mathcal G_0$. Since $B_0(\nu,y)=0$, it follows that $B_0(\nu',y)=0$. Also, each $L_z(e)$ is nonnegative, so
\[
 X_z^b+\sum_{e\in\nu'}L_z(e)
 \leq X_z^b+\sum_{e\in\nu}L_z(e)
 \leq\lambda,\qquad z\in {S_{\rm low}}.
\]
Thus $W_{S_{\rm low}}(\nu',b)=1$, proving~\eqref{eq:posterior-certificate}.

So far we have shown that a suitable certificate exists. To apply the later thinning argument, we must also choose it measurably.

For measurability, order the finite vertex set and hence its finite list of simple paths. Membership of each path in the graph and the last nonfixed step are measurable. Order excursion words first by length and then by the vertex ordering, and order their finite vectors of holding marks lexicographically. This gives a measurable ordering on each finite excursion multiset. Choose the first occurrence carrying the selected edge. Equal marks are counted with their multiplicities; one chosen copy may be assigned an auxiliary occurrence index for thinning. Only the existence of that retained copy matters to the graph. Set $\kappa=\varnothing$ on $A=0$. The selector is part of this probability argument, and is not added to the observed exploration history.
\end{proof}

The selected excursion retains all its repeated visits and holding marks. Deletion acts on the residual excursion representation, whose conditional Poisson law is Theorem~\ref{thm:root-root}; the already revealed cyclic loops remain part of the fixed history.

\paragraph{Thinning and Gamma comparison.}

The next lemma shows how preserving an event after deleting all but a bounded number of Poisson points gives a bound on the conditional Gamma law. We include the finite marked-process proof of the thinning step; see also \cite[Section 5.3, Theorem 5.8 and Corollary 5.9]{LastPenrose}.

The distributional mechanism is classical: multiplication of a Gamma density by a power of its variable increments the shape parameter after normalization. Here the certificate makes the remaining likelihood factor decreasing, which gives stochastic domination by that size-biased Gamma law. Thus the shape increment is familiar; the work needed for the application is to verify the certificate and its compatibility with the conditioning. Poisson thinning itself is the standard operation described in~\cite[Section 5.3]{LastPenrose}.

\begin{lemma}[Poisson thinning and a Gamma bound]\label{lem:posterior-thinning}
Let $\rho$ be a finite measure on a standard Borel mark space and let $Y$ have law $\Pi$. Suppose the joint reference law is
\[
 \GammaDist(\alpha,\text{rate }g^{-1})(dt)\,
 \Poi(t\rho)(d\nu)\,\Pi(dy),\qquad \alpha,g>0.
\]
Let $A(\nu,y)\in\{0,1\}$ be measurable, with no direct $t$ dependence. Suppose that, whenever $A=1$, a measurable certificate of at most $r$ occurrences guarantees $A=1$ for every subconfiguration containing it, where $r\geq0$ is an integer. Define the kernel version
\begin{equation}\label{eq:posterior-ell-kernel}
 \ell(t)=\int A(\nu,y)\,\Poi(t\rho)(d\nu)\,\Pi(dy),\qquad t\geq0.
\end{equation}
Then $\ell$ is measurable, and for all $0<s<t$,
\begin{equation}\label{eq:posterior-thinning-r}
 \ell(s)\geq(s/t)^r\ell(t).
\end{equation}
If $m=\EE[A]>0$, the conditional law of $T$ given $A=1$ is stochastically dominated by $\GammaDist(\alpha+r,\text{rate }g^{-1})$. If $m=0$, every nonnegative weighted expectation $\EE[f(T)A]$ is zero.
\end{lemma}
\begin{proof}
We first make the event probability explicit at each Poisson intensity. We shall then compare two intensities by thinning. When the conditioning event has positive probability, the resulting monotonicity will give the Gamma bound.

Put $b=\rho(1)$. The Poisson expansion gives the explicit version
\begin{equation}\label{eq:posterior-ell-series}
 \ell(t)=e^{-bt}\sum_{n=0}^{\infty}\frac{t^n}{n!}
 \int\!\int A\!\left(\sum_{j=1}^n\delta_{e_j},y\right)
        \rho^{\otimes n}(de_1\cdots de_n)\,\Pi(dy).
\end{equation}
For $n=0$ the configuration is empty and the product measure has mass one. The $n$th coefficient lies between zero and $b^n$, so the series converges locally uniformly and defines a continuous, hence measurable, version for every $t\geq0$. It also covers $b=0$.

We now couple the Poisson laws at two intensities by deleting occurrences from the larger process.

Fix $0<s<t$ and put $p=s/t$. Given $\nu_t\sim\Poi(t\rho)$, retain each occurrence independently with probability $p$, keeping its entire mark. For any measurable $0\leq h\leq1$, the generating functional of the retained configuration is
\[
 \EE\prod_{e\in\nu_s}h(e)
 =\EE\prod_{e\in\nu_t}(1-p+ph(e))
 =\exp\{-tp\rho(1-h)\}
 =\exp\{-s\rho(1-h)\}.
\]
Thus $\nu_s\sim\Poi(s\rho)$. Use the same independent $Y$ at both times. On $A(\nu_t,Y)=1$, retaining all selected certificate occurrences has probability $p^k\geq p^r$, where $k\leq r$ is its size. An empty certificate survives with probability one. Conditional on its survival, arbitrary other deletions leave $A=1$ by hypothesis. Equal marks are thinned as separate occurrences, so selecting one copy is legitimate. The other occurrences need not be controlled: the certificate hypothesis covers every intermediate configuration that retains the selected ones. Averaging first over thinning and then over $(\nu_t,Y)$ proves~\eqref{eq:posterior-thinning-r} for the explicit version~\eqref{eq:posterior-ell-kernel}.

We have finished the Poisson comparison. The next step is to express its monotonicity as a comparison between Gamma distributions.

It follows that $h(t)=t^{-r}\ell(t)$ is nonincreasing on $(0,\infty)$. Let $f_{\gamma,g}$ be the density of $\GammaDist(\gamma,\text{rate }g^{-1})$. Then
\begin{equation}\label{eq:posterior-density-ratio}
 f_{\alpha,g}(t)\ell(t)
 =\frac{g^r\Gamma(\alpha+r)}{\Gamma(\alpha)}
       f_{\alpha+r,g}(t)h(t).
\end{equation}
When $m>0$, the conditional law is therefore a decreasing nonnegative tilt of $U\sim\GammaDist(\alpha+r,\text{rate }g^{-1})$. Its tilting function has finite positive mean, by integrating~\eqref{eq:posterior-density-ratio}.

A decreasing tilt gives relatively less weight to larger values. The following covariance identity makes that observation precise.

For a bounded nonnegative increasing $f$ and an independent copy $U'$,
\[
 2\Cov(f(U),h(U))
 =\EE[(f(U)-f(U'))(h(U)-h(U'))]\leq0.
\]
The identity is integrable since $f$ is bounded and $\EE h(U)<\infty$; alternatively first truncate the decreasing function $h$ and then take a limit. Dividing by $\EE h(U)$ gives
\[
 \EE[f(T)\mid A=1]
 =\frac{\EE[f(U)h(U)]}{\EE h(U)}
 \leq\EE f(U).
\]
Taking upper-tail indicators proves stochastic domination. If $m=0$, the nonnegative indicator $A$ vanishes almost surely, proving the remaining assertion without division.
\end{proof}

In the application, the selected excursion preserves a connection when the root activation is forced to one. Deletion cannot create a connection when the activation is forced to zero, and it preserves all upper-occupation conditions. This is exactly the hypothesis used in the thinning inequality.

\paragraph{The conditional bound for root occupation.}

For $\gamma>0$ and $x\geq0$, write
\[
 Q_\gamma(x)=\PP\bigl(\GammaDist(\gamma,\text{rate }1)>x\bigr).
\]
Expectations denoted by $\EE_{\mathcal H,\mathrm{coins}}$ use the law of the remaining marked configuration conditional on the fixed exploration history, together with independent future Bernoulli variables.

We now combine the classical root law, the history density and the one-excursion certificate. The resulting bound is stated for the full conditioning used by the finite comparison. Its Gamma law alone is not the distinguishing feature: the useful conclusion is that the bound survives this particular history and pivotal event, so it can be applied at successive queries.

\begin{theorem}[Conditional bound for root occupation]\label{thm:posterior-posterior}
For the fixed history, fresh root, target, and incoming-activation graph specified above, set
\[
 \ell(t)=\EE_{R,\mathrm{coins}}[W_{S_{\rm low}}\Delta\mid T=t]
\]
under the unweighted reference law, with the kernel version~\eqref{eq:posterior-ell-kernel}. Then
\begin{equation}\label{eq:ell}
 \ell(s)\geq\frac{s}{t}\ell(t),\qquad 0<s<t.
\end{equation}
Whenever $\EE_{\mathcal H,\mathrm{coins}}\Delta>0$, the law of $T$ under this history and the further conditioning $\Delta=1$ is stochastically dominated by $\GammaDist(\alpha+1,\text{rate }g^{-1})$. For every $a\geq0$, including the zero-mass case,
\begin{equation}\label{eq:posterior-history-tail}
 \EE_{\mathcal H,\mathrm{coins}}[\ind_{\{T>a\}}\Delta]
 \leq Q_{\alpha+1}(a/g)\,
       \EE_{\mathcal H,\mathrm{coins}}\Delta.
\end{equation}
In particular, at intensity $\alpha=1/2$ and $a=\lambda$, this is the required shape-$3/2$ bound.
\end{theorem}
\begin{proof}
We first apply the one-excursion comparison under the unweighted reference law. We then transfer the bound to the law conditioned on the fixed exploration history.

Under~\eqref{eq:posterior-reference}, the background and all future coins have a joint law $\Pi$ independent of $(T,\mathcal E)$. Equations~\eqref{eq:posterior-graphs}--\eqref{eq:posterior-weight-event} and Lemma~\ref{lem:root-observables} identify the same $W_{S_{\rm low}}\Delta$ as in the original residual marked source. Lemma~\ref{lem:posterior-certificate} gives the certificate required by Lemma~\ref{lem:posterior-thinning}, with $r=1$. This proves~\eqref{eq:ell} and the reference-law estimate
\[
 \EE_{R,\mathrm{coins}}[\ind_{\{T>a\}}W_{S_{\rm low}}\Delta]
 \leq Q_{\alpha+1}(a/g)\,
       \EE_{R,\mathrm{coins}}[W_{S_{\rm low}}\Delta].
\]

Both reference-law expectations already contain the low-occupation weight. The conditional-law identity divides them by the same positive normalizer, without applying that weight a second time.

Corollary~\ref{cor:history-root-interface} gives, with $Z=\mathcal Z_R({S_{\rm low}})>0$,
\begin{align*}
 \EE_{\mathcal H,\mathrm{coins}}[\ind_{\{T>a\}}\Delta]
 &=Z^{-1}\EE_{R,\mathrm{coins}}[\ind_{\{T>a\}}W_{S_{\rm low}}\Delta],\\
 \EE_{\mathcal H,\mathrm{coins}}\Delta
 &=Z^{-1}\EE_{R,\mathrm{coins}}[W_{S_{\rm low}}\Delta].
\end{align*}
Divide the reference estimate by the same $Z$. This proves~\eqref{eq:posterior-history-tail}, with the source weight used once. If $m=\EE_{R,\mathrm{coins}}[W_{S_{\rm low}}\Delta]>0$, further conditioning on $\Delta=1$ under the history gives exactly the reference density proportional to $W_{S_{\rm low}}\Delta$, so the full stochastic comparison also transfers. If $m=0$, both displayed history expectations are zero. No lower bound on $Z$ beyond positivity is needed.
\end{proof}

The comparison bounds the root occupation after conditioning on the observation history and on the pivotal event. The retained low-occupation conditions enter the density exactly once. The scale $g$ of the comparison Gamma law is the Green value of the current principal residual network. At intensity $1/2$, any uniform bound $g\leq G_*$ gives the sufficient upper tail $Q_{3/2}(\lambda/G_*)$. Corollary~\ref{cor:continuation-lattice} will prove $g\leq4$ in the original massive normalization and $\widehat g\leq8$ in the rate-one normalization for $m\leq1$; Section~\ref{sec:routes} will then verify the numerical tail condition.

\subsection{Finite exploration and continuation}\label{subsec:finite-continuation}

We now use the conditional root-occupation bound to control the exploration. At each stage, we study the conditional probability that the graph with the unobserved activations connects the revealed set to the target.

\paragraph{The exploration and its conditional success probability.}

Continue with the finite lattice network, intensity $\alpha>0$ and threshold $\lambda>0$ specified above. Fix a query domain $J\subset D$, a seed set $S\subset D$, a target set $\mathcal T\subset D$, and a common activation parameter $\varepsilon\in[0,1]$. The later Green estimate will specialize to $D=B_{2N}$.

Use the observation rule of Section~\ref{sec:conditional-comparison}. Initially reveal the complete original marked union $K_0$ of all trace clusters meeting $S$, and set ${S_{{\rm low},0}}=\varnothing$. At history $h$ with revealed union $K$ and retained low set ${S_{\rm low}}$, define the eligible set
\begin{equation}\label{eq:continuation-eligible}
 \mathcal U(K,{S_{\rm low}})=\{u\in J\setminus(K\cup {S_{\rm low}}):
       \text{some }v\in K\text{ is star-adjacent to }u\}.
\end{equation}
Stop if $K\cap\mathcal T\ne\varnothing$ or $\mathcal U(K,{S_{\rm low}})=\varnothing$. Otherwise choose $u$ by a fixed ordering of the eligible set, or any other measurable rule using only the existing exploration history. Query $X_u>\lambda$. On a low answer, add $u$ to ${S_{\rm low}}$. On a high answer, reveal its complete original trace cluster $C_u$ with all cyclic marks and separate trivial clocks, replace $K$ by $K\cup C_u$, and replace ${S_{\rm low}}$ by ${S_{\rm low}}\setminus C_u$. Denote by $\tau$ the number of queries before stopping and freeze the record thereafter.

After a low answer, the incoming activation at the queried vertex is set to zero, but the vertex and its trace edges remain in the unexplored marked configuration. It can still be reached through a later trace connection.

Adjoin a single family $(\eta^{\rm aux}_w)_{w\in J}$ of iid Bernoulli$(\varepsilon)$ variables, independent of the entire original marked soup and all its clocks. These auxiliary variables are not observed or used by the exploration. Their restrictions consequently remain product variables after conditioning on any of its full histories. Let $\Pi_\varepsilon$ denote their product law.

For $R=D\setminus K$, a residual configuration $\omega_R$, and a coin configuration $\eta^{\rm aux}$, define a directed graph $\mathcal G_h(\omega_R,\eta^{\rm aux})$ on all of $D$. It contains the known trace edges in $K$ and the residual trace edges in $R$, each used in both directions. In addition, it contains
\begin{equation}\label{eq:continuation-arcs}
 v\longrightarrow w
 \quad\text{if }w\in J\setminus(K\cup {S_{\rm low}}),\quad
 \eta^{\rm aux}_w=1,\quad v\in D,\quad |v-w|_\infty=1.
\end{equation}
All of $K$ is seeded. In particular, arc tails in $K\setminus J$ are retained. Destinations in ${S_{\rm low}}$ or $D\setminus J$ have activation zero; arcs entering the seeded set can be omitted. Set
\begin{equation}\label{eq:continuation-potential}
 b_h(\omega_R,\eta^{\rm aux})=\ind_{\{K\leadsto\mathcal T\text{ in }\mathcal G_h(\omega_R,\eta^{\rm aux})\}},
 \qquad
 \mathcal V(h)=\int b_h(\omega_R,\eta^{\rm aux})\,\PP_R^{S_{\rm low}}(d\omega_R)\Pi_\varepsilon(d\eta^{\rm aux}).
\end{equation}
Thus $0\leq \mathcal V(h)\leq1$. The conditional probability is measurable: the vertex sets are finite, reachability is a finite union of path events, and Theorem~\ref{thm:history-source} gives a measurable conditional law for the remaining marked configuration. Write $h_j$ for the realized ordered record generating $\mathcal H_j$ and set $\mathcal V_j=\mathcal V(h_j)$, with the frozen convention after $\tau$.

Finally, let $\mathcal G^{\mathrm{in}}(\omega_D,\eta^{\rm aux})$ contain the full original trace graph on $D$ and all incoming arcs into activated destinations of $J$, without any exploration labels. Define
\begin{equation}\label{eq:continuation-initial-law}
 \Theta^{\mathrm{in}}_{D,J,S,\mathcal T}(\varepsilon)
 =\PP_{D,\Pi_\varepsilon}
       (S\leadsto\mathcal T\text{ in }\mathcal G^{\mathrm{in}}).
\end{equation}
This is a trace graph with independent incoming site activations. It is distinct from the Lupu-augmented graph in~\eqref{eq:theta}.

\paragraph{Enlarging the set of starting vertices.}

\begin{lemma}\label{lem:continuation-absorption}
In a finite directed graph $G$, let every vertex of $H\supset K$ be reachable from $K$. Replace the seed set by $H$ and change or delete any arcs with head in $H$, retaining every arc whose head lies outside $H$. The resulting reachable set is exactly the original reachable set from $K$.

At a nonterminal exploration history, fix the entire residual source and all coins except the one at the fresh root $u$. Let $B_i$ be the reachability indicator with the root activation forced to $i\in\{0,1\}$, and set $\Delta=B_1-B_0$. These are exactly the forced graphs and increment in~\eqref{eq:posterior-graphs}. If $C$ is the complete original trace cluster of $u$, then:
\begin{enumerate}
\item $B_0$ is the continuation indicator for the low graph with labels ${S_{\rm low}}\cup\{u\}$, before imposing its additional source constraint $X_u\leq\lambda$.
\item $B_1$ is the continuation indicator after seeding $K\cup C$, retaining labels ${S_{\rm low}}\setminus C$, and forgetting all destination coins in $C$. The identity holds for every value of those internal coins, including swallowed old zero coins. All arcs from $K\cup C$ to activated residual destinations are retained.
\end{enumerate}
\end{lemma}
\begin{proof}
For the general graph assertion, take a path from $K$ in $G$ to a vertex outside $H$. Its suffix after the last visit to $H$ has all heads outside $H$, so it remains a path from the new seed set in the modified graph. Conversely, every new seed is already reachable in $G$. Along a path in the modified graph, a step with head outside $H$ is an old arc, while a step with head in $H$ reaches a vertex already reachable in $G$. Induction along the path proves the reverse inclusion. This also covers targets in $H$ and paths of length zero.

For the low branch, forcing the root activation to zero and adding $u$ to the list of zero-activation destinations produce the same graph. No vertex or trace edge is deleted.

For the high graph, eligibility supplies a fixed arc from some vertex of $K$ to $u$ when the root activation is one. Every vertex of $C$ is then reached using the bidirectional trace edges of its connected component, regardless of all other destination coins in $C$. Apply the first assertion with $H=K\cup C$. Incoming activation arcs with heads in $C$ can be discarded because their heads are now seeds. Arcs with tails in $C$ and heads at activated destinations outside $H$ stay in the graph, since their coin belongs to the outside destination.

There are no trace edges between $C$ and $R\setminus C$, since $C$ is a complete original trace component. The remaining trace graph is therefore exactly the known graph on $K\cup C$ plus the residual trace graph on $R\setminus C$. Outside $K\cup C$, the old retained labels and available coins agree with the new ones. The identity is pointwise in all the source data and coins. It applies before testing whether the root is high, and hence also on that branch.
\end{proof}

A Bernoulli variable at a destination controls only arcs entering that vertex. Once the vertex is in the starting set, this variable can be averaged out. Arcs entering other destinations are unchanged.

\paragraph{Conditional expectations for low and high observations.}

\begin{lemma}[Conditional expectations after low and high observations]\label{lem:continuation-branches}
Fix an almost-everywhere nonterminal history $h$, with $R,{S_{\rm low}},u$ as above. Put $T=X_u$ and average the nonroot coins in $\bar B_i(\omega_R)=\EE_{\Pi_\varepsilon}B_i(\omega_R,\eta^{\rm aux})$. Then
\begin{align}\label{eq:continuation-low-tower}
 \EE[\ind_{\{T\leq\lambda\}}\mathcal V_{j+1}\mid\mathcal H_j=h]
 &=\EE_R^{S_{\rm low}}[\ind_{\{T\leq\lambda\}}\bar B_0],\\
 \EE[\ind_{\{T>\lambda\}}\mathcal V_{j+1}\mid\mathcal H_j=h]
 &=\EE_R^{S_{\rm low}}[\ind_{\{T>\lambda\}}\bar B_1].\label{eq:continuation-high-tower}
\end{align}
Consequently the actual continuation drift is
\begin{equation}\label{eq:branch}
 \EE[\mathcal V_{j+1}\mid\mathcal H_j]-\mathcal V_j
 =\EE_{\mathcal H_j,\mathrm{coins}}
       [(\ind_{\{T>\lambda\}}-\varepsilon)\Delta].
\end{equation}
The histories in these identities include the entire ordered numerical marked record specified above.
\end{lemma}
\begin{proof}
We compute the next conditional success probability separately for low and high observations. A low answer adds one threshold condition, while a high answer also reveals an entire marked cluster. This difference determines which conditional law we use.

The coin at $u$ is independent of the source, of the history, and of every other coin. Hence the current potential satisfies
\begin{equation}\label{eq:continuation-current-mixture}
 \mathcal V(h)=\EE_R^{S_{\rm low}}[(1-\varepsilon)\bar B_0+\varepsilon\bar B_1].
\end{equation}
This identity also holds at $\varepsilon=0$ and $\varepsilon=1$.

On the low branch the next history records only the threshold answer and its predictable root. Lemma~\ref{lem:continuation-absorption} identifies its graph with $B_0$. By the low kernel~\eqref{eq:history-low-kernel}, the branch probability is
\[
 p_0(h)=\frac{\mathcal Z_R({S_{\rm low}}\cup\{u\})}{\mathcal Z_R({S_{\rm low}})},
\]
and the next source law is $\PP_R^{{S_{\rm low}}\cup\{u\}}$. Therefore
\[
 p_0(h)\mathcal V(h^{\mathrm{low}})
 =p_0(h)\EE_R^{{S_{\rm low}}\cup\{u\}}\bar B_0
 =\EE_R^{S_{\rm low}}[\ind_{\{T\leq\lambda\}}\bar B_0],
\]
which proves~\eqref{eq:continuation-low-tower}.

The high observation records the complete numerical marks of the revealed cluster, not just the threshold answer. We therefore use the joint high-branch measure rather than conditioning on a particular value of those marks.

For the high branch, fix a deterministic candidate $C\subset R$ containing $u$, put $E=R\setminus C$, and write $y$ for the complete internal marked datum, including every numerical holding occurrence and trivial clock. The joint measure~\eqref{eq:history-high-kernel} is
\[
 \nu_{R,{S_{\rm low}},u,C}(dy)\,\PP_E^{{S_{\rm low}}\setminus C}(dz).
\]
Let $h^+(C,y)$ be the updated full history. Lemma~\ref{lem:continuation-absorption} gives, on this branch, a pointwise identity between $B_1(y\oplus z,\eta^{\rm aux})$ and the reachability indicator at $h^+(C,y)$ using $z$ and the surviving outside coins. Here $y\oplus z$ denotes the original marked union on the no-crossing event. The indicator is independent of every destination coin absorbed into $C$.

The sets $C$ and the marked data $y$ are functions of the source and the previous source history, never of the auxiliary coins. For each deterministic $C,y$, integrating the discarded coin coordinates gives one, and the outside coin coordinates retain their original product law. It follows that
\begin{align}\label{eq:continuation-high-integral}
 \EE_R^{S_{\rm low}}[\ind_{\{T>\lambda\}}\bar B_1]
 &=\sum_{\substack{C\subset R\\u\in C}}
   \int\nu_{R,{S_{\rm low}},u,C}(dy)
   \int b_{h^+(C,y)}(z,\eta^{\rm aux})\,
          \PP_E^{{S_{\rm low}}\setminus C}(dz)\Pi_\varepsilon(d\eta^{\rm aux})\nonumber\\
 &=\sum_{\substack{C\subset R\\u\in C}}
   \int \mathcal V(h^+(C,y))\,\nu_{R,{S_{\rm low}},u,C}(dy).
\end{align}
This is the left side of~\eqref{eq:continuation-high-tower}. In particular, the no-crossing factor, all earlier low-occupation conditions at vertices of $C$, and the exterior normalizer remain inside the branch measure $\nu$; none is dropped before the integration. Multiplying the integrands by any bounded measurable function of the complete $(h,C,y)$ gives the corresponding conditional identity on the full enlarged history. Thus the argument treats numerical revelations by identities of measures, without conditioning on a positive-probability numerical atom.

Adding the two branch identities gives
\[
 \EE[\mathcal V_{j+1}\mid\mathcal H_j=h]
 =\EE_R^{S_{\rm low}}[\bar B_0+\ind_{\{T>\lambda\}}(\bar B_1-\bar B_0)].
\]
Subtract~\eqref{eq:continuation-current-mixture}. The result is~\eqref{eq:branch}, with exactly the pivotal increment of~\eqref{eq:posterior-graphs} and its independent nonroot coins. Theorem~\ref{thm:history-source} supplies the common conditional source law throughout. There is no additional posterior weight introduced in this calculation.
\end{proof}

\paragraph{Initialization and finite stopping.}

Let $\mathsf{High}_{J,S,\mathcal T}(\lambda)$ denote the event that there is a star path inside $J$, starting in $S\cap J$ and ending in $\mathcal T\cap J$, all of whose vertices satisfy $X_v>\lambda$. Let $\mathsf{Succ}=\{K_\tau\cap\mathcal T\ne\varnothing\}$.

\begin{theorem}[Finite incoming-activation comparison]\label{thm:continuation-finite}
Under the exploration above, suppose that, at every executed stage and almost every nonterminal history,
\begin{equation}\label{eq:continuation-history-condition}
 \varepsilon\geq Q_{\alpha+1}(\lambda/g_h),
 \qquad g_h=(\mathsf Q_D[R(h),R(h)]^{-1})_{u(h),u(h)}.
\end{equation}
The parameter $\varepsilon$ is the same at every site and stage. Then $(\mathcal V_{j\wedge\tau})_{j\geq0}$ is a bounded supermartingale for the specified history, $\tau\leq|J|$, and
\begin{equation}\label{eq:continuation-comparison}
 \PP(\mathsf{High}_{J,S,\mathcal T}(\lambda))
 \leq\PP(\mathsf{Succ})
 =\EE \mathcal V_\tau
 \leq\EE \mathcal V_0
 =\Theta^{\mathrm{in}}_{D,J,S,\mathcal T}(\varepsilon).
\end{equation}
These probabilities are for the original marked source, and the rightmost graph uses iid incoming activations independent of that source.
\end{theorem}
\begin{proof}
We first prove that the conditional success probability cannot increase on average. We then use finite termination and the initial and terminal values to obtain the crossing comparison.

\emph{Conditional drift.} At an executed stage, Lemma~\ref{lem:continuation-branches} and Theorem~\ref{thm:posterior-posterior} give
\begin{align}\label{eq:continuation-drift-bound}
 \EE[\mathcal V_{j+1}\mid\mathcal H_j]-\mathcal V_j
 &\leq\bigl[Q_{\alpha+1}(\lambda/g_{\mathcal H_j})-\varepsilon\bigr]
          \EE_{\mathcal H_j,\mathrm{coins}}\Delta
 \leq0.
\end{align}
Only configurations on which the root activation changes reachability contribute to the drift. The conditional root-occupation bound controls precisely these configurations.

If the pivotal mass is zero, both terms involving $\Delta$ vanish. After stopping the process is constant. The null sets arising from the finitely many conditional-law identities can be combined into one, so these are valid conditional inequalities throughout the finite exploration.

\emph{Finite termination.} Initially $K_0$ is a union of complete original trace clusters. A high step adds another complete cluster; a low step leaves that union unchanged. Thus no trace edge crosses between $K_j$ and its complement at any stage. Every retained label records a true inequality $X_z\leq\lambda$ for the original full occupation. The sets $K_j\cup {S_{{\rm low},j}}$ are increasing, and every queried root is a new member of that union. Since all roots lie in $J$, at most $|J|$ queries are made. If the record has not already stopped, after that many queries there is no fresh eligible root, so $\tau\leq|J|$.

\emph{Terminal values.} If $K_\tau$ meets $\mathcal T$, the seeded target gives $\mathcal V_\tau=1$. Otherwise the stop is genuine exhaustion: $\mathcal U(K_\tau,{S_{{\rm low},\tau}})=\varnothing$. For every residual source configuration in the conditional kernel, no trace edge exits $K_\tau$. An outgoing activation arc from $K_\tau$ to its complement would have to end at a fresh activated vertex of $J\setminus(K_\tau\cup {S_{{\rm low},\tau}})$ star-adjacent to $K_\tau$. Such a vertex would be eligible, contrary to exhaustion. Old low destinations have activation zero, and destinations outside $J$ are not activated. Hence no directed path leaves the seeded set, so $\mathcal V_\tau=0$. We have proved the exact identity $\mathcal V_\tau=\ind_{\mathsf{Succ}}$.

\emph{Initial potential.} For each complete original configuration, every vertex of $K_0$ is reachable from $S$ by trace edges. Lemma~\ref{lem:continuation-absorption} therefore permits seeding all of $K_0$ and deleting its incoming destination arcs, without changing reachability in $\mathcal G^{\mathrm{in}}$. This is a pointwise identity for every coin configuration. Conditional on the full initial marked record, equation~\eqref{eq:history-seed-kernel} gives exterior law $\PP_{D\setminus K_0}$, and the unused coins retain their product law. Thus
\begin{equation}\label{eq:continuation-initial-tower}
 \mathcal V_0=\EE\!\left[
 \ind_{\{S\leadsto\mathcal T\text{ in }\mathcal G^{\mathrm{in}}\}}
 \mid\mathcal H_0\right].
\end{equation}
Averaging proves the last equality in~\eqref{eq:continuation-comparison}. The seed union may have several components; each meets $S$, which is exactly what this argument requires.

\emph{High-path inclusion.} Suppose a high star path as defined above exists and the exploration terminates without success. Its initial vertex is in $S\subset K_0\subset K_\tau$, while its final vertex is not in $K_\tau$. Take the first vertex of the path outside $K_\tau$. Its preceding vertex lies in $K_\tau$, it lies in $J$, and it cannot belong to ${S_{{\rm low},\tau}}$ because it is high. It is therefore an eligible fresh root, contradicting exhaustion. Thus $\mathsf{High}_{J,S,\mathcal T}(\lambda)\subseteq\mathsf{Succ}$. No assertion of equality is needed: full-domain trace paths can create additional successes.

Finally freeze $\mathcal V_j$ at $\tau$ and iterate the conditional drift inequality up to the deterministic time $|J|$. This gives $\EE \mathcal V_\tau\leq\EE \mathcal V_0$ by a finite tower calculation. Combining it with the two terminal values, initialization, and high-path inclusion proves~\eqref{eq:continuation-comparison}. No infinite-time stopping theorem or uniform-integrability argument is required.
\end{proof}

\subsection{A Green function bound at the next queried vertex}

To apply~\eqref{eq:continuation-history-condition} at every query, we need a uniform bound on the Green function at the queried vertex. Star adjacency to the revealed set supplies the short lattice routes used in the next proof. The energy principle is classical; compare \cite[Section 2.4, Exercise 2.13]{LyonsPeres}. We give the finite-matrix proof to fix the occupation normalization and the constants.

For comparison, the avoiding-loop occupation estimate in Gao--Nolin--Qian~\cite[Lemma 3.15]{GNQ} also uses a nearby excluded vertex to obtain domain-uniform local control. Here the needed Green bound is obtained directly from the finite-network energy principle, with the conductances and residual network kept explicit. This is a related local-control mechanism, with the constants computed for the present normalization.

\begin{lemma}[Path energy and the square-domain bound]\label{lem:continuation-green}
Let $K\subset D$, $R=D\setminus K$, and $u\in R$. In a finite precision network with nonnegative killing, suppose a path from $u$ to $K$ has positive edge conductances $c_1,\ldots,c_l$. Then
\begin{equation}\label{eq:continuation-path-green}
 (\mathsf Q_D[R,R]^{-1})_{uu}\leq\sum_{i=1}^l c_i^{-1}.
\end{equation}
Suppose more specifically that $D=B_{2N}$, every nearest-neighbor edge in $D$ has conductance at least $c_{\min}>0$, and $u\in B_N$ is star-adjacent to some vertex of $K$. Then the sharper bound is
\begin{equation}\label{eq:continuation-box-green}
(\mathsf Q_D[R,R]^{-1})_{uu}\leq c_{\min}^{-1}.
\end{equation}
No vertices of the low-label set ${S_{\rm low}}$ are removed in these inequalities.
\end{lemma}
\begin{proof}
We control the Green value through the energy of a function that vanishes on the revealed set $K$. The change from its value at $u$ to zero along a path into $K$ is controlled by the energy and the variational formula.

For a vector $f$ on $R$, extend it by zero on $K$, writing $\widetilde f$ for the extension. With each unordered edge counted once,
\begin{equation}\label{eq:continuation-energy}
 f^{\mathsf T}\mathsf Q_D[R,R]f
 =\sum_{\{v,w\}\subset D}c_{vw}(\widetilde f(v)-\widetilde f(w))^2
    +\sum_{v\in D}k_v\widetilde f(v)^2.
\end{equation}
The sum includes all edges into $K$ and all nonnegative killing. For any positive definite matrix $A$, Cauchy--Schwarz in its inner product gives
\begin{equation}\label{eq:continuation-variational}
 (A^{-1})_{uu}
 =\sup_{f\ne0}\frac{f(u)^2}{f^{\mathsf T}Af}.
\end{equation}
Indeed, $f(u)=\langle A^{-1}e_u,f\rangle_A$ gives the upper bound, and equality holds at $f=A^{-1}e_u$.

Truncate the stated path at its first vertex in $K$ and, if needed, erase loops. Its endpoint value is zero. Weighted Cauchy--Schwarz along the path gives
\[
 f(u)^2
 \leq\left(\sum_{i=1}^l c_i^{-1}\right)
       \sum_{i=1}^l c_i(\widetilde f(x_{i-1})-\widetilde f(x_i))^2
 \leq\left(\sum_{i=1}^l c_i^{-1}\right)f^{\mathsf T}\mathsf Q_D[R,R]f.
\]
Using~\eqref{eq:continuation-variational} proves~\eqref{eq:continuation-path-green}.

For the square-domain improvement, choose a star-neighbor $k\in K$ of $u$. If $u$ and $k$ are nearest neighbors, their edge alone contributes at least $c_{\min}f(u)^2$ to~\eqref{eq:continuation-energy}.

If they are diagonal neighbors, let $z_1,z_2$ be the two distinct nearest-neighbor intermediate vertices. Both are in the square $D$, including when $N=1$ or $k\notin B_N$. The two length-two routes use four distinct unordered edges. Their contribution to the energy is at least
\begin{align*}
 c_{\min}\bigl[(f(u)-\widetilde f(z_1))^2+\widetilde f(z_1)^2
         +(f(u)-\widetilde f(z_2))^2+\widetilde f(z_2)^2\bigr]
 \geq c_{\min}f(u)^2,
\end{align*}
because $(a-b)^2+b^2\geq a^2/2$ for each route. If an intermediate vertex already lies in $K$, its extended value is zero and the inequality still holds. All other energy terms are nonnegative. Formula~\eqref{eq:continuation-variational} now gives~\eqref{eq:continuation-box-green}.
\end{proof}

\begin{corollary}[Explicit lattice comparison and massive normalization]\label{cor:continuation-lattice}
Take $D=B_{2N}$, $J=B_N$, $S=\operatorname{TOP}_N$, $\mathcal T=\operatorname{BOT}_N$, and $\alpha=1/2$. For the original precision $I-P_D$, or the unnormalized massive precision $I-P_D+m^2I$ with $m\geq0$, every executed root has
\begin{equation}\label{eq:continuation-green-4}
 g_h\leq4.
\end{equation}
Hence the finite comparison~\eqref{eq:continuation-comparison} holds for every $N\geq1$ whenever
\begin{equation}\label{eq:continuation-epsilon-4}
 \lambda>0,\qquad Q_{3/2}(\lambda/4)\leq\varepsilon\leq1.
\end{equation}
For the rate-one normalized massive network with $a=1+m^2$ and precision
\[
 \widehat{\mathsf Q}_D=\frac{I-P_D+m^2I}{a},
\]
the residual scale satisfies $\widehat g_h=a g_h\leq4a$. If $0\leq m\leq1$, then $\widehat g_h\leq8$, and the sufficient condition for a normalized occupation threshold $\widehat\lambda$ is
\begin{equation}\label{eq:continuation-epsilon-8}
 \widehat\lambda>0,\qquad Q_{3/2}(\widehat\lambda/8)\leq\varepsilon\leq1.
\end{equation}
\end{corollary}
\begin{proof}
In the original and unnormalized massive networks, every nearest-neighbor conductance is $1/4$ and all killing is nonnegative. Apply Lemma~\ref{lem:continuation-green} with $c_{\min}=1/4$. This is a deterministic statement for every revealed union and every eligible root, so it holds uniformly over all legal histories. Since $Q_{3/2}$ is decreasing and $g_h\leq4$, condition~\eqref{eq:continuation-epsilon-4} implies~\eqref{eq:continuation-history-condition}. Theorem~\ref{thm:continuation-finite} proves the comparison.

The normalized conductances are $1/(4a)$, and its principal residual precision is exactly the unnormalized principal residual divided by $a$. Its inverse is multiplied by $a$. Thus $\widehat g_h\leq4a\leq8$ for $m\leq1$. Applying the same theorem at the normalized threshold proves the last assertion. Under the time rescaling of Section~\ref{sec:finite-networks}, occupations satisfy $\widehat X=aX$, so a corresponding threshold is $\widehat\lambda=a\lambda$ and the ratio $\widehat\lambda/\widehat g_h=\lambda/g_h$ is unchanged. The two sufficient uniform bounds are stated separately to avoid mixing their scales.
\end{proof}

\section{Independent routes and the massive crossing comparison}\label{sec:routes}\label{sec:route-construction}

Theorem~\ref{thm:continuation-finite} gives a comparison involving the nontrivial loop trace and independent Bernoulli variables at the destination vertices. We now construct these variables inside Lupu's augmented graph. For each destination, we use separate pieces of the Gamma occupation clocks and Poisson edge-opening processes to open short routes in a small block. The resulting Bernoulli variables are independent of the complete nontrivial marked soup, but not of the trivial clocks used to construct them. We therefore substitute this construction only after integrating the conditional comparison in Section~\ref{sec:comparison}. The augmented graph retains its exact law, and planar crossing duality then gives the massive crossing estimate.

\subsection{The normalized graph and its law}

Throughout this section let $N\geq1$, $D=B_{2N}$, $J=B_N$, and $0<m\leq1$. Put
\[
 a=1+m^2,\qquad \mathsf Q_m=I-P_D+m^2I,\qquad
 \widehat{\mathsf Q}=\mathsf Q_m/a=I-P_D/a.
\]
The normalized field $\chi=\sqrt a\,\psi_{N,m}$ has covariance $\widehat{\mathsf Q}^{-1}$. Every nearest-neighbor conductance is $c_e=1/(4a)\geq1/8$, every diagonal rate is one, and killing is nonnegative. Let $\mathcal L^{\mathrm{nt}}$ be the complete normalized nontrivial marked soup at intensity $1/2$, and let $Y_v$ be its occupation at $v$. Its independent trivial occupations have law $\GammaDist(1/2,\mathrm{rate}\ 1)$. The full occupation is $X_v=Y_v+\xi_v$.

Write $E_{\mathrm{nn}}(D)$ and $E_*(D)$ for the unordered nearest-neighbor and star edges with both endpoints in $D$. In particular, an edge belongs to $E_*(D)$ when its distinct endpoints have $\ell^\infty$ distance one. The assignments of independent variables and processes to vertices, edges, and blocks are fixed before sampling.

Independently of the entire $\mathcal L^{\mathrm{nt}}$, sample mutually independent variables and processes
\[
 \xi_{v,k}\sim\GammaDist(1/128,\mathrm{rate}\ 1),
 \quad v\in D,\ 1\leq k\leq64,
\]
and unit-rate Poisson point processes $\Pi_{e,k}$ on $[0,\infty)$, for $e\in E_{\mathrm{nn}}(D)$ and $1\leq k\leq64$. All clock pieces and all Poisson processes are independent of one another. Set
\begin{equation}\label{eq:routes-primitive-augmentation}
 \xi_v=\sum_{k=1}^{64}\xi_{v,k},\qquad X_v=Y_v+\xi_v,\qquad
 J_e=2c_e\sqrt{X_xX_y}\quad(e=\{x,y\}),
\end{equation}
and define
\begin{equation}\label{eq:routes-added-edge}
 A_e=\ind_{\{\sum_{k=1}^{64}\Pi_{e,k}([0,J_e/64])\geq1\}},\qquad
 \eta_e=\ind_{\{n_e>0\}}\vee A_e,
\end{equation}
where $n_e$ is the traversal count of $e$ in the nontrivial soup.

The next verification uses Gamma additivity and Poisson superposition to realize the augmentation law already supplied by Lupu~\cite[Theorem 1 bis and Corollary 3.6]{Lupu}. Splitting the primitive variables is a way to preserve that exact marginal while reserving disjoint randomness for different destinations. The identity of the full graph law is therefore a consistency check on the construction, not a separate loop-soup isomorphism.

\begin{lemma}[Exact full graph marginal]\label{lem:routes-marginal}
The marked source $(\mathcal L^{\mathrm{nt}},(\xi_v)_{v\in D})$ is the normalized intensity-$1/2$ source of Section~\ref{sec:soups}. Conditional on that entire source, the variables $A_e$ are independent and
\[
 \PP(A_e=1\mid\mathcal L^{\mathrm{nt}},(\xi_v))=1-e^{-J_e}.
\]
Consequently $\eta$ has exactly the graph law $\mathsf{Aug}(\widehat{\mathsf Q})$ from Section~\ref{sec:finite-networks}, and $X\overset d=\chi^2/2$.
\end{lemma}
\begin{proof}
At each vertex, the Laplace transform of the sum of its sixty-four pieces is
\[
 \EE e^{-t\xi_v}=(1+t)^{-64/128}=(1+t)^{-1/2},\qquad t\geq0.
\]
Thus these sums are independent $\GammaDist(1/2,\text{rate }1)$ variables, independent of the complete nontrivial marked soup. This is precisely the original trivial-occupation law in the normalized network.

Condition even on the finer information consisting of $\mathcal L^{\mathrm{nt}}$ and every individual $\xi_{v,k}$. Each $J_e$ is then fixed. The sum of the sixty-four independent counts in~\eqref{eq:routes-added-edge} has Poisson law of mean $J_e$, and sums for different edges are independent. Their product conditional probabilities depend on the clock pieces only through $(\xi_v)$, so the same product law holds given the full original source. Union with the traversed edges is the exact augmentation rule of~\eqref{eq:finite-augmentation}, supplied by \cite[Theorem 1 bis and Corollary 3.6]{Lupu}. The occupation and graph normalization is that of~\eqref{eq:finite-augmentation}.
\end{proof}

The full graph law uses all sixty-four Gamma variables and Poisson processes in~\eqref{eq:routes-primitive-augmentation}--\eqref{eq:routes-added-edge}, including those not assigned to any block. Subsequent extensions add auxiliary variables while keeping $X$ and $\eta$ fixed.

\subsection{Independent activations from small blocks}

\paragraph{Independent bonds and independent destination variables.}
Le Jan and Lemaire~\cite[Lemma 4.1]{LeJanLemaire} give an earlier and useful way to extract independent percolation variables from a soup: loops whose primitive support consists of two adjacent vertices produce independent Bernoulli bonds. Their treatment also relates loop clusters to coalescence and percolation, and Proposition~4.3 studies the effect of killing, including absence of percolation at sufficiently large killing for fixed intensity on $\mathbb Z^d$, $d\geq2$. The bond construction is a clear precedent for obtaining local product structure from Poisson loop classes. Its bond indicators are functions of the nontrivial soup itself, however. Independence across those bonds is different from independence of an additional family from the \emph{whole} nontrivial marked soup, which is the joint-law requirement here.

Chang and Sapozhnikov~\cite{ChangSapozhnikov} develop the phase-transition picture for loop percolation further. In particular, they prove a nontrivial massless intensity transition on $\mathbb Z^d$ for $d\geq3$ and analyze one-arm and cluster-size behavior. Their polynomial lower bounds in the massless subcritical regime show why a generic appeal to subcriticality is not an exponential-decay argument for loop models. This phase context concerns raw loop clusters; the present construction concerns a finite augmented graph and an independence property needed for the crossing comparison. Also, sufficiently large fixed killing does not cover the whole range $m=K/N$, where the lattice killing can tend to zero.

A different representation is supplied by Lupu and Werner~\cite{LupuWernerFK}. Their current-plus-Bernoulli coupling combines the support of an Ising random current with an independent edge configuration to give the associated Fortuin--Kasteleyn (FK) random-cluster law. They relate this identity to loop soups and to the conditional law of GFF signs given the field amplitudes. This gives an additional explanation of the augmentation's connection with sign geometry. The variables constructed below serve a separate purpose: they activate all permitted incoming routes to a destination and have the product law required by the integrated comparison.

Here the required independence comes from the allocation of disjoint trivial-clock pieces and auxiliary Poisson processes to destinations, even when their spatial blocks overlap. These primitives are independent of the entire nontrivial marked soup. The resulting activations can depend on the trivial occupations, so the construction is used only after the conditional comparison has been averaged. The next lemma checks both that independence and the pathwise inclusion of the resulting routes in the exact augmented graph.

Fix a bijection
$\kappa:\{-1,0,1\}^2\to\{1,\ldots,9\}$, once and for all. For
$v\in J$ define the $3\times3$ block $Q_v=v+\{-1,0,1\}^2$ and let
$F_v$ be its twelve internal nearest-neighbor edges. Every such block lies
in $D$, including when $N=1$. For an edge $e=\{x,y\}$ define its anchor
$b(e)=(\min(x_1,y_1),\min(x_2,y_2))$. Recall the independent primitive
variables in~\eqref{eq:routes-primitive-augmentation}.

\begin{lemma}[Independent destination blocks inside the exact graph]\label{lem:routes-blocks}
Put $\delta=1/256$ and
\[
 q_*=\PP\{\GammaDist(1/128,\text{rate }1)\geq1\},\qquad
 p_*=1-e^{-\delta},\qquad \varepsilon_0=2^{-440}.
\]
For $v\in J$, let
\begin{equation}\label{eq:routes-block-event}
 D_v=\bigcap_{w\in Q_v}\{\xi_{w,\kappa(w-v)}\geq1\}
 \ \cap\!\!\bigcap_{e\in F_v}
   \{\Pi_{e,\kappa(b(e)-v)}([0,\delta])\geq1\}.
\end{equation}
The events $(D_v)_{v\in J}$ are mutually independent, are independent of the
complete nontrivial marked soup $\mathcal L^{\mathrm{nt}}$, and each has
probability
\[
 s_*=q_*^9p_*^{12}>2^{-216}>\varepsilon_0.
\]
On $D_v$, every edge of $F_v$ is open in the exact augmented graph $\eta$.
Consequently, on an extension of the primitive probability space there
are iid $\eta_v^{\rm res}\sim\operatorname{Ber}(\varepsilon_0)$, $v\in J$,
independent of $\mathcal L^{\mathrm{nt}}$, such that whenever
$\eta_v^{\rm res}=1$, every incoming star arc $w\to v$ with
$w\in D$ and $|w-v|_\infty=1$ lifts to an open nearest-neighbor route
within $Q_v$.
\end{lemma}
\begin{proof}
We verify independence, the block-success probability, and the open routes
forced by success, before thinning to the required Bernoulli parameter.
For a fixed lattice vertex $w$, different destinations using it have
different offsets $w-v$, and therefore use different Gamma-layer labels.
For a fixed nearest-neighbor edge $e$, different destinations using it have
different anchors $b(e)-v$, and therefore use different Poisson-layer labels.
These offsets lie in $\{-1,0,1\}^2$ whenever the corresponding Gamma variable or Poisson process appears in~\eqref{eq:routes-block-event}. Thus the $21$ independent variables and processes used by one block are distinct from those used by every other block. This proves the product law, including independence from
$\mathcal L^{\mathrm{nt}}$.

We next bound the probability of this event to verify the required activation parameter.
Writing $s=1/128$, we have
$\Gamma(s)<128+e^{-1}<129$; integrating the Gamma density over $[1,2]$
gives
\[
 q_*>\frac1{18\cdot129}>2^{-12},\qquad
 p_*\geq\frac{\delta}{1+\delta}=\frac1{257}>2^{-9}.
\]
(The Gamma estimate uses $x^{s-1}\geq1/2$ and $e^{-x}>1/9$ on $[1,2]$.)
There are nine vertices and twelve edges per block, so
$s_*=q_*^9p_*^{12}>2^{-108-108}=2^{-216}$.

On $D_v$, the occupation $X_w$ is at least one at each $w\in Q_v$.
For $e=\{x,y\}\in F_v$ we have $c_e\geq1/8$, whence
\[
 J_e/64=2c_e\sqrt{X_xX_y}/64\geq1/256=\delta.
\]
The Poisson process assigned to this edge and destination has a point in $[0,\delta]$, so
$A_e=1$ in the \emph{original} exact augmentation rule.
Every star neighbor of $v$ has a nearest-neighbor route of one or two
edges to $v$ inside $Q_v$.

Since $s_*>\varepsilon_0$, thinning is possible and does not change the
open-route implication.

Finally take an independent uniform variable $V_v$ for each destination
and set
\begin{equation}\label{eq:routes-site-thinning}
 \eta_v^{\rm res}=\ind_{D_v}\ind_{\{V_v\leq\varepsilon_0/s_*\}}.
\end{equation}
Then the site variables are iid Bernoulli$(\varepsilon_0)$ and independent
of $\mathcal L^{\mathrm{nt}}$, and activity implies the claimed open routes.
The site field need not be independent of the trivial clock pieces; such
independence is neither asserted nor used.
\end{proof}

In particular the construction proves the joint-law identity
\begin{equation}\label{eq:routes-joint-law}
 \mathcal L(\mathcal L^{\mathrm{nt}},\eta^{\rm res})
 =\mathcal L(\mathcal L^{\mathrm{nt}})
   \otimes\operatorname{Ber}(\varepsilon_0)^J.
\end{equation}
This is exactly the product law of the \emph{integrated} comparison in
Section~\ref{sec:comparison}; the site coins are not substituted into
its conditional exploration before averaging the numerical histories.

\subsection{From the incoming graph to massive crossings}

\paragraph{Embedding the integrated graph.}

The following embedding is the direct geometric consequence of the independent block construction. Each added incoming arc is replaced by an open nearest-neighbor route, while original trace edges are retained. Its use depends on the product law with the complete nontrivial soup established above; it is not an appeal to a general enhancement or Bernoulli-domination theorem.

\begin{theorem}[Independent incoming graph inside the exact augmentation]\label{thm:routes-enhancement}
For the normalized network above, all deterministic seed and target sets $S,\mathcal T\subset D$ satisfy
\begin{equation}\label{eq:routes-enhancement}
 \Theta^{\mathrm{in}}_{D,J,S,\mathcal T}(\varepsilon_0)
 \leq\PP_{\mathsf{Aug}(\widehat{\mathsf Q})}(S\longleftrightarrow\mathcal T).
\end{equation}
The left side is exactly the integrated original-trace-plus-iid-incoming-site probability defined in~\eqref{eq:continuation-initial-law}. Paths on either side may travel throughout $D$.
\end{theorem}
\begin{proof}
Use the block construction on the same nontrivial marked soup as $\eta$. By~\eqref{eq:routes-joint-law}, its trace together with $\eta^{\rm res}$ has exactly the joint law used on the left side. Every original trace edge is open in $\eta$. For each additional directed arc $w\to v$ with active destination $v\in J$, Lemma~\ref{lem:routes-blocks} provides an open nearest-neighbor path in $Q_v$ from $w$ to $v$.

Replace each directed arc of a successful path by the block route supplied above, and retain the traversed trace edges. This gives an open walk in $\eta$ with the same endpoints. Erasing repeated portions yields a connection. The implication holds pointwise and simultaneously for all seed and target sets. Lemma~\ref{lem:routes-marginal} identifies the right-hand graph law.
\end{proof}

Apply this theorem after Theorem~\ref{thm:continuation-finite} has averaged its initial potential. Throughout the adaptive argument, its auxiliary coins are independent of the entire source and all clocks. The integrated event then uses only the nontrivial trace and those coins, so~\eqref{eq:routes-joint-law} permits its realization by the block variables. Arcs with tails outside $J$ are retained and routed in the same way as the other incoming arcs.

\paragraph{The numerical tail condition.}

\begin{lemma}[A uniform sufficient tail bound]\label{lem:routes-tail}
Let $\lambda_*=32768$. Then
\begin{equation}\label{eq:routes-tail}
 Q_{3/2}(\lambda_*/8)
 \leq Q_{3/2}(512)<2^{-506}<\varepsilon_0.
\end{equation}
Thus the common activation parameter $\varepsilon_0$ satisfies the hypothesis of Corollary~\ref{cor:continuation-lattice} in the normalized massive network, for all $N\geq1$ and $0<m\leq1$.
\end{lemma}
\begin{proof}
We check that the Gamma upper tail is smaller than the activation parameter
$\varepsilon_0$, as required for the conditional comparison.

Integration by parts, followed by $t^{-1/2}\leq x^{-1/2}$ for $t\geq x>0$, gives
\begin{align*}
 Q_{3/2}(x)
 &=\frac{2}{\sqrt\pi}\int_x^\infty t^{1/2}e^{-t}\,dt\\
 &\leq\frac{2\sqrt x+x^{-1/2}}{\sqrt\pi}\,e^{-x}.
\end{align*}
At $x=512$, $\sqrt{512}<23$, $512^{-1/2}<1$, and $\sqrt\pi>1$, so the prefactor is less than $47<64$. Since $e>2$,
\[
 Q_{3/2}(512)<64e^{-512}<2^6 2^{-512}=2^{-506}<2^{-440}.
\]
Finally $\lambda_*/8=4096\geq512$, and the Gamma upper tail is decreasing. Corollary~\ref{cor:continuation-lattice} gives $\widehat g_h\leq8$ at every executed root, so its historywise posterior tail is no larger than $Q_{3/2}(\lambda_*/8)$.
\end{proof}

\paragraph{Finite crossing duality.}

This is the classical square-lattice matching-pair alternative; see Penrose--Rosoman~\cite[Section 3]{PenroseRosoman} and Ganesan~\cite[Theorem 2(i)]{Ganesan}. Its role is deterministic: the complement of a low nearest-neighbor crossing is a high star crossing. No independence of the field-induced colors is required. We retain the short verification to fix the rectangle sides, endpoint convention and allocation of equality to the low set.

\begin{lemma}[Nearest-neighbor low crossing versus star high crossing]\label{lem:routes-duality}
For any deterministic real function $f$ on $B_N$ and any $A\geq0$, exactly one of the following holds: there is a nearest-neighbor left-right path in $B_N$ with $|f|\leq A$ everywhere, or there is a star top-bottom path in $B_N$ with $|f|>A$ everywhere. Consequently, for the crossing functional~\eqref{eq:H},
\begin{equation}\label{eq:routes-crossing-duality}
 \{H_N(f)>A\}
 =\{\text{a star top-bottom path in }B_N\text{ has }|f|>A\}.
\end{equation}
\end{lemma}
\begin{proof}
The deterministic square-lattice matching-pair alternative is stated in \cite[Section 3, first paragraph]{PenroseRosoman}; a direct proof for arbitrary finite cell rectangles is \cite[Theorem 2(i), with the definitions on pp.\ 5--6]{Ganesan}. Apply it by coloring a site low when $|f(v)|\leq A$ and high otherwise. To match the cell formulation exactly, send $v=(i,j)\in B_N$ to
\[
 [i+N,i+N+1]\times[j+N,j+N+1]\subset[0,2N+1]^2.
\]
Side-sharing cells correspond to nearest neighbors, and cells sharing a corner correspond to star neighbors. Cells on each rectangle side correspond exactly to sites on the associated side of $B_N$. Hence the cited occupied-plus left-right versus vacant-star top-bottom alternative gives the assertion.

An unrestricted crossing can be trimmed after its last visit to its starting side before its first visit to the opposite side. This yields the endpoint-only convention used in the cited formulation, without increasing the maximum of $|f|$. Conversely every such trimmed path is an unrestricted crossing. There are finitely many simple paths, so $H_N(f)\leq A$ precisely when a low nearest-neighbor crossing exists. Taking its complement proves~\eqref{eq:routes-crossing-duality}. The low/high partition is $\leq A$ versus $>A$; equality is assigned to the low side, so no atomlessness argument or boundary adjustment is required.
\end{proof}

The alternative is deterministic, so it applies directly to the correlated colors defined by a field.

\paragraph{The massive comparison.}

The next estimate combines the preceding ingredients in their required order. Matching-pair duality identifies crossing failure with a high star event; the conditional exploration bounds that event by the integrated incoming-activation model; and the independent block routes embed this model in Lupu's exact augmented graph. The numerical Gamma-tail check makes the comparison uniform in the parameters stated below. This organization separates the classical duality and coupling inputs from the particular finite comparison proved here.

\begin{theorem}[Uniform massive comparison M1]\label{thm:routes-M1}
For every integer $N\geq1$ and every $0<m\leq1$, with the exact full-domain augmented graph used in~\eqref{eq:theta},
\begin{equation}\label{eq:M1}
\PP\bigl(H_N(\psi_{N,m})>256\bigr)\leq\theta_N(m).
\end{equation}
More precisely the normalized field satisfies
\begin{equation}\label{eq:routes-normalized-M1}
 \PP\bigl(H_N(\chi)>256\bigr)\leq\theta_N(m).
\end{equation}
\end{theorem}
\begin{proof}
We use crossing duality, the conditional comparison with independent coins,
and the block embedding into the exact augmented graph, in that order.
The last step returns to the original field normalization.

For $S=\operatorname{TOP}_N$ and $\mathcal T=\operatorname{BOT}_N$, Lemma~\ref{lem:routes-duality}, $X\overset d=\chi^2/2$, and $\sqrt{2\lambda_*}=256$ give
\[
 \PP\bigl(H_N(\chi)>256\bigr)
 =\PP\bigl(\operatorname{High}_{J,S,\mathcal T}(\lambda_*)\bigr).
\]
The underlying marked soup has exactly the normalized law by Lemma~\ref{lem:routes-marginal}. Apply Corollary~\ref{cor:continuation-lattice} with an auxiliary iid coin family independent of that entire source. Lemma~\ref{lem:routes-tail} verifies its common parameter condition. After integrating, Theorem~\ref{thm:routes-enhancement} gives
\begin{align*}
 \PP\bigl(H_N(\chi)>256\bigr)
 &\leq\Theta^{\mathrm{in}}_{D,J,S,\mathcal T}(\varepsilon_0)\\
 &\leq\PP_{\mathsf{Aug}(\widehat{\mathsf Q})}
       (\operatorname{TOP}_N\longleftrightarrow\operatorname{BOT}_N)\\
 &=\theta_N(m).
\end{align*}
The last equality is the time-rescaling invariance proved in Section~\ref{sec:finite-networks}: dividing conductances and killing by $a$ multiplies all occupation times by $a$, preserves the nontrivial traces, and leaves $2c_e\sqrt{X_xX_y}$ unchanged. Thus the entire normalized augmented graph has the same law as the unnormalized graph defining $\theta_N(m)$.

Finally the minimax definition gives $H_N(\chi)=\sqrt a\,H_N(\psi_{N,m})$ pointwise when $\chi=\sqrt a\,\psi_{N,m}$. Since $a\geq1$,
\[
 \{H_N(\psi_{N,m})>256\}\subseteq\{H_N(\chi)>256\}.
\]
This proves~\eqref{eq:M1}. All constructions are finite and valid for $N=1$; all eight destination neighbors remain in $D$. No lower bound on $m$ beyond $m>0$ and no sufficiently-large-$N$ qualification is introduced.
\end{proof}

\section{Probes and massive connectivity}\label{sec:probes}

We estimate the augmented top--bottom connection probability
$\theta_N(m)$ by adding two auxiliary vertices, one near each side of the
crossing. Given a crossing selected from the original augmented graph, loops
visiting the new vertices give independent opportunities to join suitable
portions of that crossing. A bound on the Gaussian correlation between the
two auxiliary vertices then gives an upper bound for their connection
probability. Together, these estimates will give
\[
 \theta_N(K/N)\leq
 \min\left\{1,\frac{2^{28}}3K^3e^{-K/2}\right\},
 \qquad K\geq8,\quad N\geq K.
\]
We first couple the original and enlarged networks so that loops avoiding
the new vertices have the original law. We then estimate the probabilities
of the two connections and compute the enlarged Gaussian covariance. This
argument does not use Theorem~\ref{thm:routes-M1}; we combine their conclusions
only afterward. Connections may use all of $D=B_{2N}$.

The massive-loop literature gives a useful perspective on this estimate. For a full-plane massive Brownian loop soup with mass bounded below by a positive constant and intensity at most $1/2$, Camia proves an exponential diameter tail for the filled cluster containing a fixed deterministic point, uniformly in that point; see \cite[Theorem 2.4, second bullet]{Camia}. Papon couples massive CLE$_4$ to the massive planar GFF and identifies it with the outermost boundaries of the massive Brownian half-soup in her domain and mass setting, which includes bounded simply connected domains with bounded continuous nonnegative mass; see \cite[Theorems 1.3 and 1.8]{Papon}.

After rescaling distances by $N$, the lattice mass $m=K/N$ corresponds to a mass of order $K$, so these continuum results are directly relevant to the decay mechanism. The estimate below keeps track of the particular finite augmented graph, its boundary connection event, and the dependence on both $N$ and $K$. We prove that estimate within the finite network; a transfer from the continuum statements with this uniformity and these explicit constants is not used.

\subsection{The original and enlarged networks}\label{sec:probe-extension}

\subsubsection{Extending the network without changing the original holding rates}

The construction uses the classical Poisson restriction principle~\cite{LastPenrose,LeJan} and the augmented-edge coupling of Lupu~\cite[Theorem 1 bis]{Lupu}. It preserves the original marked loop law for loops avoiding the new vertices and couples the two augmented graphs using the same uniforms on old edges.

Le Jan's loop-measure framework and the generator-based treatment of Chang and Le Jan explain restriction at the level of loops with their holding times; see \cite{LeJan,ChangLeJan}. This viewpoint is useful here because preserving the old principal precision preserves both the old jump weights and the total holding rates. The conductances to the new vertices are paid for by an equal reduction in the old killing. Thus the loop measure restricted to paths that avoid the new vertices has exactly the original marked law. The construction below records this classical restriction principle together with the shared-uniform realization of the two augmented graphs.

Start with a finite network $(V,c,k,\mathsf Q)$ as in Section~\ref{sec:finite-networks}. Add a finite set $W$ of new vertices and deterministic conductances $h_{vw}\geq0$ from $v\in V$ to $w\in W$. Keep old-old conductances unchanged and put
\begin{equation}\label{eq:budget}
 k'_v=k_v-\sum_{w\in W}h_{vw}\geq0.
\end{equation}
Choose nonnegative killing and conductances among the new vertices, with zero diagonal conductances, so that the enlarged precision $\mathsf Q^+$ is positive definite. Its old diagonal rates are still $a_v$, and
\begin{equation}\label{eq:principal}
 \mathsf Q^+_{VV}=\mathsf Q.
\end{equation}

In an enlarged half-soup, let $\mathcal L_0$ be its nontrivial loops avoiding $W$ and let $\mathcal L_W$ be its nontrivial loops visiting $W$. Retain all old trivial occupations $\xi_v$ from that same soup. Let $X^0$ be the occupation from $\mathcal L_0$ and the old $\xi_v$, and let $X^+$ be the full enlarged occupation. Write $n_e^0,n_e^+$ for the respective traversal counts of an old edge.

Take one family of independent uniforms $U_e$ on all enlarged edges, independent of the entire enlarged soup. On old edges use those same uniforms to define
\begin{align}
 \eta_e^0&=\ind_{\{n_e^0>0\}}\vee\ind_{\{U_e\leq p_e(X^0)\}},\label{eq:eta0}\\
 \eta_e^+&=\ind_{\{n_e^+>0\}}\vee\ind_{\{U_e\leq p_e(X^+)\}},\label{eq:etaplus}\\
 p_e(x)&=1-\exp(-2c_{uv}\sqrt{x_ux_v}),\qquad e=\{u,v\}.\nonumber
\end{align}
On new edges only $\eta^+$ is defined, by~\eqref{eq:finite-augmentation}. Set
\begin{equation}\label{eq:finite-F0}
 \mathcal F_0=\sigma\bigl(\mathcal L_0,(\xi_v)_{v\in V},
                    (U_e)_{e\text{ an old edge}}\bigr).
\end{equation}

\begin{theorem}[Exact restriction and shared-uniform embedding]\label{thm:finite-coupling}
Under the preceding finite-network hypotheses, the construction has the following properties.
\begin{enumerate}
\item The marked soup $\mathcal L_0$, together with the old trivial occupations, has exactly the original marked half-soup law on $V$. The same assertion holds with the entire trivial-loop clouds retained.
\item $\eta^0\sim\mathsf{Aug}(\mathsf Q)$ and $\eta^+\sim\mathsf{Aug}(\mathsf Q^+)$. These are the full graph marginals, not only occupation or partition marginals.
\item For every old edge, $\eta_e^0\leq\eta_e^+$ almost surely. In particular, every baseline open path survives in the enlarged graph.
\item $\mathcal F_0$ is independent of $\mathcal L_W$, the new trivial occupations, and the new-edge uniforms. Any deterministically selected measurable baseline path has the same independence from those variables.
\item If an enlarged loop visits a new vertex $w$ and a vertex $v$ of a baseline open path, then $w$ is connected to that path in $\eta^+$ through the loop's traversed edges. Two loops attaching to vertices of the same baseline component connect their new vertices in $\eta^+$.
\end{enumerate}
\end{theorem}
\begin{proof}
\textbf{Exact marked restriction.} Every term in~\eqref{eq:finite-marked} for a loop whose vertices lie in $V$ uses only old-old conductances and old diagonal rates. Both are unchanged. Therefore, as measures on the full cyclic marked loop space,
\begin{equation}\label{eq:finite-measure-restriction}
 \mu_{\mathsf Q^+}^{\mathrm{nt}}\big|_{\{\ell:\ell\cap W=\varnothing\}}
 =\mu_{\mathsf Q}^{\mathrm{nt}}.
\end{equation}
The old trivial measures $e^{-a_vt}\,dt/t$ are also unchanged. Restricting the enlarged Poisson soup to these classes proves the first conclusion, including all holding marks and repeated-visit factors. This is an equality of measures before any conditioning on occupations.

\textbf{Both graph marginals.} The first conclusion gives the exact input soup for $\eta^0$. Its uniforms are independent and have the correct product law even conditional on that soup. Thus~\eqref{eq:eta0} has law $\mathsf{Aug}(\mathsf Q)$. The full enlarged soup has its original law, and the entire uniform family is independent of it; hence~\eqref{eq:etaplus}, together with the new edges, has law $\mathsf{Aug}(\mathsf Q^+)$. Sharing some uniforms with another graph changes the joint coupling, not either marginal.

\textbf{Pathwise inclusion.} For each old vertex and edge,
\[
 X_v^0\leq X_v^+,\qquad n_e^0\leq n_e^+.
\]
The opening function $p_e$ is increasing in each endpoint occupation. If an old edge is traversed in the baseline, it is traversed in the full soup. If it is opened instead by $U_e\leq p_e(X^0)$, then $U_e\leq p_e(X^+)$. Both cases give $\eta_e^0\leq\eta_e^+$.

\textbf{Independence.} Avoiding and visiting $W$ are deterministic disjoint loop classes. Their Poisson restrictions are independent. Trivial-loop clouds at distinct vertices are independent of these nontrivial restrictions. All uniforms were chosen independently of all loop variables, and the old and new uniform families are independent. These facts give the asserted joint independence of $\mathcal F_0$ from the visiting soup, new clocks, and new uniforms. A measurable function of $\mathcal F_0$, including a selected path, inherits it.

\textbf{Attachment.} Every edge traversed by an enlarged loop is open in $\eta^+$ by definition. A finite segment of that loop connects $w$ to the visited vertex $v$. The baseline path and every baseline connection remain present by edgewise inclusion. Concatenating these connections proves the fifth conclusion.
\end{proof}

The original augmented graph uses the marked loops that avoid the added vertices, the trivial occupations at old vertices, and the original edge uniforms. The theorem preserves this graph law and realizes it as a subgraph of the enlarged augmented graph.

\subsubsection{Conditional independence after path selection}

The independence used here is the independence of Poisson restrictions to disjoint deterministic classes, as in \cite{LastPenrose,LeJan}. The order of revealing information is essential to its application. We first select the path using only the baseline sigma-field $\mathcal F_0$, which is independent of the visiting-loop processes. Conditional on this sigma-field, the selected vertex sets are fixed, and each visiting process can be restricted to the loops that hit its selected set. This gives the conditional product used below. Averaging that product over the selected path is a separate step.

Let $\mathcal A_1,\mathcal A_2$ be disjoint deterministic measurable subclasses of the nontrivial probe-visiting loops. Let $L_1,L_2\subseteq V$ be random vertex sets measurable with respect to $\mathcal F_0$. Write
\[
 \mathcal A_i(L_i)=\{\ell\in\mathcal A_i:\ell\text{ visits }L_i\},\qquad
 \Lambda_i(L_i)=\tfrac12\mu_{\mathsf Q^+}^{\mathrm{nt}}(\mathcal A_i(L_i)).
\]
The means are finite by Lemma~\ref{lem:finite-marked}. Empty selected sets are allowed and give mean zero.

\begin{corollary}[Exact conditional attachment formula]\label{cor:finite-attachment}
Conditional on $\mathcal F_0$, the numbers of loops in $\mathcal A_1(L_1)$ and $\mathcal A_2(L_2)$ are independent Poisson variables with means $\Lambda_1(L_1)$ and $\Lambda_2(L_2)$. If $E\in\mathcal F_0$ ensures that $L_1,L_2$ lie in a single baseline component, and every loop in $\mathcal A_i$ visits a specified new vertex $w_i$, then
\begin{equation}\label{eq:finite-transfer}
 \PP(w_1\longleftrightarrow w_2\text{ in }\eta^+)
 \geq\EE\!\left[\ind_E
       (1-e^{-\Lambda_1(L_1)})(1-e^{-\Lambda_2(L_2)})\right].
\end{equation}
\end{corollary}
\begin{proof}
The restrictions to $\mathcal A_1$ and $\mathcal A_2$ are independent Poisson processes and are jointly independent of $\mathcal F_0$. In particular, selecting $L_1,L_2$ from $\mathcal F_0$ does not reveal any of the probe-visiting loops.
After conditioning on $\mathcal F_0$, the two sets $L_i$ are fixed; restricting each process to loops hitting its fixed set gives the stated Poisson counts. Their conditional probability of both being nonzero is the product in~\eqref{eq:finite-transfer}. On $E$ and this event, the last conclusion of Theorem~\ref{thm:finite-coupling} connects $w_1$ to $w_2$. Integrating the conditional probability proves the inequality.
\end{proof}

The product in~\eqref{eq:finite-transfer} is a conditional probability given $\mathcal F_0$. We average over the original graph and selected path only afterward. The two loop classes remain part of the enlarged soup, and the connection event refers to its augmented graph.

\subsubsection{The two-probe network}

We will use the following deterministic network extension in the connection estimate of Section~\ref{subsec:probe-connections}. Let $D=B_{2N}$, $N\geq1$, and $0<m\leq1$. On $D$ use old conductances $1/4$ and killing
\[
 k_v=m^2+1-\deg_D(v)/4.
\]
Then every old diagonal rate is $a_v=1+m^2$ and the old precision is $\mathsf Q_m=\mathsf Q_D+m^2I$. Define the deterministic lattice rectangles
\begin{align*}
 U_+&=\bigl([-3N/2,3N/2]\times[N/2,3N/2]\bigr)\cap\mathbb Z^2,\\
 U_-&=\{(x,-y):(x,y)\in U_+\}.
\end{align*}
Real endpoints followed by intersection with the lattice specify rounding exactly. Both sets are nonempty, disjoint, and strictly inside $D$: they contain $(0,N)$ and $(0,-N)$, respectively, and $\lfloor3N/2\rfloor\leq2N-1$.

Put $\kappa=m^2/2$. Add probes $a,b$, join $a$ to every vertex of $U_+$ and $b$ to every vertex of $U_-$ with conductance $\kappa$, and set
\[
 k'_v=k_v-\kappa\ind_{\{v\in U_+\cup U_-\}}.
\]
There is no $a$--$b$ edge and no killing at either probe. Their diagonal rates are
\[
 d_a=\kappa|U_+|>0,\qquad d_b=\kappa|U_-|>0.
\]
Their trivial occupations have shape $1/2$ and rates $d_a,d_b$. Uniform rescaling of the enlarged network rescales these probe rates and all probe conductances together, as in Section~\ref{sec:finite-networks}. With $\mathbf1_U$ denoting the indicator column indexed by $D$, the enlarged precision is
\begin{equation}\label{eq:finite-probe-matrix}
 \mathsf Q^+=
 \begin{pmatrix}
 \mathsf Q_m&-\kappa\mathbf1_{U_+}&-\kappa\mathbf1_{U_-}\\
 -\kappa\mathbf1_{U_+}^{T}&d_a&0\\
 -\kappa\mathbf1_{U_-}^{T}&0&d_b
 \end{pmatrix}.
\end{equation}

\begin{lemma}[All probe-extension hypotheses hold]\label{lem:finite-probes}
For every $N\geq1$ and $0<m\leq1$, this deterministic extension satisfies Theorem~\ref{thm:finite-coupling}. In particular $\mathsf Q^+_{DD}=\mathsf Q_m$, the enlarged network is transient, and the baseline graph has exactly the massive graph law used in $\theta_N(m)$.
\end{lemma}
\begin{proof}
We check the two remaining conditions needed for the extension theorem: the old diagonal rates are unchanged, and the enlarged precision is positive definite.
Disjointness of the rectangles gives $k'_v\geq m^2/2>0$ at every old vertex. The added conductance at $v$ equals the amount removed from its killing, so its total rate remains $1+m^2$. This proves the principal-block identity.

For a real function $z$ on $D\cup\{a,b\}$, its energy is
\begin{align*}
 z^T\mathsf Q^+z={}&\frac14\sum_{\substack{\{u,v\}\subset D\\|u-v|_1=1}}
                   (z_u-z_v)^2+\sum_{v\in D}k'_vz_v^2\\
 &+\kappa\sum_{v\in U_+}(z_v-z_a)^2
  +\kappa\sum_{v\in U_-}(z_v-z_b)^2,
\end{align*}
where the old-edge sum is over unordered edges. All terms are nonnegative. Zero energy forces every old $z_v=0$ because $k'_v>0$; the nonempty probe neighborhoods then force $z_a=z_b=0$. Thus $\mathsf Q^+$ is positive definite. The extension theorem now gives the exact baseline graph law. The probe degrees may grow with $N$, but each network is finite. The construction uses $m>0$, as does the intended choice $m=K/N$ with $K>0$.
\end{proof}

\paragraph{Loop classes for the two attachments.}
Partition the enlarged nontrivial loop space into four deterministic classes:
\begin{enumerate}
\item $\mathcal C_0$: loops avoiding both probes;
\item $\mathcal C_a$: loops with exactly one cyclic occurrence of $a$, all other vertices lying in $U_+$;
\item $\mathcal C_b$: loops with exactly one cyclic occurrence of $b$, all other vertices lying in $U_-$;
\item $\mathcal C_{\mathrm{rest}}$: all remaining nontrivial loops.
\end{enumerate}
They are pairwise disjoint, and their Poisson restrictions are independent. A loop in $\mathcal C_a$ avoids $b$, and one in $\mathcal C_b$ avoids $a$. The holding marks are retained. ``One occurrence'' refers to the cyclic discrete word, with its terminal repeated endpoint omitted from the count.

\paragraph{Selecting the crossing path.}
Let $E_0$ be the event that $\operatorname{TOP}_N$ connects to $\operatorname{BOT}_N$ in the baseline graph $\eta^0$, allowing travel through all of $D$. On $E_0$, select the first simple open crossing path in a fixed deterministic ordering of the finite path family. The selection uses the full baseline graph and is $\mathcal F_0$-measurable. On $E_0^c$, use an empty dummy path.

For later geometric estimates, define
\[
 C_+=\bigl([-5N/4,5N/4]\times[3N/4,5N/4]\bigr)\cap\mathbb Z^2
\]
and its reflection $C_-$. These contain the top and bottom sides, respectively, and lie in $U_+$ and $U_-$. Orient the selected path from top to bottom. Its initial segment strictly \emph{before} its first exit from $C_+$ defines $L_+$; reverse the path and do the same for $C_-$ to define $L_-$. Thus $L_\pm\subset C_\pm\subset U_\pm$ literally, including rounding. Their coordinate-span and distance-to-boundary estimates are proved in Lemma~\ref{lem:potential-geometry}.

With
\[
 \Lambda_a(L_+)=\tfrac12\mu_{\mathsf Q^+}^{\mathrm{nt}}(\{\ell\in\mathcal C_a:\ell\cap L_+\ne\varnothing\}),
\]
and the analogous $\Lambda_b(L_-)$, Corollary~\ref{cor:finite-attachment} gives
\begin{equation}\label{eq:finite-DGFF-transfer}
 \PP_{\mathsf Q^+}(a\longleftrightarrow b)
 \geq\EE\!\left[\ind_{E_0}
 (1-e^{-\Lambda_a(L_+)})(1-e^{-\Lambda_b(L_-)})\right],
 \qquad \PP(E_0)=\theta_N(m).
\end{equation}
On $E_0^c$, set $L_+=L_-=\varnothing$, so both intensities in this unconditional expectation are defined everywhere. The left side uses the exact enlarged augmented graph. Theorem~\ref{thm:potential-inputs} supplies the quantitative potential estimates; Lemmas~\ref{lem:probes-attachment} and~\ref{lem:probes-arcsine} prove the lower and upper bounds for this exact interface.

\subsection{Potential estimates for the probe rectangles}\label{subsec:potential-estimates}

The next estimates apply to every path segment selected from an original crossing, regardless of which crossing is chosen. They will bound the probabilities that loops from the two auxiliary vertices hit these segments. The underlying semigroup and potential-theory tools are classical. Lawler and Limic give the random-walk and Green-function background in Section~1.2 and Chapter~4, and develop Dirichlet and capacity methods in Chapter~6; see \cite{LawlerLimic}. For the finite-network energy viewpoint, see \cite[Sections 2.4--2.5]{LyonsPeres}. We give the estimates in the killed-network normalization used here and retain the constants and pathwise quantifiers needed for the two probe attachments.

\subsubsection{The walk, killing, and occupation normalization}

Let $W_t$ be the continuous-time walk on $\mathbb Z^2$ which jumps to each of its four nearest neighbors at rate $1/4$. Its total jump rate is one and its generator is $P-I$, where $P$ is the discrete-time simple random walk transition matrix. Independently let $\zeta_m$ have exponential law of rate $m^2$, where $m>0$. For a finite nonempty set $U$, put
\[
 \tau_U=\inf\{t\geq0:W_t\notin U\},\qquad
 \sigma_U=\tau_U\wedge\zeta_m,
 \qquad \mathsf Q_{m,U}=I-P_U+m^2I.
\]
Here $P_U$ is the principal restriction of $P$, without renormalizing its rows. Write
\[
 G_{m,U}(x,y)=\EE_x\int_0^{\sigma_U}\ind_{\{W_t=y\}}\,dt,
 \qquad
 G_m(x,y)=\EE_x\int_0^{\zeta_m}\ind_{\{W_t=y\}}\,dt.
\]
The full-lattice kernel is translation invariant. For a nonempty $L\subset U$ we use
$T_L=\inf\{t\geq0:W_t\in L\}$, with $T_L=0$ if the starting point is in $L$. An infimum over an empty time set is $+\infty$.

\begin{lemma}[Exact resolvent and domain comparison]\label{lem:potential-normalization}
For $q=(1+m^2)^{-1}$,
\begin{equation}\label{eq:potential-resolvent}
 G_m=q\sum_{n\geq0}q^nP^n,\qquad
 G_{m,U}=q\sum_{n\geq0}q^nP_U^n=\mathsf Q_{m,U}^{-1}.
\end{equation}
The finite matrix $\mathsf Q_{m,U}$ is symmetric positive definite. Its Green kernel is nonnegative and symmetric, and
\begin{equation}\label{eq:potential-domain}
 0\leq G_{m,U}(x,y)\leq G_m(x,y),\qquad
 (G_{m,U}\mathbf1)(x)\leq m^{-2}.
\end{equation}
The same conclusions hold with $D$ in place of $U$.
\end{lemma}
\begin{proof}
We first obtain the Green resolvent by integrating the continuous-time walk. The same occupation formula will then give the comparison with the walk killed on leaving $U$.
Poissonizing the number of rate-one jumps gives
\[
 H_t=e^{t(P-I)}=e^{-t}\sum_{n\geq0}\frac{t^n}{n!}P^n.
\]
The killed semigroup has the same expression with $P_U$. Tonelli's theorem and independent mass killing give
\[
 \int_0^\infty e^{-m^2t}e^{-t}\frac{t^n}{n!}\,dt
 =(1+m^2)^{-n-1}=q^{n+1},
\]
which proves \eqref{eq:potential-resolvent}. The geometric series is an inverse since $q<1$ and $P_U$ is a contraction. Its symmetry follows from that of $P_U$. The quadratic form of $I-P_U$ is the nearest-neighbor edge energy of the zero extension, so adding $m^2\sum_{x\in U}f(x)^2$ makes it positive definite.

The occupation representation also gives
\[
 G_{m,U}(x,y)=\int_0^\infty e^{-m^2t}
       \PP_x(W_t=y,\ t<\tau_U)\,dt.
\]
Removing the exit restriction proves the entrywise upper bound in \eqref{eq:potential-domain}. Summing over $y$ and then dropping the survival probability bounds the row sum by $\int_0^\infty e^{-m^2t}dt=m^{-2}$.
\end{proof}

The leading factor $q$ and the $n=0$ term in \eqref{eq:potential-resolvent} are both retained. The latter is the first mean holding time before a jump or mass death. Discrete-time parity remains encoded in $P^n$; Poissonization is exact. These are the occupation conventions of Section~\ref{sec:finite-networks}. Restriction from $D$ to $U$ takes the principal killed precision.

\subsubsection{Rounded geometry for every full-domain crossing}

Recall the probe rectangles and their cores from Section~\ref{sec:probe-extension}:
\begin{align*}
 U_+&=([-3N/2,3N/2]\times[N/2,3N/2])\cap\mathbb Z^2,\\
 C_+&=([-5N/4,5N/4]\times[3N/4,5N/4])\cap\mathbb Z^2,
\end{align*}
and their reflections $U_-,C_-$ in the horizontal axis. For a vertex set $L$ define its span in coordinate $i$ by
\[
 \operatorname{span}_i(L)=\max_{v\in L}v_i-\min_{v\in L}v_i.
\]

\begin{lemma}[Strict pre-exit segments and deterministic bounds]\label{lem:potential-geometry}
Let $N\geq8$ be an integer. Every nearest-neighbor path in $D=B_{2N}$ from $\operatorname{TOP}_N$ to $\operatorname{BOT}_N$ has the strict pre-exit segments $L_+,L_-$ defined in Section~\ref{sec:probe-extension}. They are nonempty connected subsets of $C_+,C_-$, respectively, and
\begin{equation}\label{eq:potential-span}
 \max_{i=1,2}\operatorname{span}_i(L_\pm)
 \geq\lfloor N/4\rfloor\geq N/8.
\end{equation}
For either sign,
\begin{equation}\label{eq:potential-geometry}
 \operatorname{dist}_\infty(C_\pm,U_\pm^c)\geq N/4,
 \quad 3N^2\leq |U_\pm|\leq4N^2,
 \quad\operatorname{dist}_\infty(U_+,U_-)\geq N.
\end{equation}
These segments are measurable functions of the baseline path selected in Section~\ref{sec:probe-extension}. No restriction of that path to $B_N$ is required.
\end{lemma}
\begin{proof}
We first use the path's exit from the core to find a long segment that stays inside it. We then check the rectangle bounds and the measurability of the selected segment.
Write the path as $v_0,\ldots,v_s$, oriented from top to bottom, and let $t$ be its first index outside $C_+$. Its starting point is in $C_+$ and its final point is outside, so $1\leq t\leq s$. The retained segment consists of $v_0,\ldots,v_{t-1}$ and is literally contained in $C_+$. The exit step changes one coordinate. In that coordinate $v_{t-1}$ lies on the relevant lattice face.

For a horizontal-coordinate exit, the distance between that face and the starting coordinate in $[-N,N]$ is at least
\[
 \lfloor5N/4\rfloor-N=\lfloor N/4\rfloor.
\]
For an upper vertical exit the same formula applies to the starting height $N$. For a lower vertical exit the displacement is
$N-\lceil3N/4\rceil=\lfloor N/4\rfloor$. Thus one coordinate span has the stated lower bound. Reversing the path and reflecting proves the assertion for $L_-$. Later excursions anywhere in $D$ do not affect these first-exit statements.

The span estimate used the path. The distance, area, and separation bounds depend only on the rectangles.
A point of the real core rectangle is at least $N/4$ from each corresponding real outer face. A lattice point outside the outer rectangle lies strictly beyond one of those faces. This proves the first bound in \eqref{eq:potential-geometry}. Explicitly, the right-hand gap is
$\lfloor3N/2\rfloor+1-\lfloor5N/4\rfloor>N/4$; the other faces satisfy the same bound. No exit vertex is added to the core.

For odd $N$, $|U_+|=3N^2$; for even $N$, $|U_+|=(3N+1)(N+1)$. These give the area bounds for $N\geq8$. Reflection preserves area. The lowest height in $U_+$ is $\lceil N/2\rceil$ and the highest height in $U_-$ is $-\lceil N/2\rceil$, proving their separation.

The baseline path is the first available simple path in a fixed ordering of a finite family. Taking its first exit and reversing its terminal segment are deterministic operations. They preserve its $\mathcal F_0$-measurability. On the event $E_0^c$ of no baseline crossing, the empty dummy sets from Section~\ref{sec:probe-extension} remain in use; the positive lower bounds here concern the nonempty segments on $E_0$.
\end{proof}

Since the projection of a connected nearest-neighbor path contains every integer between its extrema, an eligible segment contains vertices $v_1,\ldots,v_M$ with consecutive distinct values of one coordinate, where
\begin{equation}\label{eq:potential-level-count}
 M\geq\lfloor N/4\rfloor+1>N/4.
\end{equation}
One may choose the first occurrence at each level along the selected path. The chosen vertices need not themselves form a nearest-neighbor path.

\subsubsection{A one-dimensional Green sum and a spatial bound}

Michta and Slade give a detailed analysis of the massive lattice Green function, including sharp long-distance asymptotics, an exact one-dimensional formula, and a Bessel-function representation; see \cite[Theorem 1.3, Section 1.4.3 and Section 2]{MichtaSlade}. Their definitions (1.1)--(1.3), in dimension two with their parameter $a=m$, give the same infinite-lattice resolvent $G_m$ used here. This supplies useful context for the elementary coordinate calculation below. We retain the latter to obtain the stated inequalities throughout the finite mass and distance ranges used by the probe argument.

The Green estimate can be obtained directly from the two independent coordinate walks. This avoids a short-time Gaussian approximation and produces an explicit constant.

\begin{lemma}[Axial envelope and its exact sum]\label{lem:potential-axis}
For $m>0$ and integer $j\geq0$, put $g_m(j)=G_m(0,(j,0))$. Then $g_m$ is nonnegative and nonincreasing, and
\begin{equation}\label{eq:potential-envelope}
 \sup_{|z|_\infty\geq j}G_m(0,z)=g_m(j).
\end{equation}
Moreover,
\begin{equation}\label{eq:potential-axis-sum}
 \sum_{j\in\mathbb Z}g_m(|j|)
 =\frac1{m\sqrt{1+m^2}}\leq\frac1m.
\end{equation}
Consequently,
\begin{equation}\label{eq:potential-radial-sum}
 \sum_{j\geq0}\sup_{|z|_\infty\geq j}G_m(0,z)\leq\frac1m,
 \qquad g_m(j)\leq\frac1{m(j+1)}.
\end{equation}
\end{lemma}
\begin{proof}
The proof has two parts. We first use one-dimensional heat-kernel monotonicity to bound the Green function by its values on a coordinate axis. We then compute the exact sum along that axis, which will control the Green energy of selected path vertices.
Each coordinate of $W$ jumps to either neighbor at rate $1/4$, and the two coordinates are independent. Let $p_t(k)$ be the one-dimensional kernel from zero. Thus
\[
 H_t(0,(x,y))=p_t(x)p_t(y),\qquad p_t(-k)=p_t(k).
\]
We first prove that $p_t(k)$ decreases for $k\geq0$. Set $d_k(t)=p_t(k)-p_t(k+1)$. The heat equation gives
\[
 d'_0=(d_1-3d_0)/4,\qquad
 d'_k=(d_{k-1}+d_{k+1}-2d_k)/4\quad(k\geq1),
 \qquad d(0)=(1,0,0,\ldots).
\]
Let $\mathcal B$ denote this bounded operator on $\ell^\infty(\mathbb Z_{\geq0})$. Every matrix entry of $I+\mathcal B$ is nonnegative, and its row sums are at most one. Therefore
\[
 e^{t\mathcal B}=e^{-t}\sum_{n\geq0}\frac{t^n}{n!}(I+\mathcal B)^n
\]
preserves nonnegativity. Uniqueness for the bounded linear differential equation identifies $d(t)$ with this semigroup applied to $d(0)$. Hence $d_k(t)\geq0$.

If $|z|_\infty\geq j$, symmetry and this monotonicity imply
\[
 p_t(z_1)p_t(z_2)\leq p_t(j)p_t(0).
\]
The fixed point $(j,0)$ attains the bound at every time. Integration against $e^{-m^2t}$ proves \eqref{eq:potential-envelope}; there is no exchange of an unattained supremum with an integral.

Tonelli's theorem and $\sum_{j\in\mathbb Z}p_t(j)=1$ now give
\[
 \sum_{j\in\mathbb Z}G_m(0,(j,0))
 =\int_0^\infty e^{-m^2t}p_t(0)\,dt.
\]
To evaluate this one-dimensional return resolvent, write
$R_m(k)=\int_0^\infty e^{-m^2t}p_t(k)\,dt$. It is bounded and satisfies
\[
 (m^2+1/2)R_m(k)-\tfrac14\bigl(R_m(k-1)+R_m(k+1)\bigr)
 =\ind_{\{k=0\}}.
\]
Put $\rho=(\sqrt{1+m^2}-m)^2$. Then $0<\rho<1$ and
$\rho+\rho^{-1}=2+4m^2$. Direct substitution, including the equation at zero, gives the bounded solution
\[
 R_m(k)=\frac{\rho^{|k|}}{m\sqrt{1+m^2}}.
\]
It is the unique bounded solution: a homogeneous bounded solution $f$ would satisfy
$\|f\|_\infty\leq(1+2m^2)^{-1}\|f\|_\infty$ and hence vanish. Evaluating at zero proves \eqref{eq:potential-axis-sum}. Its nonnegative half-axis sum is at most the full sum. Finally,
$(j+1)g_m(j)\leq\sum_{i=0}^jg_m(i)\leq1/m$, which proves \eqref{eq:potential-radial-sum}.
\end{proof}

The exact two-sided sum in \eqref{eq:potential-axis-sum} will be used in the path-energy estimate. Retaining it avoids an unnecessary factor of two from separately bounding the two sides of a selected vertex.

\subsubsection{Exponential decay from a coordinate hitting time}

\begin{lemma}[A uniform separated Green bound]\label{lem:potential-decay}
For $0<m\leq1$ and $r=|z|_\infty\geq1$,
\begin{equation}\label{eq:potential-decay}
 G_m(0,z)\leq\frac4{mr}e^{-mr/2}.
\end{equation}
For any $R>0$ and $\alpha>0$ with $mR\geq\alpha$,
\begin{equation}\label{eq:potential-separated}
 \sup_{|z|_\infty\geq R}G_m(0,z)\leq\frac4\alpha e^{-mR/2}.
\end{equation}
In particular, for fixed $c_0>0$ and $mN\geq1$, the right-hand side at $R=c_0N$ is bounded by $(4/c_0)e^{-c_0mN/2}$. All these upper bounds persist after domain killing.
\end{lemma}
\begin{proof}
For a single coordinate started at zero, let $T^{(1)}_\ell$ be its first hitting time of the nonnegative integer $\ell$. If $\ell=0$ this time is zero. For $\ell>0$, set
\[
 \nu_m=2\operatorname{arsinh}(m),\qquad
 \tfrac12(\cosh\nu_m-1)=m^2.
\]
The independent-increment formula for the rate-$1/4$ jumps in each direction shows that
$\exp\{\nu_m W_t^{(1)}-m^2t\}$ is a martingale. Stopped at $T^{(1)}_\ell\wedge t$, it is bounded by $e^{\nu_m\ell}$, since the nearest-neighbor coordinate cannot pass $\ell$ before hitting it. On the event of no hit, its limit is zero. Optional stopping and dominated convergence give
\begin{equation}\label{eq:potential-hit-transform}
 \EE_0[e^{-m^2T^{(1)}_\ell};T^{(1)}_\ell<\infty]
 =e^{-\nu_m\ell}.
\end{equation}
No assertion of almost-sure hitting is needed. Also
$\nu_m=2\int_0^m(1+s^2)^{-1/2}ds\geq m$ for $m\leq1$.

Reflect and exchange coordinates so that $z_1=r$. Put $\ell=\lfloor r/2\rfloor$. This split lets us apply the hitting-time transform to the first part of the displacement and the axial Green bound to the distance that remains.
The walk cannot occupy $z$ before $T^{(1)}_\ell$. By the strong Markov property and the occupation representation,
\[
 G_m(0,z)
 =\EE_0\!\left[e^{-m^2T^{(1)}_\ell}
       G_m(W_{T^{(1)}_\ell},z);T^{(1)}_\ell<\infty\right].
\]
At this time the first-coordinate distance from $z$ is $r-\ell$. Lemma~\ref{lem:potential-axis} and \eqref{eq:potential-hit-transform} therefore yield
\begin{align*}
 G_m(0,z)&\leq e^{-\nu_m\ell}g_m(r-\ell)
 \leq\frac{e^{-m\lfloor r/2\rfloor}}{m(r-\lfloor r/2\rfloor+1)}\\
 &\leq\frac{2\sqrt e}{mr}e^{-mr/2}
 \leq\frac4{mr}e^{-mr/2}.
\end{align*}
Here $\lfloor r/2\rfloor\geq(r-1)/2$, $m\leq1$, and $e<4$ were used. This includes $r=1$, when the stopping time is zero. If $r\geq R$ and $mR\geq\alpha$, the bound implies \eqref{eq:potential-separated}. Domain comparison follows from Lemma~\ref{lem:potential-normalization}.
\end{proof}

This proof controls every time scale in the exact resolvent. In particular, it does not apply a Gaussian bound of the form $e^{-c|z|^2/t}/(1+t)$ to the small-time continuous-time kernel.

\subsubsection{Equilibrium measure and capacity}

The equilibrium representation and reciprocal minimum-energy characterization in the next lemma are classical finite-network identities. Sznitman gives them in \cite[author notes, Propositions 1.8--1.9]{SznitmanTopics}; the last-exit representation also gives the hitting identity after averaging over the starting point. These references provide an exact source for the underlying identities, in addition to the broader random-walk and electrical-network accounts cited above.

For the convention dictionary, put conductance $c_{xy}=1/4$ on each nearest-neighbor edge within $U$ and killing $\kappa_x=(4-\deg_U(x))/4+m^2$. The resulting precision is $\mathsf Q_{m,U}$. In the variable-rate chain of \cite[author notes, Section 1.6]{SznitmanTopics}, its inverse is the expected raw occupation time at a vertex, exactly the Green convention used below. The equilibrium charge is therefore $e_L=\mathsf Q_{m,U}h_L$, and capacity is its total mass. If $U$ is disconnected, these identifications apply on each component. We retain the short proof to make this normalization explicit; the subsequent geometric estimates still have to control every eligible segment selected from the crossing.

\begin{lemma}[Equilibrium, energy, and hitting identities]\label{lem:potential-equilibrium}
Let $m>0$, let $U$ be finite and nonempty, and let $\varnothing\ne L\subset U$. Define
\[
 h_L(x)=\PP_x(T_L<\sigma_U),\qquad e_L=\mathsf Q_{m,U}h_L.
\]
Then $e_L$ is supported on $L$, is nonnegative, and in fact $e_L(v)\geq m^2$ for $v\in L$. With
\[
 \operatorname{cap}_{m,U}(L)=\sum_{v\in L}e_L(v),
\]
one has
\begin{equation}\label{eq:potential-equilibrium}
 h_L=G_{m,U}e_L,\qquad
 e_L^TG_{m,U}e_L=\operatorname{cap}_{m,U}(L).
\end{equation}
For every probability measure $\nu$ supported on $L$,
\begin{equation}\label{eq:potential-energy}
 \operatorname{cap}_{m,U}(L)\geq
 \frac1{\nu^TG_{m,U}\nu}.
\end{equation}
Equality is attained by $\nu=e_L/\operatorname{cap}_{m,U}(L)$. If $\mu_U$ is uniform on $U$, then
\begin{equation}\label{eq:potential-uniform-identity}
 \PP_{\mu_U}(T_L<\sigma_U)
 =\frac1{|U|}\sum_{v\in L}e_L(v)(G_{m,U}\mathbf1)(v).
\end{equation}
\end{lemma}
\begin{proof}
We first identify the charge supported on $L$. The Green inner product will then give the capacity inequality, and averaging the hitting probability will give the final identity.
Extend $h_L$ by zero off $U$. It takes values in $[0,1]$ and equals one on $L$, including the starting-in-$L$ case. For $x\in U\setminus L$, the first jump or death gives
\[
 (1+m^2)h_L(x)=\tfrac14
       \sum_{\substack{y\in U\\|y-x|_1=1}}h_L(y).
\]
Thus $e_L=0$ off $L$. For $v\in L$,
\[
 e_L(v)=1+m^2-\tfrac14
       \sum_{\substack{y\in U\\|y-v|_1=1}}h_L(y)
 \geq m^2+1-\deg_U(v)/4\geq m^2.
\]
Inverting the precision gives the first identity in \eqref{eq:potential-equilibrium}. Since $e_L$ is supported on $L$ and $h_L=1$ there, its Green energy is its total mass, proving the second identity.

For a probability measure $\nu$ on $L$, $\nu^TG_{m,U}e_L=\nu^Th_L=1$. Cauchy--Schwarz in the positive definite Green inner product gives
\[
 1\leq (\nu^TG_{m,U}\nu)(e_L^TG_{m,U}e_L),
\]
which proves \eqref{eq:potential-energy}. The stated normalized equilibrium measure is a probability measure and attains equality. Finally,
\[
 \PP_{\mu_U}(T_L<\sigma_U)=\mu_U^Th_L
 =\mu_U^TG_{m,U}e_L.
\]
Symmetry of $G_{m,U}$ and $\mu_U=|U|^{-1}\mathbf1$ give \eqref{eq:potential-uniform-identity}.
\end{proof}

The definition $e_L=\mathsf Q_{m,U}h_L$ fixes capacity in continuous-time conductance units. No discrete-time escape-probability factor is omitted. It also avoids an ambiguous ``positive return time'' during the initial holding interval. For example, when $U=L$ is a singleton, $G_{m,U}=(1+m^2)^{-1}$ and the capacity is $1+m^2$, consistently with all three identities.

\subsubsection{Capacity of every eligible path segment}

\begin{lemma}[Coordinate levels give a uniform capacity bound]\label{lem:potential-capacity}
Let $N\geq8$, $m=K/N>0$, and let $L\subset C\subset U$ be either eligible segment from Lemma~\ref{lem:potential-geometry}. Then
\begin{equation}\label{eq:potential-capacity}
 \operatorname{cap}_{m,U}(L)\geq K/4.
\end{equation}
The same conclusion holds for any nonempty connected $L\subset C$ with a coordinate span at least $\lfloor N/4\rfloor$.
\end{lemma}
\begin{proof}
Choose $v_1,\ldots,v_M\in L$ at consecutive values of that coordinate, with $M$ as in \eqref{eq:potential-level-count}, and let $\nu=M^{-1}\sum_{i=1}^M\delta_{v_i}$. For every $i,j$,
$|v_i-v_j|_\infty\geq|i-j|$. Domain comparison and the axial envelope imply
\[
 G_{m,U}(v_i,v_j)\leq g_m(|i-j|).
\]
For fixed $i$, the sum over $j$ is at most the full-axis sum in \eqref{eq:potential-axis-sum}. Consequently,
\begin{equation}\label{eq:potential-path-energy}
 \nu^TG_{m,U}\nu
 \leq\frac1{mM\sqrt{1+m^2}}\leq\frac1{mM}\leq\frac4K.
\end{equation}
The energy inequality \eqref{eq:potential-energy} proves \eqref{eq:potential-capacity}.
\end{proof}

This estimate is deterministic and uniform over all eligible paths. Folding, revisiting coordinate levels, and excursions of the full baseline path outside $B_N$ cause no loss. It is therefore available after any permitted baseline-measurable selection, without enumerating the possible paths.

\subsubsection{Expected time and hitting probabilities}

\begin{lemma}[Expected time until exit or killing]\label{lem:potential-occupation}
Let $K\geq8$, let $N\geq\max\{8,K\}$ be an integer, and put $m=K/N$. For either $C=C_\pm$ and $U=U_\pm$,
\begin{equation}\label{eq:potential-occupation}
 \frac1{4m^2}\leq (G_{m,U}\mathbf1)(v)\leq\frac1{m^2},
 \qquad v\in C.
\end{equation}
The upper bound holds on all of $U$.
\end{lemma}
\begin{proof}
The core lies a distance at least $N/4$ from the outside of $U$. We use this distance to keep the walk inside $U$ with substantial probability up to the mass-killing time scale $m^{-2}$.
For the rate-one walk started at $v$, the process $|W_t-v|_2^2-t$ is a martingale: each jump has squared length one and zero mean direction. Stopping at $t\wedge\tau_U$ is legitimate since the stopped positions lie in the finite set consisting of $U$ and its nearest neighbors. Hence
\[
 \EE_v|W_{t\wedge\tau_U}-v|_2^2=\EE_v(t\wedge\tau_U)\leq t.
\]
On $\{\tau_U\leq t\}$, the displacement is at least $N/4$ by \eqref{eq:potential-geometry}. Therefore
\[
 \PP_v(\tau_U\leq t)\leq\frac{16t}{N^2}.
\]
For $0\leq t\leq m^{-2}$ and $K\geq8$, this is at most $16/K^2\leq1/4$. Integrating the survival probability against the independent mass clock gives
\[
 (G_{m,U}\mathbf1)(v)
 \geq\frac34\int_0^{m^{-2}}e^{-m^2t}\,dt
 =\frac{3(1-e^{-1})}{4m^2}\geq\frac1{4m^2}.
\]
The upper bound is \eqref{eq:potential-domain}.
\end{proof}

\begin{theorem}[Capacity and hitting estimates for path segments]\label{thm:potential-inputs}
Let $K\geq8$, let $N\geq\max\{8,K\}$ be an integer, and put $m=K/N$. For every nearest-neighbor $\operatorname{TOP}_N$--$\operatorname{BOT}_N$ path in $D=B_{2N}$, form $L_+,L_-$ strictly before first exit from $C_+,C_-$ as above. For either pair $(L,U)=(L_+,U_+)$ or $(L_-,U_-)$,
\begin{equation}\label{eq:potential-main}
\begin{aligned}
 \operatorname{cap}_{m,U}(L)&\geq\frac K4,\qquad
 \inf_{v\in L}(G_{m,U}\mathbf1)(v)\geq\frac1{4m^2},\\
 \PP_{\mu_U}(T_L<\sigma_U)&\geq\frac1{64K}.
 \end{aligned}
\end{equation}
In addition,
\begin{equation}\label{eq:potential-probe-separation}
 \sup_{\substack{x\in U_+\\y\in U_-}}G_{m,D}(x,y)
 \leq\frac4K e^{-K/2},\qquad
 G_{m,D}\mathbf1\leq m^{-2}\mathbf1.
\end{equation}
All estimates hold simultaneously for every eligible path, and hence for the $\mathcal F_0$-measurable selected path on $E_0$ from Section~\ref{sec:probe-extension}.
\end{theorem}
\begin{proof}
The capacity and occupation bounds are Lemmas~\ref{lem:potential-capacity} and~\ref{lem:potential-occupation}. The nonnegative equilibrium identity, area bound, and $mN=K$ give
\begin{align*}
 \PP_{\mu_U}(T_L<\sigma_U)
 &=\frac1{|U|}\sum_{v\in L}e_L(v)(G_{m,U}\mathbf1)(v)\\
 &\geq\frac{\operatorname{cap}_{m,U}(L)}{4m^2|U|}
 \geq\frac{K/4}{4m^2\,4N^2}
 =\frac1{64K}.
\end{align*}
The rectangles are separated by at least $N$. Lemma~\ref{lem:potential-decay}, followed by domain comparison, proves the first bound in \eqref{eq:potential-probe-separation}; the row upper bound is Lemma~\ref{lem:potential-normalization}. The pathwise quantifiers were retained at each step, so the baseline selection introduces no new averaging or counting argument.
\end{proof}

\subsection{Probe attachments and Gaussian comparison}\label{subsec:probe-connections}

We use the capacity and hitting estimates to bound the intensities of loops that join each auxiliary vertex to the selected crossing path.

\subsubsection{The enlarged network and the selected crossing}

Throughout this subsection, let $K\geq8$, let $N\in\mathbb N$ satisfy $N\geq\max\{8,K\}$, and put
\[
 m=K/N,\qquad \mathsf Q_m=I-P_D+m^2I,\qquad
 G=G_{m,D}=\mathsf Q_m^{-1},\qquad \kappa=m^2/2.
\]
Use the rectangles $U_\pm$, cores $C_\pm$, and enlarged precision $\mathsf Q^+$ from Section~\ref{sec:probe-extension}, without changing their lattice endpoints. Old nearest-neighbor conductances are $1/4$. The probe $a$ joins every vertex of $U_+$ with conductance $\kappa$, and $b$ joins every vertex of $U_-$ with that conductance. The old killing is reduced by $\kappa$ there, so every old diagonal remains $1+m^2$. The probe diagonals are
\[
 d_a=\kappa|U_+|,\qquad d_b=\kappa|U_-|.
\]
There is no probe killing and no $a$--$b$ edge. Lemma~\ref{lem:finite-probes} proves positive definiteness and transience of this exact finite network.

The baseline graph $\eta^0$ and enlarged graph $\eta^+$ are coupled by Theorem~\ref{thm:finite-coupling}. In particular, their laws are $\mathsf{Aug}(\mathsf Q_m)$ and $\mathsf{Aug}(\mathsf Q^+)$, respectively, and every baseline edge survives in $\eta^+$. Selection uses exactly
\[
 \mathcal F_0=\sigma\bigl(\mathcal L_0,(\xi_v)_{v\in D},
                 (U_e)_{e\text{ an old edge}}\bigr)
\]
from \eqref{eq:finite-F0}. On the event $E_0$ of a baseline top--bottom connection, the deterministically selected path and its strict pre-exit segments $L_\pm$ are $\mathcal F_0$-measurable. They lie in one baseline component. On $E_0^c$ they are empty, and $\PP(E_0)=\theta_N(m)$.

Retain the deterministic classes $\mathcal C_a,\mathcal C_b$ from Section~\ref{sec:probe-extension}. A loop in $\mathcal C_a$ has exactly one cyclic occurrence of $a$ and all other vertices in $U_+$; the definition for $b$ is reflected. For $L\subseteq U_+$, define
\[
 \Lambda_a(L)=\tfrac12\mu_{\mathsf Q^+}^{\mathrm{nt}}
       \{\ell\in\mathcal C_a:\ell\text{ visits }L\},
\]
and define $\Lambda_b$ analogously. The factor $1/2$ is the soup intensity. The count of cyclic occurrences excludes the repeated terminal endpoint.

\subsubsection{Loop intensity and the killed Green function}

For one probe, write $U$ for its rectangle and $\mathbf1$ for the vector of ones indexed by $U$. Let $W$ be the rate-one lattice walk defined above, let $\sigma_U$ be its first exit from $U$ or independent rate-$m^2$ death, and set
\[
 f_U=G_{m,U}\mathbf1,\qquad T_L=\inf\{t\geq0:W_t\in L\}.
\]
Thus $f_U(x)=\EE_x\sigma_U$ and $T_L=0$ when $x\in L$. Write $\widetilde G_{m,U\setminus L}$ for the Green matrix killed on leaving $U\setminus L$, extended by zero to $U$. If $U\setminus L$ is empty, this extended matrix is zero.

\begin{lemma}[One cyclic probe occurrence has an exact intensity]\label{lem:probes-intensity}
For every $L\subseteq U_+$,
\begin{align}
 \Lambda_a(L)
 &=\frac{\kappa}{2|U_+|}\,
    \mathbf1^T\bigl(G_{m,U_+}
               -\widetilde G_{m,U_+\setminus L}\bigr)\mathbf1
       \label{eq:probes-intensity}\\
 &=\frac\kappa2\,\EE_{\mu_{U_+}}
       \bigl[f_{U_+}(W_{T_L});\,T_L<\sigma_{U_+}\bigr],
       \label{eq:probes-intensity-hit}\\
 0&\leq\Lambda_a(L)\leq\frac14. \label{eq:probes-intensity-upper}
\end{align}
The expectation is defined as zero when $L=\varnothing$. The same statements hold for $b,U_-$.
\end{lemma}
\begin{proof}
We prove the upper-probe case and abbreviate $U=U_+$. The calculation has three steps: root each loop at its unique visit to $a$, sum the resulting words as a Green-function difference, and express that difference using the first hit of $L$.
Begin with the marked cyclic measure in Lemma~\ref{lem:finite-marked}. Each loop in $\mathcal C_a$ has one cyclic position at $a$. Rotating both the vertex word and its holding marks to that position gives a unique rooted representative. A word with a unique $a$ cannot repeat a shorter word. Therefore the sum over its cyclic starting positions cancels exactly the cyclic length factor in \eqref{eq:finite-marked}. Reversal is still not an additional identification.

Write this rooted word as
\[
 a,x_0,x_1,\ldots,x_n,a,\qquad n\geq0,\quad x_i\in U.
\]
It has $n+2$ jumps, one holding interval at $a$, and $n+1$ old holding intervals. Let their lengths be $s_a,s_0,\ldots,s_n$. With the old-edge product interpreted as one when $n=0$, its density after the cyclic cancellation is
\[
 \kappa^2
 \left(\prod_{i=0}^{n-1}P_U(x_i,x_{i+1})\right)
 e^{-d_as_a}\,ds_a
 \prod_{i=0}^{n}e^{-(1+m^2)s_i}\,ds_i.
\]
Here $P_U(x,y)=1/4$ on old nearest-neighbor edges. It is the principal submatrix of the original transition matrix, without another rate normalization. Integrating every holding interval, and setting $q=(1+m^2)^{-1}$, gives
\[
 \frac{\kappa^2}{d_a}\,q^{n+1}
       \prod_{i=0}^{n-1}P_U(x_i,x_{i+1}).
\]
In particular, $n=0$ is the two-jump word $a,x_0,a$ and supplies the leading $q$ term in the resolvent.

All weights are nonnegative. Summing by Tonelli over $n$ and the old vertices gives the kernel
$q\sum_{n\geq0}q^nP_U^n=G_{m,U}$ from Lemma~\ref{lem:potential-normalization}. The words avoiding $L$ are exactly those with all $x_i\in U\setminus L$, whose sum gives the zero-extended killed kernel. Consequently
\[
 \mu_{\mathsf Q^+}^{\mathrm{nt}}\{\ell\in\mathcal C_a:\ell\text{ visits }L\}
 =\frac{\kappa^2}{d_a}\,
       \mathbf1^T\bigl(G_{m,U}
          -\widetilde G_{m,U\setminus L}\bigr)\mathbf1.
\]
The factor $1/2$ enters only now: it is the soup intensity, not an additional factor from choosing the cyclic root.
Multiplication by the soup intensity $1/2$ and the identity $d_a=\kappa|U|$ prove \eqref{eq:probes-intensity}.

For nonempty $L$, decompose the occupation of the killed walk at its first hit of $L$. Before the hit, its occupation is exactly that of the process killed on leaving $U\setminus L$. If the starting point belongs to $L$, this pre-hit occupation is zero, consistently with the zero extension. Strong Markov for the killed process gives, for $x,y\in U$,
\begin{equation}\label{eq:probes-green-hit}
 G_{m,U}(x,y)-\widetilde G_{m,U\setminus L}(x,y)
 =\EE_x\bigl[G_{m,U}(W_{T_L},y);\,T_L<\sigma_U\bigr].
\end{equation}
This also follows from strong Markov for $W$ and the memoryless independent mass clock. Sum over $y$ and average over $x\in U$. Equation~\eqref{eq:probes-intensity} becomes \eqref{eq:probes-intensity-hit}. All expectations are finite because $f_U\leq m^{-2}$ by Lemma~\ref{lem:potential-normalization}. For $L=\varnothing$, both the Green difference and the intensity are zero, which gives the stated convention. For $L=U$, $T_L=0$ and the avoided kernel is zero.

Finally, \eqref{eq:probes-intensity-hit} and the row upper bound imply
\[
 0\leq\Lambda_a(L)\leq\frac\kappa{2m^2}=\frac14.
\]
Reflection gives the lower-probe statements.
\end{proof}

The formula can also be understood by following one excursion from the auxiliary vertex. The first departure from $a$
is uniform on $U_+$, since $\kappa/d_a=1/|U_+|$. Inside this rectangle,
return to $a$ and the remaining ground killing each have rate $\kappa$;
their sum is $m^2$, and exit is an additional failure. The excursion kernel
is therefore $G_{m,U_+}$ and its return probability from $x$ is
$\kappa f_{U_+}(x)$, in agreement with the marked-word factors above.
Excursions that leave $U_+$ and later return belong to other loop classes.

\subsubsection{Conditional attachments to the selected path}

The potential estimates now turn into connection opportunities for the two probes. Their uniformity over every eligible segment is what permits the path to be chosen from the original augmented graph. Conditional on that graph, the two visiting-loop classes remain independent, and the lower bounds on their intensities give a lower bound for connecting both probes to the selected component. Thus the following lemma is the quantitative application of the restriction and equilibrium principles explained above to the specified probe rectangles.

\begin{lemma}[Quantitative attachments and the lower connection bound]\label{lem:probes-attachment}
For every eligible segment $L\subset C_+$,
\begin{equation}\label{eq:probes-attachment-bounds}
 \frac1{1024K}\leq\Lambda_a(L)\leq\frac14,
 \qquad
 1-e^{-\Lambda_a(L)}\geq\frac1{2048K}.
\end{equation}
The analogous estimates hold for $b$ and eligible segments in $C_-$. For the exact enlarged augmented graph,
\begin{equation}\label{eq:probes-connection-lower}
 \PP_{\mathsf Q^+}(a\longleftrightarrow b)
 \geq\frac{\theta_N(K/N)}{2048^2K^2}.
\end{equation}
\end{lemma}
\begin{proof}
Theorem~\ref{thm:potential-inputs} gives, for every eligible $L$ in its rectangle,
\[
 \PP_{\mu_U}(T_L<\sigma_U)\geq\frac1{64K},
 \qquad
 \inf_{v\in L}f_U(v)\geq\frac1{4m^2}.
\]
Use these bounds in \eqref{eq:probes-intensity-hit}, with $\kappa/2=m^2/4$, to get
\[
 \Lambda_a(L)\geq
 \frac{m^2}{4}\,\frac1{64K}\,\frac1{4m^2}
 =\frac1{1024K}.
\]
The upper bound is \eqref{eq:probes-intensity-upper}. If $0\leq t\leq1/4$, then
$e^{-s}\geq1-s\geq3/4$ for $0\leq s\leq t$, and hence
$1-e^{-t}=\int_0^t e^{-s}\,ds\geq t/2$.
This proves \eqref{eq:probes-attachment-bounds}.

The preceding bounds hold for every eligible segment, so on $E_0$ they also hold for the segments selected from the original graph.
Conditional on $\mathcal F_0$, the segments $L_\pm$ are fixed. The deterministic classes $\mathcal C_a,\mathcal C_b$ are disjoint, and their Poisson restrictions are jointly independent of $\mathcal F_0$. Corollary~\ref{cor:finite-attachment} therefore gives independent conditional Poisson counts of the loops hitting $L_+$ and $L_-$, with the respective means $\Lambda_a(L_+)$ and $\Lambda_b(L_-)$.

On $E_0$, both segments are eligible and belong to the same baseline path. If both counts are positive, loop-traversed edges attach $a$ and $b$ to this path, which remains open in $\eta^+$. Thus \eqref{eq:finite-DGFF-transfer} and \eqref{eq:probes-attachment-bounds} give
\[
 \PP_{\mathsf Q^+}(a\longleftrightarrow b)
 \geq\EE\!\left[
       \ind_{E_0}(1-e^{-\Lambda_a(L_+)})
                   (1-e^{-\Lambda_b(L_-)})\right]
 \geq\frac{\PP(E_0)}{2048^2K^2}.
\]
Since $\PP(E_0)=\theta_N(m)$, this proves \eqref{eq:probes-connection-lower}. On $E_0^c$, the dummy sets and both intensities are zero by Lemma~\ref{lem:probes-intensity}; no positive lower bound is asserted there.
\end{proof}

The product is conditional on the baseline, and that is sufficient. It does not assert independence after averaging over a random selected path. There is no conditioning on enlarged occupation, probe failures, or the visiting soup. The other loop classes remain in the full enlarged soup and graph; omitting them from this sufficient connection event does not change either graph marginal.

\subsubsection{The Schur complement and Gaussian correlation}

Eliminating coordinates from a positive quadratic form gives the classical Schur complement. For a network Dirichlet form, this is also its trace on the retained vertices: minimize the energy over the eliminated coordinates while keeping the retained values fixed. Sznitman's trace-network formula gives this interpretation in the finite killed setting; see \cite[author notes, Proposition 1.11]{SznitmanTopics}. Here it identifies the precision of the two-probe Gaussian marginal and makes its correlation computable. The connection event will continue to refer to the full enlarged augmented graph, whose law was constructed above.

Let $u_\pm=\mathbf1_{U_\pm}$ now denote indicator columns indexed by $D$, and put
\[
 H=\kappa\begin{pmatrix}u_+&u_-\end{pmatrix},
 \qquad D_p=\operatorname{diag}(d_a,d_b).
\]
Then \eqref{eq:finite-probe-matrix} is
$\mathsf Q^+=\left(\begin{smallmatrix}\mathsf Q_m&-H\\-H^T&D_p\end{smallmatrix}\right)$.
Write $\Gamma=(\mathsf Q^+)^{-1}$ for the covariance on the entire enlarged vertex set.

\begin{lemma}[Exact Schur matrix and explicit correlation bound]\label{lem:probes-schur}
The covariance block of the two probes is $S^{-1}$, where
\begin{align}
 S&=D_p-H^TGH
   =\begin{pmatrix}S_a&-B\\-B&S_b\end{pmatrix}>0,
       \label{eq:probes-schur}\\
 S_a&=d_a-\kappa^2u_+^TGu_+,\qquad
 S_b=d_b-\kappa^2u_-^TGu_-,\qquad
 B=\kappa^2u_+^TGu_-\geq0. \nonumber
\end{align}
Its entries satisfy
\begin{equation}\label{eq:probes-schur-bounds}
 S_a,S_b\geq\frac{3K^2}{4},
 \qquad
 B\leq16K^3e^{-K/2}.
\end{equation}
The actual Gaussian probe correlation obeys
\begin{equation}\label{eq:probes-rho}
 \rho=\frac{\Gamma_{ab}}{\sqrt{\Gamma_{aa}\Gamma_{bb}}}
      =\frac{B}{\sqrt{S_aS_b}}\in[0,1),
 \qquad
 \rho\leq\frac{64}{3}Ke^{-K/2}.
\end{equation}
\end{lemma}
\begin{proof}
The attachment argument gives a lower bound on the two-probe connection probability. We now prepare the complementary upper bound by expressing their Gaussian correlation through the Schur complement.
For old coordinates $z$ and probe coordinates $w$, complete the square:
\[
 \begin{pmatrix}z\\w\end{pmatrix}^{\!T}
 \mathsf Q^+\begin{pmatrix}z\\w\end{pmatrix}
 =(z-GHw)^T\mathsf Q_m(z-GHw)+w^TSw.
\]
Taking $z=GHw$ shows $S>0$, because $\mathsf Q^+>0$. Solving
$\mathsf Q^+(z,w)^T=(0,v)^T$ first gives $z=GHw$ and then $Sw=v$. Therefore the probe block of $\Gamma$ is exactly $S^{-1}$. Expanding $H^TGH$ gives \eqref{eq:probes-schur}, with $B\geq0$ by entrywise nonnegativity of $G$.

The row bound in Lemma~\ref{lem:potential-normalization} gives
\[
 u_+^TGu_+\leq |U_+|m^{-2}.
\]
Since $\kappa=m^2/2$ and $|U_+|\geq3N^2$,
\[
 S_a\geq d_a-\kappa^2|U_+|m^{-2}
       =d_a/2\geq3K^2/4.
\]
The same calculation gives the lower bound for $S_b$. The separated Green estimate and area bounds from Theorem~\ref{thm:potential-inputs} give
\begin{align*}
 B&\leq\kappa^2|U_+||U_-|
        \sup_{x\in U_+,\,y\in U_-}G(x,y)\\
  &\leq\frac{m^4}{4}(4N^2)^2\frac4K e^{-K/2}
   =16K^3e^{-K/2}.
\end{align*}
Every power of $N$ cancels through $mN=K$.

Put $\Delta=S_aS_b-B^2>0$. Direct inversion yields
\[
 S^{-1}=\Delta^{-1}
             \begin{pmatrix}S_b&B\\B&S_a\end{pmatrix}.
\]
The determinant cancels in the correlation, giving
$\rho=B/\sqrt{S_aS_b}$. Positivity gives $0\leq\rho<1$, and \eqref{eq:probes-schur-bounds} yields its displayed upper bound.
\end{proof}

The matrix $S$ is the precision of the two-probe Gaussian marginal;
$G=\mathsf Q_m^{-1}$ is the inverse of the old principal precision, rather
than the old-coordinate covariance block of $\Gamma$. All connection events
remain in the full graph $\mathsf{Aug}(\mathsf Q^+)$: the Schur calculation
does not assert that deleting old vertices produces a graph with law
$\mathsf{Aug}(S)$.

\subsubsection{The augmented connection identity}

Lupu's geometric coupling identifies each augmented loop cluster with the original vertices of a sign cluster of the metric-graph GFF; see \cite[Proposition 2.1, Theorem 1 bis and Corollary 3.6]{Lupu}. Its two-point consequence is exactly the arcsine identity used here, as obtained for the extended clusters in the proof of \cite[Proposition 5.2]{Lupu}. The raw loop-trace connection event is smaller and receives an inequality in that proposition. Applying the equality to the full augmented network, with the covariance computed by the Schur complement, gives the complementary upper bound on probe connection.

There is a related electrical viewpoint in Lupu and Werner's study of the random pseudo-metric defined by the zero set of the metric-graph GFF~\cite{LWmetric}. They measure a path through the local time accumulated at zero and obtain a two-boundary-point law determined by effective resistance. This illustrates another way in which Gaussian geometry and network reduction interact. Here the Schur calculation is used for the two-probe Gaussian marginal, and the connection event remains the one in the full enlarged augmented graph.

\begin{lemma}[Augmented connectivity and the scalar upper bound]\label{lem:probes-arcsine}
For the enlarged graph specified in Section~\ref{sec:probe-extension},
\begin{equation}\label{eq:probes-arcsine}
 \PP_{\mathsf Q^+}(a\longleftrightarrow b)
 =\frac2\pi\arcsin\rho
 \leq\rho
 \leq\frac{64}{3}Ke^{-K/2},
\end{equation}
where $\rho$ is the actual correlation in \eqref{eq:probes-rho}.
\end{lemma}
\begin{proof}
We apply Lupu's coupling \cite[Theorem 1 bis and Corollary 3.6]{Lupu} and its augmented-cluster connection identity \cite[proof of Proposition 5.2]{Lupu} to the full enlarged network. We verify its hypotheses and normalization here. The old square is connected, and both probe neighborhoods are nonempty, so the enlarged graph is connected. It is finite, simple, and undirected, with positive conductances on its present edges. Every degree is finite. Old killing is at least $m^2/2>0$, probe killing is zero, and $\mathsf Q^+>0$ gives transience. Finite state space and finite rates preclude explosion. The cited coupling requires neither planarity nor a degree bound uniform in $N$.

The generator is $-\mathsf Q^+$ and its occupation Green matrix is $\Gamma$. The half-soup contains all nontrivial marked loops and all trivial occupations, including those of rates $d_a,d_b$ at the probes. Conditional on that soup, the rule \eqref{eq:finite-augmentation} retains traversed edges and opens unused edges independently with probability
$1-\exp(-2c_{xy}\sqrt{X_xX_y})$.
This is the augmented partition in Lupu's coupling. Sharing old uniforms with the baseline does not change its enlarged marginal.

By that construction, assigning independent fair signs to the augmented clusters and multiplying them by $\sqrt{2X}$ gives the Gaussian field with covariance $\Gamma$. Conditional on the partition and occupations, the product of the signs at $a,b$ has expectation one if they share a cluster and zero otherwise. The two Gaussian variances are positive, so neither coordinate is zero with positive probability. Consequently
\[
 \PP_{\mathsf Q^+}(a\longleftrightarrow b)
 =\EE[\operatorname{sign}(\phi_a)\operatorname{sign}(\phi_b)]
 =\frac2\pi\arcsin\rho.
\]
This is the extended-cluster equality in the proof of Proposition 5.2;
the proposition's raw-loop connection statement is an inequality.

For the elementary Gaussian computation, divide each coordinate by its positive standard deviation and represent the pair by
$(Z_1,\rho Z_1+\sqrt{1-\rho^2}Z_2)$, with independent standard Gaussians. Rotational symmetry makes the angle uniform; the sign-disagreement probability is $\arccos(\rho)/\pi$. The expected sign product is therefore $1-2\arccos(\rho)/\pi=(2/\pi)\arcsin\rho$.

Finally, concavity of $\sin t$ on $[0,\pi/2]$ gives $\sin t\geq2t/\pi$. Set $t=\arcsin\rho$ to obtain $(2/\pi)\arcsin\rho\leq\rho$, and then apply Lemma~\ref{lem:probes-schur}.
\end{proof}

The scalar inequality is applied to the actual $\rho\in[0,1)$ before
using $(64/3)Ke^{-K/2}$, which may exceed one near $K=8$. Thus the stated
parameter range requires no larger cutoff.

\subsubsection{Massive connectivity decay}

The preceding steps now meet at a single observable: the connection probability of the two probes. The selected-path argument bounds it from below in terms of the original crossing probability; the classical augmented two-point identity bounds it from above through the Gaussian correlation. The Schur complement and the massive Green estimate make the latter bound explicit. The theorem records the resulting dependence on the mass scale, with the attachment construction and its conditional independence supplying the link to the crossing event.

\begin{theorem}[Uniform massive connectivity decay]\label{thm:probes-M2}
For every real $K\geq8$ and integer $N\geq\max\{8,K\}$, the full-domain augmented crossing probability in \eqref{eq:theta} satisfies
\begin{equation}\label{eq:probes-M2}
 \theta_N(K/N)\leq
 \min\left\{1,\frac{2^{28}}3K^3e^{-K/2}\right\}.
\end{equation}
\end{theorem}
\begin{proof}
Combine Lemmas~\ref{lem:probes-attachment} and~\ref{lem:probes-arcsine}:
\[
 \frac{\theta_N(K/N)}{2048^2K^2}
 \leq\PP_{\mathsf Q^+}(a\longleftrightarrow b)
 \leq\frac{64}{3}Ke^{-K/2}.
\]
Since $2048^2\cdot64/3=2^{28}/3$, this gives the second entry in the minimum. The first entry is the probability bound.
\end{proof}

\begin{corollary}[A uniform massive crossing-height bound]\label{cor:probes-massive-height}
In the same range, with $\psi_{N,K/N}$ as in \eqref{eq:massive},
\begin{equation}\label{eq:probes-massive-height}
 \PP\bigl(H_N(\psi_{N,K/N})>256\bigr)
 \leq\min\left\{1,\frac{2^{28}}3K^3e^{-K/2}\right\}.
\end{equation}
\end{corollary}
\begin{proof}
Theorem~\ref{thm:routes-M1} bounds the left side by $\theta_N(K/N)$. Apply Theorem~\ref{thm:probes-M2}. Theorem~\ref{thm:routes-M1} applies for $N\geq1$, $0<m\leq1$, which contains the present range.
\end{proof}

The factor $4/K$ in the separated Green estimate yields the polynomial factor $K^3$. All constants are explicit sufficient values. The bound may be trivial for some small permitted $K$, but it decays to zero as $K$ increases and is uniform in every integer $N\geq\max\{8,K\}$. The endpoint $m=1$ is included.

\section{Gaussian remainder and the uniform upper tail}\label{sec:gaussian}\label{sec:gaussian-estimates}

We pass from the massive field to the massless field by adding an
independent Gaussian field. An explicit bound on the covariance difference
shows that the maximum of this remainder is of order $\log(2K)$
on $B_N$, with constants uniform in the scale $K=mN$. We then choose $K$
according to the desired crossing height and handle smaller boxes
separately. Throughout, all logarithms are natural.

The decomposition used here belongs to a standard strategy for comparing Gaussian fields: keep the massive component and recover the massless field by adding an independent Gaussian correction. Rodriguez and Zhang use a closely related construction in their study of planar GFF level-set percolation; see \cite[version 2, equations (3.12)--(3.13) and Lemma 3.7]{RZ}. They verify positivity of the covariance difference and control the correction through estimates on its spatial increments. Their argument provides a close precedent for the mass-removal step considered here.

The parameter conventions explain the estimate needed below. Write $M$ for the inverse-mass parameter in the cited argument, whose killing is $M^{-2}$. Lemma~3.7 there bounds a one-sided supremum with probability at least one half when the domain radius $L$ satisfies $2\leq L\leq100M$. Matching the killing to $m=K/N$ gives $M=N/K$: the ratio of our box size to this length scale grows with $K$. We need an estimate for the absolute maximum throughout $1\leq K\leq N$, because $K$ will be chosen according to the desired crossing height. The point and increment estimates below retain that dependence using our own box and Green conventions.

\subsection{Finite Gaussian inputs}

We state the three finite Gaussian inputs with their full parameter ranges
and constants. Their complete proofs are collected in
Appendix~\ref{app:gaussian}. Throughout this section put
\[
 D=B_{2N},\quad \mathsf Q_D=I-P_D,\quad L=4N+2,\quad m=K/N,
 \qquad \ell_K=\log(2K).
\]
Here $L$ denotes the denominator in the finite sine basis, not a path segment. For $1\leq K\leq N$,
write $\zeta_{N,K}$ for a centered Gaussian field with covariance
$\mathsf Q_D^{-1}-(\mathsf Q_D+K^2N^{-2}I)^{-1}$.
The next lemma proves positivity and supplies its independent coupling to
the massive field.

\paragraph{Scalar Gaussian tails.}
For a centered Gaussian $X$ of variance at most $\sigma^2>0$, completing the
square gives $\EE e^{uX}\leq e^{u^2\sigma^2/2}$. Markov's inequality with
$u=v/\sigma^2$, applied to $X$ and $-X$, proves
\begin{equation}\label{eq:gaussian-scalar-tail}
 \PP(|X|>v)\leq2e^{-v^2/(2\sigma^2)},\qquad v\geq0.
\end{equation}
A zero-variance variable causes no failure. We will use this estimate both
in the finite grid argument and in the direct bound for small boxes.

\paragraph{Spectrum and covariance.}
The sine basis for a rectangle is standard; see
\cite[author version, Section 8.2]{LawlerLimic}. The next lemma fixes its
normalization and the covariance decomposition, including the smallest box.

The positivity of the covariance difference and the independent Gaussian sum are elementary consequences of diagonalizing the two resolvents in the same basis. We include this calculation to identify the correction precisely and to retain the exact massive covariance required by the loop-soup comparison. The subsequent estimates are for this correction field, with its full mass dependence. This is the covariance-splitting step of the established Gaussian comparison strategy described above.

\begin{lemma}[Full-range finite spectrum and covariance]\label{lem:gaussian-spectrum}
For every integer $N\geq1$ and every $m>0$, the functions
\[
 e_{k\ell}(x)=\frac2L
 \sin\frac{\pi k(x_1+2N+1)}L\,
 \sin\frac{\pi\ell(x_2+2N+1)}L,
 \qquad 1\leq k,\ell\leq L-1,
\]
form a counting-measure orthonormal eigenbasis of $\mathsf Q_D$, with
\[
 \lambda_{k\ell}
 =\sin^2\frac{\pi k}{2L}+\sin^2\frac{\pi\ell}{2L},
 \qquad
 \frac{k^2+\ell^2}{L^2}\leq\lambda_{k\ell}
 \leq\frac{\pi^2(k^2+\ell^2)}{4L^2}.
\]
The covariance $R=\mathsf Q_D^{-1}-(\mathsf Q_D+m^2I)^{-1}$ is positive definite, with eigenvalues
\begin{equation}\label{eq:gaussian-coefficients}
 a_{k\ell}=\frac{m^2}{\lambda_{k\ell}(\lambda_{k\ell}+m^2)}.
\end{equation}
There are independent centered Gaussian fields $\psi$ and $\zeta$ with covariances $(\mathsf Q_D+m^2I)^{-1}$ and $R$, respectively, whose sum has exactly the original massless law.
\end{lemma}

The exact finite-mode normalization and independent Gaussian construction are proved in Appendix~\ref{app:gaussian-spectrum}.

\paragraph{Point and increment estimates.}

To control a maximum, we need an estimate at each point and an estimate for spatial differences. Both follow from the same full spectral coefficient; the second estimate also retains the oscillation of each sine product.

\begin{lemma}[Uniform logarithmic variance and short-distance increments]\label{lem:gaussian-increments}
For every integer $N\geq1$, real $1\leq K\leq N$, and $x\in D$,
\begin{equation}\label{eq:gaussian-point}
 \Var\zeta_{N,K}(x)\leq64\ell_K.
\end{equation}
For $x,y\in D$, put $\rho=\|x-y\|_\infty/N$. Whenever $0<K\rho\leq1$,
\begin{equation}\label{eq:gaussian-increment}
 \EE\bigl(\zeta_{N,K}(x)-\zeta_{N,K}(y)\bigr)^2
 \leq2048(K\rho)^2\log\frac e{K\rho}.
\end{equation}
When $x=y$, the increment is zero. These constants do not depend on $K$ or $N$.
\end{lemma}

The frequency sums, including their real cutoffs and endpoint cases, are proved in Appendix~\ref{app:gaussian-increments}.

\paragraph{A grid at the mass scale.}

The passage from point-variance and spatial-increment bounds to an expected absolute maximum follows Dudley's classical chaining method~\cite{Dudley}. Borell--TIS concentration then controls deviations from that expectation~\cite{Borell,TIS}. For the formulations used here, see \cite[Theorem 6.8 and Corollary 6.6]{Biskup} and \cite[Theorem 8.1.3 and Remark 8.1.6]{Vershynin}. The point-variance and increment estimates together determine the size of the expected maximum, while the maximal pointwise variance also determines its Gaussian deviation scale.

We use finite integer-grid chaining and Gaussian concentration, standard methods for Gaussian maxima;
see, for example, \cite[Theorem 8.1.3 and Remark 8.1.6]{Vershynin} for the
general principle. The complete integer-grid construction and the concentration calculation
are given in Appendix~\ref{app:gaussian-net}. The short-distance estimate applies when $K\rho\leq1$,
so the initial integer mesh is no larger than $N/K$. This includes the case
in which the initial grid already contains every lattice point.

\begin{lemma}[Numerical finite Gaussian bound]\label{lem:gaussian-net}
Let $N\geq1$ be an integer, let $1\leq K\leq N$ be real, and let
$(Z_x)_{x\in B_N}$ be a centered jointly Gaussian vector. Suppose
\[
 \Var Z_x\leq64\ell_K,
 \qquad
 \Var(Z_x-Z_y)\leq2048(K\rho)^2\log\frac e{K\rho}
 \quad\text{if }0<K\rho\leq1,
 \qquad \rho=\|x-y\|_\infty/N.
\]
Then for every $t\geq0$,
\begin{equation}\label{eq:gaussian-net-tail}
 \PP\!\left(\max_{x\in B_N}|Z_x|
 >800\{\ell_K+t\sqrt{\ell_K}\}\right)
 \leq\frac{25}{126}e^{-8t^2}\leq2e^{-t^2}.
\end{equation}
The field values and increments need not be independent.
\end{lemma}

The proof in Appendix~\ref{app:gaussian-net} combines the exact integer
parent grids with a bound for the expected maximum and the classical
Borell--TIS concentration inequality~\cite[Corollary 6.6]{Biskup}.
The numerical specialization is the one needed for Theorem~\ref{thm:gaussian-remainder}.

\subsection{The remainder bound and the crossing estimate}

\begin{theorem}[Uniform logarithmic remainder]\label{thm:gaussian-remainder}
For every integer $N\geq1$, real $1\leq K\leq N$, and $t\geq0$, the Gaussian remainder with covariance
\[
 R_{N,K}=\mathsf Q_D^{-1}-(\mathsf Q_D+K^2N^{-2}I)^{-1}
\]
satisfies
\begin{equation}\label{eq:gaussian-remainder-tail}
 \PP\!\left(\|\zeta_{N,K}\|_{\infty,B_N}
 >800\{\ell_K+t\sqrt{\ell_K}\}\right)\leq2e^{-t^2},
 \qquad \ell_K=\log(2K).
\end{equation}
In fact the right side may be replaced by $(25/126)e^{-8t^2}$.
\end{theorem}

\begin{proof}
Lemma~\ref{lem:gaussian-spectrum} constructs the remainder throughout the stated range. Lemma~\ref{lem:gaussian-increments} supplies exactly the variance constant $64$ and increment constant $2048$ required by Lemma~\ref{lem:gaussian-net}. Its numerical bound yields both assertions, uniformly in $K$ and $N$.\end{proof}

\begin{corollary}[A two-parameter crossing estimate]\label{cor:gaussian-massless}
For every real $K\geq8$, integer $N\geq K$, and $t\geq0$,
\begin{equation}\label{eq:gaussian-massless-bound}
 \PP\!\left(H_N>256+800\{\ell_K+t\sqrt{\ell_K}\}\right)
 \leq\min\left\{1,\frac{2^{28}}3K^3e^{-K/2}+2e^{-t^2}\right\}.
\end{equation}
\end{corollary}

\begin{proof}
Use Lemma~\ref{lem:gaussian-spectrum} to realize the original field as the sum of its massive marginal and an independent remainder. The same path can be used for the sum, with its height increased by at most the maximum absolute value of the remainder. The deterministic inequality \eqref{eq:deterministic-H} contains the event on the left in the union of the massive-height event above $256$ and the remainder event in \eqref{eq:gaussian-remainder-tail}. Corollary~\ref{cor:probes-massive-height}, Theorem~\ref{thm:gaussian-remainder}, and the union bound give \eqref{eq:gaussian-massless-bound}. The common range has $N\geq8$, so all earlier massive hypotheses hold. No new independence condition on a selected crossing is required.
\end{proof}

We next choose $K$ and $t$ in terms of the height $A$. Boxes with $N<K$
are estimated directly, so the resulting tail bound will hold for every
positive integer $N$.

\paragraph{Constants for the upper tail.}
Recall $F(a)=a^2/\log(2+a)$ for $a\geq0$. Put
\begin{equation}\label{eq:tails-constants}
 A_*=2^{16},\qquad c_*=2^{-29}.
\end{equation}
These values are sufficient for the estimates below; no optimality is asserted.

\subsection{A finite-box estimate with no size exception}

The massive comparison in Corollary~\ref{cor:gaussian-massless} requires $N\geq K$. To cover the other boxes, we first estimate the original field directly.

\begin{lemma}[Diagonal Green bound for every box]\label{lem:tails-green}
For every integer $N\geq1$,
\begin{equation}\label{eq:tails-green}
 \max_{x\in D}G_D(x,x)\leq32\log(N+2).
\end{equation}
Consequently, for every $A\geq0$,
\begin{equation}\label{eq:tails-box-maximum}
 \PP(H_N>A)
 \leq2(2N+1)^2
       \exp\!\left\{-\frac{A^2}{64\log(N+2)}\right\}.
\end{equation}
\end{lemma}
\begin{proof}
Use the full-range sine basis of Lemma~\ref{lem:gaussian-spectrum}, with $L=4N+2$ and $M=L-1=4N+1$. Its amplitude and eigenvalue bounds give
\[
 G_D(x,x)
 =\sum_{k,\ell=1}^M\frac{e_{k\ell}(x)^2}{\lambda_{k\ell}}
 \leq4\sum_{k,\ell=1}^M\frac1{k^2+\ell^2}
 \leq8\sum_{n=1}^M\frac1n
 \leq8\{1+\log(4N+1)\}.
\]
The middle estimate groups pairs by $n=\max(k,\ell)$, with $2n-1\leq2n$ pairs in each shell. Since $4N+1\leq(N+2)^2$ and $\log(N+2)\geq\log2\geq1/2$, the final expression is at most $32\log(N+2)$. This proves~\eqref{eq:tails-green}, including $N=1$.

The scalar estimate~\eqref{eq:gaussian-scalar-tail} gives
\[
 \PP(|\phi_N(x)|>A)
 \leq2\exp\!\left\{-\frac{A^2}{64\log(N+2)}\right\}.
\]
There is a nearest-neighbor left-right path in every $B_N$, so deterministically $H_N\leq\max_{B_N}|\phi_N|$. A union bound over its $(2N+1)^2$ vertices proves~\eqref{eq:tails-box-maximum}. No independence between vertices is used.
\end{proof}

\subsection{The uniform upper tail}\label{subsec:uniform-tails}

\begin{theorem}[Uniform upper tail]\label{thm:tails-upper}
For every real $A\geq0$,
\begin{equation}\label{eq:tails-upper}
 \sup_{N\geq1}\PP(H_N>A)\leq4e^{-c_*F(A)},
 \qquad c_*=2^{-29}.
\end{equation}
For $A\geq A_*=2^{16}$, the following stronger estimate also holds:
\begin{equation}\label{eq:tails-upper-large}
 \sup_{N\geq1}\PP(H_N>A)
 \leq3\exp\!\left\{-\frac{F(A)}{7680000}\right\}.
\end{equation}
The supremum is over all positive integers $N$.
\end{theorem}
\begin{proof}
We choose $K$ to control the massive-height probability and $t$ to pay the Gaussian-remainder fluctuation within height $A$. The resulting estimate applies directly for $N\geq K$; smaller boxes use the finite-box Gaussian bound instead.

\emph{The height budget.}
Fix real $A\geq A_*$ and set
\[
 K=\lceil A^2\rceil,\qquad \ell=\log(2K),
 \qquad t=\frac{A}{1600\sqrt\ell}.
\]
The choice of $t$ spends exactly $800t\sqrt\ell=A/2$ of the height threshold. It remains to show that $256+800\ell$ costs at most the other half.
Then $A^2\leq K\leq2A^2$, $K\geq8$, and, since $A\geq4$,
\begin{equation}\label{eq:tails-log-budget}
 \ell\leq\log(4A^2)\leq3\log A\leq3\log(2+A).
\end{equation}
The function $A/2-256-2400\log A$ has positive derivative on $[A_*,\infty)$. At $A_*$ its value is at least
\[
 32768-256-2400\cdot16\cdot\frac34=3712>0,
\]
using $\log2\leq3/4$. Hence
\begin{equation}\label{eq:tails-height-budget}
 256+800\ell\leq A/2,\qquad
 256+800\{\ell+t\sqrt\ell\}\leq A.
\end{equation}
In particular $t>0$, and the strict event $\{H_N>A\}$ is contained in the event at the comparison threshold, also when the thresholds coincide.

\emph{Large boxes: $N\geq K$.}
All hypotheses of Corollary~\ref{cor:gaussian-massless} hold. Equations~\eqref{eq:tails-log-budget}--\eqref{eq:tails-height-budget} give
\begin{equation}\label{eq:tails-large-before-absorb}
 \PP(H_N>A)
 \leq\frac{2^{28}}3K^3e^{-K/2}
       +2\exp\!\left\{-\frac{F(A)}{7680000}\right\}.
\end{equation}
To pay the first prefactor explicitly, for real $K\geq256$ the function
\[
 q(K)=K/4-28\log2-3\log K
\]
has positive derivative. At $K=256$ it is $64-52\log2\geq25$. Thus
\[
 \frac{2^{28}}3K^3e^{-K/2}\leq e^{-K/4}\leq e^{-A^2/4}.
\]
Our chosen $K\geq2^{32}$ lies in this range. Since $\log(2+A)\geq1/2$, one has $A^2/4\geq F(A)/7680000$. Therefore~\eqref{eq:tails-large-before-absorb} is at most the right side of~\eqref{eq:tails-upper-large}.

\emph{Small boxes: $N<K$.}
The integer condition gives
\[
 N\leq\lceil A^2\rceil-1<A^2.
\]
This remains true when $A^2$ is an integer. Consequently
\[
 (2N+1)^2\leq9A^4,\qquad
 \log(N+2)\leq\log(A^2+2)\leq2\log(A+2).
\]
Lemma~\ref{lem:tails-green} now gives
\begin{equation}\label{eq:tails-small-before-absorb}
 \PP(H_N>A)\leq18A^4e^{-F(A)/128}.
\end{equation}
For $A\geq16$, $\sqrt A-\log A$ is increasing and has value at least $1$ at $16$. Hence $\log A\leq\sqrt A$. For $A\geq2$,
\[
 \log(A+2)\leq\log(2A)\leq2\log A.
\]
Since our $A\geq A_*\geq2560$ and $A_*\geq18$, these estimates imply
\[
 \log18+4\log A
 \leq5\sqrt A
 \leq\frac{A^{3/2}}{512}
 \leq\frac{F(A)}{256}.
\]
Thus~\eqref{eq:tails-small-before-absorb} is at most $e^{-F(A)/256}$, which is no larger than the right side of~\eqref{eq:tails-upper-large}. The two box ranges cover every positive integer $N$, proving that estimate.

\emph{All remaining heights.}
The function $F$ is continuous with $F(0)=0$, and for $a>0$,
\begin{equation}\label{eq:tails-F-derivative}
 F'(a)=
 \frac{a\{2\log(2+a)-a/(2+a)\}}{\{\log(2+a)\}^2}>0.
\end{equation}
Indeed $2\log(2+a)\geq2\log2\geq1>a/(2+a)$. Thus $F$ is strictly increasing. Moreover
\[
 F(A_*)\leq\frac{2^{32}}{16\log2}\leq2^{29}.
\]
For $0\leq A\leq A_*$, we have $c_*F(A)\leq1$ and
\[
 1\leq4e^{-1}\leq4e^{-c_*F(A)},
\]
where $\log4=2\log2\geq1$ proves the first inequality. The probability bound by one therefore gives~\eqref{eq:tails-upper} on this entire interval. For $A\geq A_*$, use $7680000<2^{29}$ in~\eqref{eq:tails-upper-large}. This proves the assertion for every $A\geq0$ without exchanging limits.
\end{proof}

\subsection{Explicit uniform moment consequences}

The upper-tail estimate holds for every box size and height. We can therefore obtain the moment bounds by integrating this tail.

\begin{corollary}[Uniform exponential and polynomial moments]\label{cor:tails-moments}
For every $0<\eta<c_*=2^{-29}$,
\begin{equation}\label{eq:tails-moment-general}
 \sup_{N\geq1}\EE e^{\eta F(H_N)}
 \leq1+\frac{4\eta}{c_*-\eta}.
\end{equation}
In particular,
\begin{equation}\label{eq:tails-moment-five}
 \sup_{N\geq1}\EE\exp\!\left\{
 2^{-30}\frac{H_N^2}{\log(2+H_N)}\right\}\leq5.
\end{equation}
For every real $p>0$,
\begin{equation}\label{eq:tails-polynomial}
 \sup_{N\geq1}\EE H_N^p
 \leq2^p+4\Gamma(p+1)2^{29p}<\infty.
\end{equation}
Here $\Gamma(p)=\int_0^\infty u^{p-1}e^{-u}\,du$ is the ordinary gamma function.
\end{corollary}
\begin{proof}
Equation~\eqref{eq:tails-F-derivative} proves strict monotonicity of $F$. Also $a-\log(a+2)$ is increasing for $a\geq2$, and at $a=2$ it is at least $2-2(3/4)=1/2$. Hence
\[
 F(a)\geq a\quad(a\geq2).
\]
Thus $F$ is continuous, strictly increasing, unbounded, and maps $[0,\infty)$ onto itself. For each $u\geq0$ let $a$ be its unique inverse value. Theorem~\ref{thm:tails-upper} gives the exact transformed-tail bound
\[
 \PP(F(H_N)>u)=\PP(H_N>a)\leq4e^{-c_*u}.
\]
The original finite Gaussian vector and the finite crossing minimum make $H_N$ finite and nonnegative.
We integrate the transformed tail for the exponential moment. For the polynomial moment, we return to the original height variable and split its tail integral at $a=2$.

For $Y=F(H_N)$, the pointwise identity
\[
 e^{\eta Y}
 =1+\int_0^\infty\eta e^{\eta u}\ind_{\{Y>u\}}\,du
\]
and Tonelli's theorem give
\[
 \EE e^{\eta Y}
 \leq1+4\eta\int_0^\infty e^{-(c_*-\eta)u}\,du
 =1+\frac{4\eta}{c_*-\eta}.
\]
All integrands are nonnegative, so no prior integrability assumption is used. Taking $\eta=c_*/2=2^{-30}$ gives~\eqref{eq:tails-moment-five}.

For $p>0$, the same nonnegative tail integration and $F(a)\geq a$ for $a\geq2$ give
\begin{align*}
 \EE H_N^p
 &=p\int_0^\infty a^{p-1}\PP(H_N>a)\,da\\
 &\leq2^p+4p\int_0^\infty a^{p-1}e^{-c_*a}\,da\\
 &=2^p+4\Gamma(p+1)c_*^{-p}.
\end{align*}
The gamma integral is finite for $p>0$, and $\Gamma(p+1)=p\Gamma(p)$ follows by integration by parts. This proves~\eqref{eq:tails-polynomial} uniformly in $N$.
\end{proof}

\section{A uniform lower bound for crossing failure}\label{sec:lower}

\subsection{A Gaussian translation and a positive crossing at level zero}

This lower bound uses the massless Gaussian field directly; it does not depend on the massive-field comparison. A positive star crossing gives a barrier to a crossing in the two-sided band. We then shift the field by a deterministic function, using a Gaussian translation inequality and a bound on its Dirichlet energy. Recall that $\Phi$ is the standard normal distribution function. The next lemma is the finite-dimensional Gaussian shift inequality of Li and Kuelbs~\cite{KuelbsLi}; its likelihood-ratio proof goes back to the Neyman--Pearson principle~\cite[Section III(a)]{NeymanPearson}. We include a short proof to fix the sign of the shift and the exact constant.

More precisely, \cite[Theorem 1, equation (1.1)]{KuelbsLi} gives the inequality below after whitening. For a Gaussian vector with precision $\mathsf Q$, the linear change of coordinates $x\mapsto\mathsf Q^{1/2}x$ makes the covariance the identity. It takes $E$ to $\mathsf Q^{1/2}E$ and the translation vector to $\mathsf Q^{1/2}h$, whose Euclidean length is $\sqrt{h^{\mathsf T}\mathsf Qh}$. This is exactly the shift cost appearing in the lemma.

Li and Kuelbs' result applies to measurable sets, so no convexity condition on the crossing event is needed. We retain the proof of this known inequality to show the halfspace comparison directly and fix the direction of translation used for the positive barrier. The later work is to choose that barrier and a deterministic shift with bounded Dirichlet energy.

\begin{lemma}[Finite Gaussian translation]\label{lem:tails-translation}
Let $\mu$ be a centered nondegenerate finite-dimensional Gaussian law with precision $\mathsf Q$. For a measurable set $E$ and a deterministic vector $h$,
\begin{equation}\label{eq:tails-translation}
 \mu(E+h)\geq
 \Phi\!\left(\Phi^{-1}(\mu(E))-\sqrt{h^{\mathsf T}\mathsf Qh}\right).
\end{equation}
The endpoint probabilities are interpreted by the extended normal quantiles.
\end{lemma}
\begin{proof}
Suppose first that $0<p=\mu(E)<1$ and $q=\sqrt{h^{\mathsf T}\mathsf Qh}>0$. For $X$ with law $\mu$, the scalar
\[
 Z=\frac{h^{\mathsf T}\mathsf QX}{q}
\]
is standard normal. Changing variables in the Gaussian density gives
\[
 \mu(E+h)=\EE\!\left[\ind_E(X)e^{-qZ-q^2/2}\right].
\]
Because the weight decreases with $Z$, the natural candidate for the least weighted mass among events of probability $p$ is an upper normal tail. The pointwise comparison below verifies that this choice is optimal.
Set $a=\Phi^{-1}(1-p)$, $T=\{Z\geq a\}$, and $w(z)=e^{-qz-q^2/2}$. Then $\PP(T)=p$, $w$ is decreasing, and pointwise
\[
 (\ind_E-\ind_T)\{w(Z)-w(a)\}\geq0.
\]
All terms are integrable, and $\EE(\ind_E-\ind_T)=0$. Integration followed by completing the square therefore gives
\begin{align*}
 \mu(E+h)
 &\geq\EE[\ind_Tw(Z)]
 =\int_a^\infty\frac{e^{-(z+q)^2/2}}{\sqrt{2\pi}}\,dz\\
 &=\Phi(-a-q)=\Phi(\Phi^{-1}(p)-q).
\end{align*}
If $q=0$, positive definiteness implies $h=0$ and the result is equality. If $p=0$ or $p=1$, equivalence of the strictly positive Gaussian densities and their translates proves the endpoint statement.
\end{proof}

\begin{lemma}[A positive star barrier]\label{lem:tails-positive}
For every integer $N\geq1$, the event $E_N$ that $\phi_N$ has a strictly positive star top-bottom crossing in $B_N$ satisfies
\begin{equation}\label{eq:tails-positive}
 \PP(E_N)\geq\frac12.
\end{equation}
For every deterministic $f=1$ on $B_N$ and every $A\geq0$,
\begin{equation}\label{eq:tails-shift-inclusion}
 E_N+Af\subseteq\{H_N>A\}.
\end{equation}
\end{lemma}
\begin{proof}
Apply Lemma~\ref{lem:routes-duality} to the binary field $g(x)=\ind_{\{\phi_N(x)>0\}}$ at threshold zero. Exactly one of a nonpositive nearest-neighbor left-right crossing for $\phi_N$ and a positive star top-bottom crossing occurs.

Every coordinate has positive variance because $\mathsf Q_D^{-1}$ is positive definite. The finite set $B_N$ therefore has no zero field value almost surely. Sign symmetry and a quarter-turn of the square show that the nonpositive nearest-neighbor left-right crossing has the same probability as a strictly positive nearest-neighbor top-bottom crossing. The rotation preserves the law because it preserves $D$ and the kernel with nearest-neighbor entries $1/4$. Writing these last two positive-crossing probabilities as $p_{\mathrm{star}}$ and $p_{\mathrm{NN}}$, the complementarity gives
\[
 p_{\mathrm{star}}+p_{\mathrm{NN}}=1.
\]
Since every nearest-neighbor path is a star path, $p_{\mathrm{NN}}\leq p_{\mathrm{star}}$. This proves~\eqref{eq:tails-positive}. The event is measurable: it is a finite union over simple star paths.

If $\phi\in E_N$, its positive star crossing has
\[
 (\phi+Af)(x)=\phi(x)+A>A
\]
at every vertex of that crossing. The same deterministic duality lemma gives~\eqref{eq:tails-shift-inclusion}. This inclusion is strict pointwise and includes $A=0$; no threshold perturbation is needed.
\end{proof}

\subsection{Square-shell capacity and the lower-tail bound}

The use of a deterministic field shift with a Dirichlet-energy cost is also central to Sznitman's disconnection lower bound for the GFF in dimension $d\geq3$; see \cite[Theorem 2.1 and its proof, equations (2.3)--(2.12)]{SznitmanDisc}. There the shifted field makes a one-sided disconnection event likely, and the resulting asymptotic cost is expressed through capacity. Here planar duality has already supplied a positive barrier with probability at least one half. The exact finite Gaussian shift inequality converts a bound on the energy of a function equal to one on $B_N$ into the stated lower bound for failure of a two-sided crossing. This explains the role of the condenser capacity introduced next.

Extend a vector $f\in\mathbb R^D$ by zero outside $D$. In the original normalization,
\begin{equation}\label{eq:tails-energy}
 f^{\mathsf T}\mathsf Q_Df
 =\frac14\sum_{\{u,v\}:u\sim v}(f(u)-f(v))^2.
\end{equation}
The sum counts each unordered nearest-neighbor lattice edge once. Edges outside $D$ contribute zero, and edges leaving $D$ are included. Expanding the squares proves the identity: each vertex has four incident edges, and the internal cross terms are precisely those of $-P_D$.

Define the massless condenser capacity by
\begin{equation}\label{eq:tails-capacity}
 \operatorname{cap}_D(B_N)
 =
 \min\{f^{\mathsf T}\mathsf Q_Df:\ f\in\mathbb R^D,\ f=1\text{ on }B_N\}.
\end{equation}
Positive definiteness makes the quadratic form coercive and strictly convex, so the nonempty closed affine constraint set has a unique minimizer. This convention also equals the infimum with $f\geq1$ on $B_N$: clipping to $[0,1]$ fixes the exterior value zero, does not increase any term in~\eqref{eq:tails-energy}, and gives equality to one on $B_N$. This is the finite-network Dirichlet variational principle; see \cite[Sections 2.4--2.5]{LyonsPeres}. The shell bound below also follows directly from the Nash--Williams disjoint-cutset inequality~\cite[Section 2.5, equation (2.13)]{LyonsPeres}.

The network ingredient is Nash--Williams' disjoint-cutset resistance inequality~\cite{NashWilliams}, in the formulation recalled in \cite[Section 2.5, equation (2.13)]{LyonsPeres}. The square shells furnish an explicit sequence of cutsets, including the last cutset to the killing exterior. Their resistance lower bound becomes an upper bound on capacity. The shell-constant test functions also attain the corresponding energy within that restricted class. This gives the energy upper bound needed here; it does not evaluate the unrestricted condenser capacity.

\begin{lemma}[Square-shell test in the killed domain]\label{lem:tails-shell}
For every integer $N\geq1$, put
\[
 R_N=\sum_{j=N}^{2N}\frac1{2j+1},
 \qquad E_N^{\square}=R_N^{-1}.
\]
Then
\begin{equation}\label{eq:tails-shell-capacity}
 \operatorname{cap}_D(B_N)\leq E_N^{\square}<\frac2{\log2}.
\end{equation}
The value $E_N^{\square}$ is the minimum energy among tests constant on each square shell, equal to one on $B_N$, and zero outside $D$. Moreover
\begin{equation}\label{eq:tails-shell-limit}
 0<\frac2{\log2}-E_N^{\square}
 <\frac4{(\log2)^2(2N+1)},
 \qquad E_N^{\square}\longrightarrow\frac2{\log2}.
\end{equation}
These statements do not identify the full condenser capacity or its limit.
\end{lemma}
\begin{proof}
Wire $B_N$ to a source and $D^c$ to a sink in the nearest-neighbor
network of conductance $1/4$. For each $j=N,\ldots,2N$, the edges
connecting the square shells of radii $j$ and $j+1$ form a cutset. There
are exactly $4(2j+1)$ such edges, so their total conductance is
$2j+1$. The cutsets are pairwise edge-disjoint, including the final
cutset from radius $2N$ to the killing exterior. Nash--Williams'
inequality~\cite[Section 2.5, equation (2.13)]{LyonsPeres} gives
\[
 R_{\rm eff}(B_N\leftrightarrow D^c)
 \geq\sum_{j=N}^{2N}\frac1{2j+1}=R_N.
\]
Effective conductance equals the finite Dirichlet capacity in
\eqref{eq:tails-capacity}. Hence
$\operatorname{cap}_D(B_N)\leq R_N^{-1}$.

We also compute the minimum energy among functions that are constant on each square shell.
If an admissible function is constant on each shell, let $\delta_j$
be its successive drops, $j=N,\ldots,2N$, including the final drop
to zero outside $D$. Their sum is one, and~\eqref{eq:tails-energy}
gives the energy $\sum_{j=N}^{2N}(2j+1)\delta_j^2$.
Cauchy--Schwarz shows that this is at least $R_N^{-1}$, with equality
at $\delta_j=1/[R_N(2j+1)]$. These drops define a shell-constant
admissible function, so its minimum energy is exactly
$E_N^{\square}=R_N^{-1}$. This equality is only within the restricted
shell-constant class; no equality with the unrestricted capacity is
claimed.

Finally, the decreasing function $(2x+1)^{-1}$ gives
\[
 \frac12\log2<R_N
 <\frac12\log2+\frac1{2N+1}.
\]
Writing $a=(\log2)/2$, we obtain
\[
 0<\frac1a-\frac1{R_N}
 =\frac{R_N-a}{aR_N}
 <\frac1{a^2(2N+1)}.
\]
This proves both bounds. In particular $N=1$ gives
$R_1=8/15$ and $E_1^{\square}=15/8$.
\end{proof}

An upper bound on this capacity is enough: it bounds the energy of an admissible deterministic shift and hence the cost in the Gaussian translation inequality. Combining that shift with the positive star barrier gives the following lower tail. This argument uses the massless field, planar duality, and the classical Gaussian and network inputs just described, independently of the massive crossing estimates.

\begin{theorem}[Independent uniform Gaussian lower tail]\label{thm:tails-lower}
For every integer $N\geq1$ and real $A\geq0$,
\begin{equation}\label{eq:tails-lower-capacity}
 \PP(H_N>A)
 \geq\Phi\!\left(-A\sqrt{\operatorname{cap}_D(B_N)}\right)
 \geq\Phi(-A\sqrt{E_N^{\square}})
 \geq\Phi\!\left(-A\sqrt{\frac2{\log2}}\right).
\end{equation}
In particular,
\begin{equation}\label{eq:tails-lower}
 \inf_{N\geq1}\PP(H_N>A)
 \geq\Phi\!\left(-A\sqrt{\frac2{\log2}}\right)>0,
 \qquad A\geq0.
\end{equation}

\end{theorem}
\begin{proof}
Let $f$ be the minimizer in~\eqref{eq:tails-capacity}. For $h=Af$, the squared translation cost is $h^{\mathsf T}\mathsf Q_Dh=A^2\operatorname{cap}_D(B_N)$. Also, Lemma~\ref{lem:tails-positive} gives $\Phi^{-1}(\PP(E_N))\geq0$. These two facts explain the last inequality in the calculation below.
Lemmas~\ref{lem:tails-positive} and~\ref{lem:tails-translation}, with $h=Af$, give
\begin{align*}
 \PP(H_N>A)
 &\geq\PP(\phi_N\in E_N+Af)\\
 &\geq
 \Phi\!\left(\Phi^{-1}(\PP(E_N))
              -A\sqrt{f^{\mathsf T}\mathsf Q_Df}\right)\\
 &\geq\Phi\!\left(-A\sqrt{\operatorname{cap}_D(B_N)}\right).
\end{align*}
The event $E_N$ has probability strictly between zero and one, since both the all-positive and the all-negative open orthants have positive Gaussian density mass. The displayed calculation also covers $A=0$ by the zero-shift case. Lemma~\ref{lem:tails-shell} and monotonicity of $\Phi$ give the remaining inequalities. Taking the infimum over all positive integers $N$ proves~\eqref{eq:tails-lower}.
\end{proof}

\begin{proof}[Completion of the main theorem]
Theorem~\ref{thm:tails-upper} proves the upper inequality in
Theorem~\ref{thm:main} for every height and every box size.
Theorem~\ref{thm:tails-lower} proves its lower inequality under the same
field normalization. Corollary~\ref{cor:tightness} is the direct
consequence already proved in the introduction, and
Corollary~\ref{cor:tails-moments} proves the moment bounds
stated in the introduction.
\end{proof}

\appendix

\section{Finite Gaussian calculations}\label{app:gaussian}

This appendix gives the complete proofs of the three finite Gaussian inputs
stated in Section~\ref{sec:gaussian}. We use $D=B_{2N}$,
$\mathsf Q_D=I-P_D$, and $L=4N+2$. The sine basis $e_{k\ell}$, eigenvalues
$\lambda_{k\ell}$, and covariance coefficients $a_{k\ell}$ are those of
Lemma~\ref{lem:gaussian-spectrum}. In the point, increment, and grid estimates,
$1\leq K\leq N$ is real, $m=K/N$, and $\ell_K=\log(2K)$.
All Gaussian fields and grids used below are finite. We use infinite frequency sums only to obtain upper bounds for the finite covariance sums.

The appendix follows the same three steps as Section~\ref{sec:gaussian}: the standard rectangular eigenbasis fixes the covariance, the frequency calculation bounds points and increments, and the Gaussian-process argument controls the maximum. The last step is an explicit finite application of Dudley chaining and Borell--TIS concentration, with the modern references given there. Retaining these calculations makes the numerical bounds and their dependence on the real parameter $K$ directly checkable, including the smallest box and the largest permitted mass. The appendix is intended as a detailed implementation of these established tools.

\subsection{The exact finite spectrum}\label{app:gaussian-spectrum}

\begin{proof}[Proof of Lemma~\ref{lem:gaussian-spectrum}]
To construct the remainder, we diagonalize the finite operator $\mathsf Q_D$. Both $\mathsf Q_D^{-1}$ and $(\mathsf Q_D+m^2 I)^{-1}$ use the same eigenbasis, so their covariance difference is a scalar subtraction in each mode.
Let $A$ on $\{1,\ldots,L-1\}$ have adjacent entries $1/2$ and zero entries otherwise, with vectors extended by zero at $0,L$. The sine addition identity gives eigenvectors $\sin(\pi ki/L)$ with eigenvalues $\cos(\pi k/L)$. The eigenvalues are distinct, so symmetry gives orthogonality. The finite geometric sum of $\exp(2\pi\mathrm i ki/L)$, $0\leq i<L$, is zero for $1\leq k<L$. Taking real parts and using $\sin^2u=(1-\cos2u)/2$ gives
\[
 \sum_{i=1}^{L-1}\sin^2(\pi ki/L)=L/2.
\]
Thus the normalized one-dimensional vectors have factor $\sqrt{2/L}$. Their products have factor $2/L$ and form a basis. The identity $P_D=(A\otimes I+I\otimes A)/2$ gives the displayed eigenvalues, including the boundary rows. The inequalities $2u/\pi\leq\sin u\leq u$ on $[0,\pi/2]$ prove the two spectral bounds for every finite mode.

In this basis the covariance difference is \eqref{eq:gaussian-coefficients}, which is strictly positive. Take two independent arrays $(U_{k\ell})$ and $(V_{k\ell})$ of independent standard normals and form
\[
 \psi=\sum_{k,\ell}\frac{U_{k\ell}}{\sqrt{\lambda_{k\ell}+m^2}}e_{k\ell},
 \qquad
 \zeta=\sum_{k,\ell}V_{k\ell}\sqrt{a_{k\ell}}e_{k\ell}.
\]
Both sums are finite. The variances add to $\lambda_{k\ell}^{-1}$ mode by mode, proving the exact massless law. The argument used only $N\geq1$ and $m>0$.
\end{proof}

\subsection{Point and increment estimates}\label{app:gaussian-increments}

\begin{proof}[Proof of Lemma~\ref{lem:gaussian-increments}]
Write $\zeta=\zeta_{N,K}$ for the remainder of Lemma~\ref{lem:gaussian-spectrum}.
The point-variance bound uses mode amplitudes, while the increment bound also exploits the change of each sine mode between two sites.

\emph{Elementary sums with real cutoffs.}
For every real $T\geq1$, the usual decreasing-function comparison gives
\begin{equation}\label{eq:gaussian-elementary-sums}
 \sum_{1\leq n\leq T}\frac1n\leq1+\log T,
 \qquad
 \sum_{n>T}\frac1{n^3}\leq T^{-2}.
\end{equation}
For the second assertion, let $q=\lfloor T\rfloor+1\geq2$. Then
\[
 \sum_{n>T}n^{-3}
 \leq q^{-3}+\int_q^\infty u^{-3}\,du
 =q^{-2}(q^{-1}+1/2)\leq q^{-2}\leq T^{-2}.
\]
Also, if $1\leq K\leq B$, then
\begin{equation}\label{eq:gaussian-middle-harmonic}
 \sum_{K<n\leq B}\frac1n
 \leq\log(B/K)+\log2.
\end{equation}
Indeed, $p=\lfloor K\rfloor\geq K/2$ and $r=\lfloor B\rfloor\geq p$, so the sum is at most $\int_p^r du/u=\log(r/p)\leq\log(2B/K)$. This also covers the empty sum $r=p$. Finally,
\begin{equation}\label{eq:gaussian-logtwo}
 \frac12\leq\log2\leq\frac34.
\end{equation}
The lower bound follows by integrating $1/u\geq1/2$ on $[1,2]$. For the upper bound, the convex function $1/u$ lies below its endpoint chord $(3-u)/2$ on that interval, whose integral is $3/4$.

\emph{A coefficient bound for every frequency.}
Put $b=6K$ and $n=\max\{k,\ell\}$. Since $N\geq1$,
\[
 4\leq L/N=4+2/N\leq6,
 \qquad mL\leq b.
\]
Lemma~\ref{lem:gaussian-spectrum} and monotonicity of $u^2/(n^2+u^2)$ give
\begin{equation}\label{eq:gaussian-shell-coefficient}
 a_{k\ell}
 \leq\frac{L^2(mL)^2}{n^2\{n^2+(mL)^2\}}
 \leq\frac{L^2b^2}{n^2(n^2+b^2)}.
\end{equation}
There are exactly $2n-1$ positive-index pairs with maximum $n$. All summands below are nonnegative, so extending the finite frequency sum to all positive indices preserves the upper bounds.

\emph{Point variances.}
Since $|e_{k\ell}(x)|\leq2/L$, the covariance expansion and the shell count give
\[
 \Var\zeta(x)
 \leq8\sum_{n\geq1}\frac{b^2}{n(n^2+b^2)}
 \leq8\left\{\sum_{n\leq b}\frac1n
                  +b^2\sum_{n>b}n^{-3}\right\}
 \leq8(2+\log b).
\]
Here $b\geq6$, so \eqref{eq:gaussian-elementary-sums} applies. Since $6K\leq(2K)^3$ for $K\geq1$, we have $\log b\leq3\ell_K$. Equation~\eqref{eq:gaussian-logtwo} gives $2\leq4\ell_K$. Thus the last bound is at most $56\ell_K\leq64\ell_K$, proving \eqref{eq:gaussian-point}.

\emph{Spatial differences.}
Taking a difference of the two sine products, and using $|\sin u-\sin v|\leq|u-v|$, gives
\[
 |e_{k\ell}(x)-e_{k\ell}(y)|
 \leq\frac{4\pi n\|x-y\|_\infty}{L^2},
 \qquad
 |e_{k\ell}(x)-e_{k\ell}(y)|\leq\frac4L.
\]
Because $L\geq4N$ and $\pi<4$, these bounds imply
\[
 |e_{k\ell}(x)-e_{k\ell}(y)|^2
 \leq\frac{16}{L^2}\min\{1,n^2\rho^2\}.
\]
Consequently \eqref{eq:gaussian-shell-coefficient} and the shell count yield
\begin{equation}\label{eq:gaussian-increment-sum}
 \Var(\zeta(x)-\zeta(y))
 \leq32\sum_{n\geq1}
 \frac{b^2\min\{1,n^2\rho^2\}}{n(n^2+b^2)}.
\end{equation}
Write $w=K\rho\in(0,1]$ and $B=1/\rho\geq K$. Split this last sum into the three ranges $n\leq K$, $K<n\leq B$, and $n>B$.
The first cutoff is the mass scale, up to the fixed factor $b/K=6$; $B$ marks where the sine-difference gain $n\rho$ saturates. The low and middle ranges use that gain, the middle also uses spectral decay, and the high range uses spectral decay alone.

In the first range, $b^2/(n^2+b^2)\leq1$, so the contribution is at most
\[
 32\rho^2\sum_{n\leq K}n\leq32K^2\rho^2=32w^2.
\]
We used $\sum_{n\leq K}n\leq K(K+1)/2\leq K^2$, valid for real $K\geq1$. In the middle range, use $n^2+b^2\geq n^2$, $b^2=36K^2$, and \eqref{eq:gaussian-middle-harmonic}. The contribution is at most
\[
 1152w^2\sum_{K<n\leq B}\frac1n
 \leq1152w^2\{\log(1/w)+\log2\}.
\]
In the final range, \eqref{eq:gaussian-elementary-sums} with $T=B\geq1$ bounds the contribution by
\[
 1152K^2\sum_{n>B}n^{-3}\leq1152K^2B^{-2}=1152w^2.
\]
The sum of the three estimates is at most
\[
 w^2\{32+1152\log2+1152+1152\log(1/w)\}
 \leq2048w^2\{1+\log(1/w)\},
\]
by \eqref{eq:gaussian-logtwo} and $1152\leq2048$. This is \eqref{eq:gaussian-increment}. Empty middle ranges, real cutoffs, and finite-spectrum endpoints have all been included. In particular the argument remains valid when $K=N$, when $N=1$, and when $w=1$.
\end{proof}

\subsection{The integer grid and Gaussian concentration}\label{app:gaussian-net}

The proof first builds the expected-maximum estimate by following each lattice site through a finite sequence of parent grids. This is a concrete chaining argument. It then applies the Gaussian concentration statement in \cite[Corollary 6.6]{Biskup} to the absolute maximum, which is Lipschitz in the sup norm. The integer meshes keep the increment bound valid at every step and cover the endpoint where no refinement is needed. Thus the grid details serve the stated uniform numerical estimate, while the expectation and deviation principles remain the classical ones cited in Section~\ref{sec:gaussian}.

\begin{proof}[Proof of Lemma~\ref{lem:gaussian-net}]
Let $Z$ satisfy the numerical variance and increment hypotheses, and
write $\ell=\ell_K=\log(2K)$. We use integer grids so that all parameter
endpoints, including $K=N$ and $N=1$, are covered.
Starting at a mesh comparable to $N/K$ puts each successive
refinement within the range of the increment bound. Halving the
integer mesh down to one recovers every lattice site.

\emph{Integer grids and their parent maps.}
Set
\[
 s=\lfloor\log_2(N/K)\rfloor\geq0,
 \qquad h_i=2^{s-i},\quad0\leq i\leq s.
\]
The meshes are positive integers, and
\begin{equation}\label{eq:gaussian-initial-mesh}
 \frac{N}{2K}<h_0\leq\frac NK,\qquad h_s=1.
\end{equation}
For $r=1,2$, define
\begin{equation}\label{eq:gaussian-projections}
 \pi_i(x)_r=-N+h_i\left\lfloor\frac{x_r+N}{h_i}\right\rfloor,
 \qquad A_i=\pi_i(B_N).
\end{equation}
The grids are nested, $A_s=B_N$, and
\begin{equation}\label{eq:gaussian-grid-count}
 |A_i|=\left(\left\lfloor\frac{2N}{h_i}\right\rfloor+1\right)^2
 \leq25K^2 4^i.
\end{equation}
Indeed $2N/h_i<4K2^i$, while $K2^i\geq1$. For $i\geq1$ the parent
$p_i(z)=\pi_{i-1}(z)$ is unique for each $z\in A_i$, and
$\|z-p_i(z)\|_\infty\leq (N/K)2^{-i}$.
Since adjacent meshes differ by a factor of two and share the
anchor $-N$, projecting a point to level $i$ does not change its
level-$(i-1)$ parent. This compatibility makes the telescoping
identity below exact.
The increment condition and monotonicity of $u^2\log(e/u)$ for
$0<u\leq1$ imply
\begin{equation}\label{eq:gaussian-level-variance}
 \Var(Z_z-Z_{p_i(z)})\leq2048(i+1)4^{-i}.
\end{equation}
The number of distinct child-parent pairs is at most $|A_i|$.

\emph{Expected maximum.}
For centered Gaussian variables $Y_1,\ldots,Y_r$, each having variance
at most $\sigma^2$, their exponential moments give
\[
 \EE\max_{1\leq j\leq r}|Y_j|
 \leq\sigma\sqrt{2\log(2r)}.
\]
Since $\log(50K^2 4^i)\leq8(i+1)\ell$ for $K\geq1$,
\eqref{eq:gaussian-grid-count} and the point variance bound give
\[
 \EE\max_{z\in A_0}|Z_z|\leq32\ell.
\]
Likewise~\eqref{eq:gaussian-level-variance} gives for every $i\geq1$
\[
 \EE\max_{z\in A_i}|Z_z-Z_{p_i(z)}|
 \leq4\sqrt{2048}\,(i+1)2^{-i}\sqrt\ell.
\]
The exact telescoping identity
$Z_x=Z_{\pi_0(x)}+\sum_{i=1}^s
(Z_{\pi_i(x)}-Z_{\pi_{i-1}(x)})$
therefore yields, with $M=\max_{x\in B_N}|Z_x|$,
\begin{equation}\label{eq:gaussian-expected-maximum}
 \EE M\leq32\ell+12\sqrt{2048}\,\sqrt\ell
 =32\ell+384\sqrt{2\ell}.
\end{equation}
This calculation includes the case $s=0$, where the sum is empty.

\emph{Gaussian concentration and the numerical slack.}
The map $z\mapsto\max_{x\in B_N}|z_x|$ is $1$-Lipschitz for
the sup norm, and $\max_x\Var Z_x\leq64\ell$.
The Gaussian concentration inequality
\cite[Corollary 6.6]{Biskup}, applied to this map, therefore gives for $u\geq0$
\begin{equation}\label{eq:gaussian-borell}
 \PP(M>\EE M+u)\leq\exp\{-u^2/(128\ell)\}.
\end{equation}
The elementary bounds $\ell\geq\log2\geq2/3$,
$\sqrt\ell\geq4/5$, and $\sqrt2\leq3/2$ show from
\eqref{eq:gaussian-expected-maximum} that
\[
 800\ell-\EE M
 \geq\sqrt\ell\,(768\sqrt\ell-384\sqrt2)
 \geq16\sqrt\ell.
\]
Taking $u=(16+800t)\sqrt\ell$ in~\eqref{eq:gaussian-borell},
we obtain, for every $t\geq0$,
\[
 \PP\bigl(M>800(\ell+t\sqrt\ell)\bigr)
 \leq\exp\{-(16+800t)^2/128\}
 \leq e^{-2}e^{-8t^2}
 <\frac{25}{126}e^{-8t^2}.
\]
This proves~\eqref{eq:gaussian-net-tail}. No independence of Gaussian
values or increments was used.
\end{proof}

\subsection*{Funding.}
X. F. is supported by the National Science and Technology Council (NSTC) of
Taiwan, grant no.~114-2115-M-A49-003-MY3.

\subsection*{AI Declaration}

\enlargethispage{\baselineskip}
The proof in this work was developed during a long chain of conversations with OpenAI ChatGPT (GPT-6) which also helped with the preparation
of the manuscript and numerical verification.  All results were written and checked by the authors who take full responsibility.

\end{document}